\documentclass[11pt]{article}
 \usepackage{color}
\usepackage{amsmath}
\usepackage{amssymb}
\usepackage{amsthm}
\usepackage{epstopdf}
\usepackage{graphicx}
\usepackage{mathtools}
\usepackage[usenames,dvipsnames, table]{xcolor}
\usepackage[hang]{caption}
\usepackage{hyperref}
\usepackage{lscape} 

\usepackage{oldgerm}  

\usepackage{pgfplotstable}
\usepackage{longtable}

\def\be{\begin{equation}}\def\ee{\end{equation}}
\def\bee{\begin{enumerate}}\def\eee{\end{enumerate}}
\def\bei{\begin{itemize}}\def\eei{\end{itemize}}
\newcommand{\icon}[1]{\includegraphics[height=40pt]{#1}}
\newcommand{\SU}[1]{\mathrm{SU(#1)}}

\newcommand{\nco}{\newcommand}

\nco{\one}{\ensuremath{\,\,\mathrm{l}\!\!\!1}} 
\nco{\NN}{\mathbb{N}}
\nco{\ZZ}{\mathbb{Z}}
\nco{\QQ}{\mathbb{Q}}
\nco{\RR}{\mathbb{R}}
\nco{\CC}{\mathbb{C}}
\nco{\HH}{\mathbb{H}}
\nco{\OO}{\mathbb{O}}

\nco{\red}{\color{red}}
\nco{\blue}{\color{blue}}
\nco{\cyan}{\color{cyan}}
\nco{\brown}{\color{Magenta}}

\nco{\magenta}{\normalcolor}

\nco{\violet}{\color{violet}}
\nco{\redend}{\normalcolor}
\nco{\magentaend}{\normalcolor}

\def\ie{{\it i.e. }}
\def\ommit#1{{}}
\def\({\left(}
\def\){\right)}

\def\ie{{\it i.e.,\/}\ }
\def\ie{{\rm i.e.,\/}\ }
\def\etc{{\rm etc.\/}\ }

\def\be{\begin{equation}}\def\ee{\end{equation}}
\def\bea{\begin{eqnarray}}\def\eea{\end{eqnarray}}

\nco{\rnc}{\renewcommand}
\rnc{\title}[1]{{\Large\bf\mbox{}\\\medskip#1\bigskip\medskip\\}}
\rnc{\author}[1]{{\large #1\smallskip\\}}
\nco{\address}[1]{{\em #1\medskip\\}}

\begin{document}

\begin{titlepage}
\begin{center}
\title{Theta functions for lattices of $\SU{3}$ hyper-roots}
\medskip
\author{Robert Coquereaux} 
\address{Aix Marseille Univ, Universit\'e de Toulon, CNRS, CPT, Marseille, France\\
Centre de Physique Th\'eorique}

\bigskip\medskip

\today

\begin{abstract}
\noindent {We recall the definition of the hyper-roots that can be associated to modules-categories over fusion categories defined by the choice of a simple Lie group $G$ together with a positive integer~$k$. This definition was proposed in 2000,  using another language,  by Adrian Ocneanu. If $G=\SU{2}$, the obtained hyper-roots coincide with the usual roots for ADE Dynkin diagrams.  We consider the
associated lattices when $G=\SU{3}$ and determine their theta functions in a number of cases; these functions can be expressed as modular forms twisted by appropriate Dirichlet characters.
} 
\end{abstract}
\end{center}

\vspace*{70mm}
\end{titlepage}

\section{Introduction}

The ADE correspondence between indecomposable module-categories of type $\SU{2}$ and simply-laced Dynkin diagrams was  first obtained by theoretical physicists in the framework of conformal field theories  (classification of modular invariant partition functions for the WZW models of type $\SU{2}$, \cite{CIZ}, \cite{YellowBook}). Its relation with subfactors was studied in \cite{Ocneanu:paths} and it was set in a categorical framework by \cite{KirilovOstrik,Ostrik}. 
In plain terms, the diagrams encoding the action of the fundamental  representation of $\SU{2}$ at level $k$ (which is classically $2$-dimensional) on the simple objects of the various module-categories existing at that level,  are the Dynkin diagrams describing the simply-laced simple Lie groups with Coxeter number $k+2$. 

At a deeper level, there is a correspondence between fusion coefficients of the $\SU{2}$ module-category described by a Dynkin diagram ${\mathcal E}$
and the inner products between  
all the roots of the simply-laced Lie group associated with the same Dynkin diagram.
In the non-ADE cases one can also define the action of an appropriate ring on modules associated with the chosen Dynkin diagrams and still obtain a correspondence between structure coefficients describing this action and inner products between roots ---one has only to introduce scaling coefficients in appropriate places.

The correspondence relating fusion coefficients for module-categories of type $\SU{2}$ and inner products between weights and/or roots of root systems was clearly stated (but not much discussed) in \cite{Ocneanu:Bariloche}; in a different context, some of its aspects were already present in the article \cite{Dorey:CoxeterElement}. The correspondence  was used and described in some detail in one section of \cite{RC:periodicquivers}. As observed in \cite{Ocneanu:Bariloche}, one can start from the $\SU{2}$ fusion categories and their modules\footnote{They are called ``quantum subgroups'' of $\SU{2}$ in the latter reference.} to recover or define the usual root systems, and associate with each of them a periodic quiver describing, in particular, the inner products between all the roots. It is also observed in the same reference that the construction can be generalized: replacing $\SU{2}$ by an arbitrary simple Lie group $G$ leads, for every choice of a module-category ${\mathcal E}$ of type\footnote{Here and below, this means that the underlying monoidal category is ${\mathcal A}_k(G)$, whose definition is briefly recalled at the beginning of section~\ref{fusioncatbasics}.} $G$, to a system of ``higher roots'' that we call ``hyper-roots of type $G$''. Usual root systems are therefore hyper-root systems of type $\SU{2}$.

A usual root system gives rise, in particular, to an Euclidean lattice. The same is true for hyper-root systems. Given a lattice, one may consider its theta series whose $n^{th}$ coefficient gives the number of vectors of given norm, or, equivalently, the number of representations of the integer $n$ by the associated quadratic form. Theta functions of root lattices are well known and are usually expressed in terms of appropriate modular forms (see for instance the book \cite{ConwaySloane}).
Our purpose, in the present paper, is not to provide a detailed account of the properties of hyper-root systems --this should be done elsewhere
-- but only to discuss some general features of their theta functions and describe those functions associated with the systems defined by fusion categories of type ${\mathcal A}_k(\SU{3})$, or their modules, for small values of the (conformal) level $k$.
The structure of the hyper-root lattice of type $\SU{3}$ obtained when $k=1$  was announced in \cite{Ocneanu:Bariloche}:  it was there recognized as a scaled version of $D_6^+$, the so-called ``shifted $D_6$ lattice''. We shall recover and comment this result below\footnote{Section \ref{secL1}, about $D_6^+$, discusses some properties of the latter; it may have an independent interest.} but we shall also obtain closed formulae in terms of modular forms, as well as the corresponding series, for lattices of type $\SU{3}$ associated with higher valued of $k$. 

As already stated, the theory of hyper-roots that we use here is due to A. Ocneanu. 
Since it is poorly documented, we had to incorporate some general discussion based on a material that is, in essence, published in \cite{Ocneanu:Bariloche}, or available on line \cite{Ocneanu:MSRI}.
A general account of the theory of hyper-roots and of other higher analogues of  Lie groups concepts, as well as their interpretation in terms of usual representation theory, has been long awaited for, and should appear one day \cite{Ocneanu:WIP}. 
Let us stress again the fact that this is {\sl not} the purpose of the present article.
The family of scalar products between SU(3) systems of hyper-roots, called ``higher roots'' in \cite{Ocneanu:MSRI}, for several choices of ${\mathcal E}$, and using a different language, was obtained long ago \cite{Ocneanu:posters} and displayed in several places using beautiful posters. Our purpose, here, which is therefore the original contribution of this article, is to determine, in various cases, one or several convenient Gram matrices for the associated lattices, and to discuss some properties of their theta functions.  These results were hitherto apparently unknown,  this is why we decided to make them available.
Most of them were obtained in March 2009,  while the author
was a guest of the Mathematical Department at the University of Luxembourg,  whose hospitality is acknowledged.
We hope that this presentation will trigger new ideas and insights.

\section{From fusion categories and module categories to lattices of hyper-roots}

\subsection{On extended fusion matrices and their periods}
\label{fusioncatbasics}
From now on $k$ denotes a positive integer called the ``level'' (or conformal level), $G$ a simple\footnote{or semi-simple, but $k$ is then a multiplet of positive integers.}, simply connected, compact Lie group, and $Lie(G)$ its complex Lie algebra.
We call ${\mathcal A}_k(G)$ the category of integrable modules of the affine Kac-Moody algebra associated with $Lie(G)$ at level $k$, see e.g. \cite{Kac:book}. 
It is equivalent (an equivalence\footnote{The proof of equivalence given in the first two references assumed a negative level. 
The fact that it holds in all cases has been part of the folklore for a long time because it could be verified on a case by case basis. Its general validity 
is now considered as a consequence of the Huang's proof of Verlinde conjecture \cite{Huang}.} of modular tensor categories), 
 \cite{Finkelberg}, \cite{Huang}, \cite{KazhdanLusztig}, to a category constructed in terms of  representations of the quantum group $G_q$ at  root of unity ${\mathfrak q} =  exp(\frac{i \pi}{g^\vee+k})$, where $g^\vee$ is the dual Coxeter number of $G$  (take the quotient of the category of tilting modules by the additive subcategory generated by indecomposable modules of zero quantum dimension).
These categories, called fusion categories of type $(G,k)$, play a key role in the Wess - Zumino - Witten models of conformal field theory.
${\mathcal A}_k(G)$ being monoidal, with a finite number of simple objects denoted $m,n,p\ldots$, we consider the corresponding Grothendieck ring and its structure coefficients, the so-called fusion coefficients $N_{mnp}$, where  $m \times n = \sum_p N_{mnp} \, p$. They are encoded by fusion matrices $N_m$ with matrix elements $(N_m)_{np} = N_{mnp}$. The fusion category being given, one may consider module-categories\footnote{This amounts to say \cite{Ostrik} that we are given a monoidal functor from ${\mathcal A}_k(G)$ to the category of endofunctors of an abelian category ${\mathcal E}_k(G)$.} ${\mathcal E}_k(G)$ associated with it. 
In the following we assume that the chosen module-categories are indecomposable.
Of course, one can take for example ${\mathcal E}_k(G) = {\mathcal A}_k(G)$. The fusion coefficients $F_{nab}$ characterize the module structure: $n \times a = \sum_b F_{nab} \, b$, where $a,b,\ldots$ denote the simple objects of ${\mathcal E}_k(G)$. They are encoded either by square matrices $F_n$, with matrix elements $(F_n)_{ab}= F_{nab}$, still called fusion matrices, or by the rectangular matrices\footnote{The $F_n$ are sometimes called ``annular matrices'' when ${\mathcal A}_k(G)$ and ${\mathcal E}_k(G)$ are distinct (if they are the same, then $F_n=N_n$), and the $\tau_a$ are sometimes called ``essential matrices''.}  $\tau_a = (\tau_a)_{nb}$, with $(\tau_a)_{nb} =  (F_n)_{ab}$. The simple objects of ${\mathcal A}_k(G)$, or irreps, are labelled by the vertices of the Weyl alcove of $G$ at level $k$. With $G = \SU{2}$, this alcove is the Dynkin diagram $A_{k+1}$, vertices are labelled $(n)$, with non-negative integers $n\in\{0,1,\ldots,k\}$, and the fusion matrices $N_{(n)}$ or $F_{(n)}$ obey the Chebyshev recursion relation 
\be
F_{(n)}=F_{(n-1)} F_{(1)} - F_{(n-2)}
\label{SU2recursion}
\ee 
where $F_{(0)}$ is the identity matrix (the weight with Dynkin component $(0)$ is the highest weight of the the trivial representation) and $F_{(1)}$ refers to the generator (the fundamental irrep of classical dimension $2$).
 In the general case these matrices still obey recursion relations that depend on the choice of the underlying Lie group $G$. 
With $G = \SU{3}$, the simple objects (irreps) are labelled by pairs $(p,q)$ of non-negative integers with $p+q \leq k$, and the recursion relations read
\begin{eqnarray}
F_{(p,q)} &=& F_{(1,0)} \, F_{(p-1,q)} - F_{(p-1,q-1)} -F_{(p-2,q+1)} \qquad \qquad \textrm{if} \; q \not= 0
 \nonumber  \label{recursion}\\
F_{(p,0)} &=& F_{(1,0)} \, F_{(p-1,0)} - F_{(p-2,1)} 
\label{recF} \\
F_{(0,q)} &=& (F_{(q,0)})^T 
\nonumber
\end{eqnarray}
$F_{(0,0)}$ is the identity matrix and $F_{(1,0)}$ and $F_{(0,1)}$ are the two generators. Also, $F_{(q,p)}=F_{(p,q)}^T$.
Expressions of the fundamental fusion matrices $F_{(1,0)}$, for all the modules ${\mathcal E}_k(G)$ considered in this paper are recalled in the appendix (sec.~\ref{AppendixFusionMatrices}) using the weight ordering $(p_1,q_1) < (p_2,q_2)$ if 
$p_1+q_1 < p_2+q_2$ or if $q_1 < q_2$ and $p_1+q_1 = p_2+q_2$.

In some applications one sets to zero the fusion matrices whose Dynkin labels do not belong to the chosen Weyl alcove. This is \underline{{\sl not}} what we do here. 
At the contrary, the idea is 
to use the same recursion relations to extend the definition of the matrices $F_n$ at level $k$ from the Weyl alcove to
the fundamental Weyl chamber of $G$ (cone of dominant weights) and to use signed reflections with respect to the hyperplanes of the affine Weyl lattice in order to extend their definition to arbitrary arguments $n \in \Lambda$, 
 the weight lattice of $G$.
By so doing one obtains an infinite family of matrices $F_n$ that we shall still (abusively) call  ``fusion matrices", and for which we keep the same notations, although their elements can be of both signs.
It is also useful to shift  (translation by the Weyl vector)  the labelling index of the these matrices to the origin of the weight lattice;
in other words, for $n\in \Lambda$, using multi-indices, we set ${\{n\}} = {(n-1)}$, where the use of parenthesis refers to the usual Dynkin labels (we hope that this brace notation will not confuse the reader --- see examples below).
The following results about SU(2) and SU(3) are known and belong to the folklore.

If $G = \SU{2}$ 
 the terms $F_{\{n\}}$ with $1 \leq n\leq k+1$ are the usual fusion matrices at level $k$,  they have non-negative integer matrix elements,  $F_{\{1\}} = F_{(0)}$ is the identity and  $F_{\{0\}} = F_{(-1)}$ is the zero matrix; more generally
 the terms $F_{\{n\}}$ with $n=0 \; \text{mod} \; N$, where $N=k+2$, vanish. Matrices $F_{\{N+m\}}=-F_{\{N-m\}}$ have non-positive integer matrix elements for $1 \leq m \leq N-1$.
The $F$ sequence is periodic of period  $2N$ and the reflection symmetries (Weyl mirrors), with sign, are centered in position $\{n\} = 0  \; \text{mod} \; N$. 
Notice that $F_{\{2N-1\}} = - F_{\{1\}} = - \one$.

\begin{figure}[htbp]
\begin{minipage}{16pc}
\centering
\includegraphics[width=15pc]{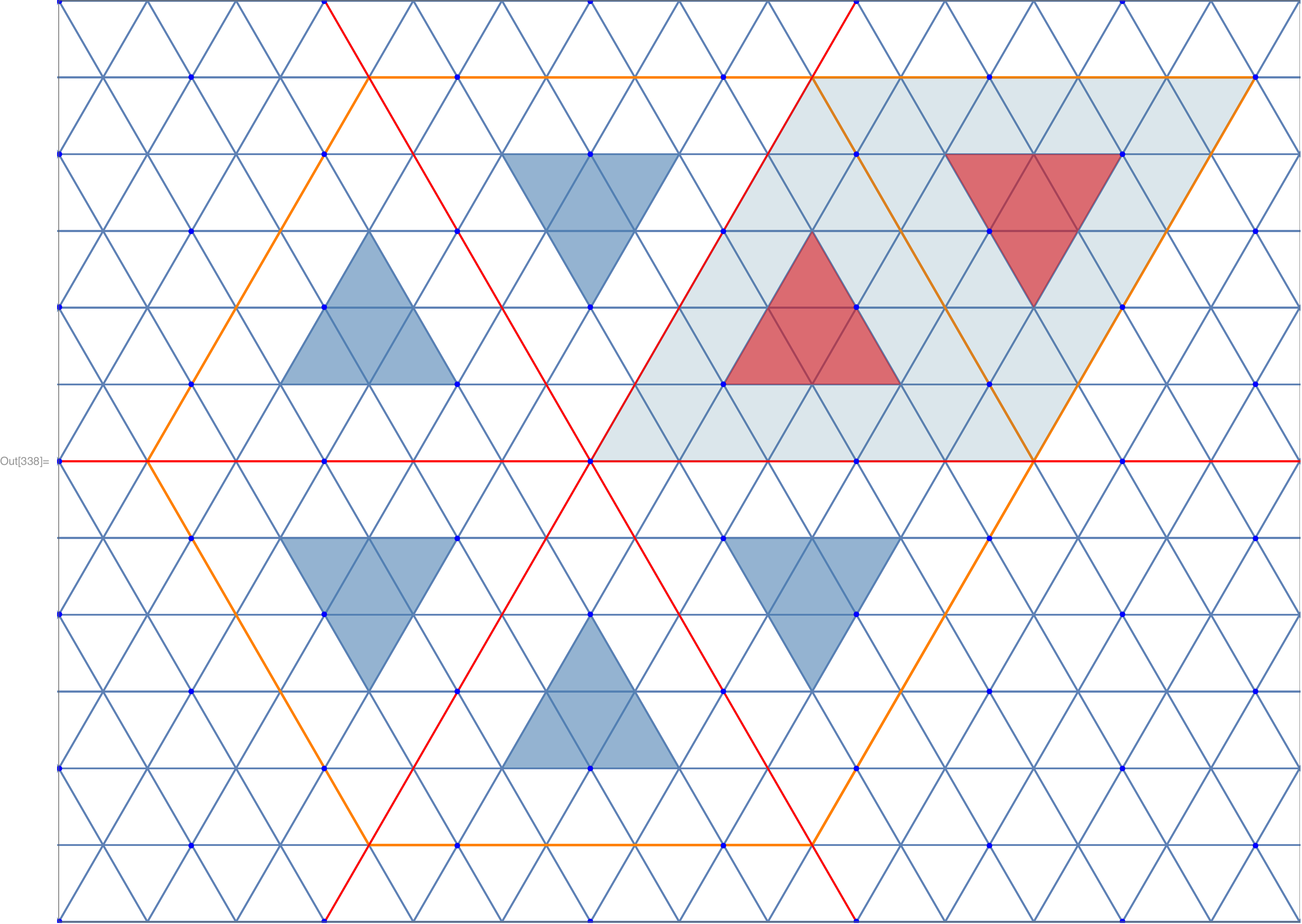}
\end{minipage}
\hspace{1.4pc}
 \begin{minipage}{16pc}
 \centering
 \includegraphics[width=22pc, height = 18pc, angle = 00]{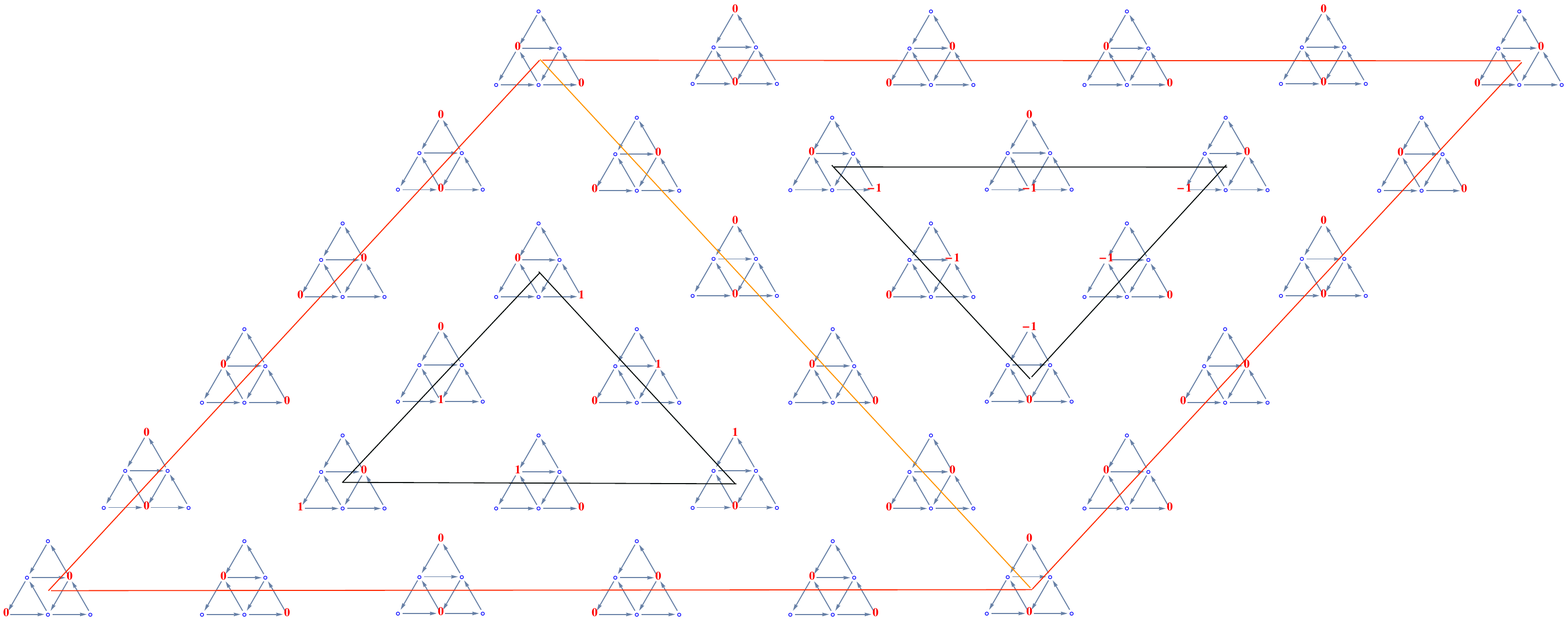}
\end{minipage}
\caption{\label{fondamentalrhombuslevel2} $\SU{3}$ at level $2$.\\
Left: The alcove (lower red triangle), its images under the Weyl group, and the periodicity rhombus.
Right: The function $\tau_a$ on ${\mathcal R}^\vee$  for $a$, the trivial irrep, with highest weight $(0,0)$;
Notice that reflection across the drawn diagonal induces a sign flip of the values of $F_n$ or of the function $\tau_a$.\\
Ideally, these figures, like those that follow, should be magnified on a terminal device.}
\end{figure}

If $G = \SU{3}$, we set $F_{\{p,q\}} = F_{(p-1,q-1)}$, so that  $F_{\{0,0\}}$ is the zero matrix and $F_{\{1,1\}}=F_{(0,0)}$ is the identity (the latter corresponding to the weight with components $(0,0)$ in the Dynkin basis, \ie to the highest weight of the trivial representation).  One has  $F_{\{p,q\}} = 0$ whenever $p=0 \; \text{mod} \; N$, $q=0 \; \text{mod} \; N$ and $p+q=0 \; \text{mod} \; N$.
One also gets immediately the following equalities:  $F_{\{p+N,q\}}= (P.F)_{\{p,q\}}$,  $F_{\{p,q+N\}}= (P^2.F)_{\{p,q\}}$ where $P=F_{\{N-2,1\}}$ is a generator of $\ZZ_3$ (with $P^3=1$) acting by rotation on the fusion graph of ${\mathcal A}_k(SU(3))$ and $F_{\{p+3N,q\}}= F_{\{p,q+3N\}} = F_{\{p+N,q+N\}}=F_{\{p,q\}}$.
The sequence $F_{\{p,q\}}$  is periodic of period  $3 N$ in each of the variables $p$ and $q$ but it is completely characterized by the values that it takes in a rhombus $(N,N)$ with $N^2$ vertices; for this reason, this rhombus will be called periodicity cell, or periodicity rhombus. We have reflection symmetries (with sign) with respect to the lines $\{p\}=0 \; \text{mod} \; N$, $\{q\}=0 \; \text{mod} \; N$ and $\{p+q\}=0 \; \text{mod} \; N$.
The $F$ matrices  labelled by vertices belonging to the Weyl alcove (which can be strictly included in the first half of a periodicity rhombus) have non-negative integer matrix elements; those with indices belonging to the other half of the inside of the rhombus have non-positive entries, those with vertices belonging to the walls of the Weyl chamber or to the second diagonal  of the rhombus vanish, and the whole structure is periodic. The Weyl group action 
on the alcove\footnote{This is the shifted Weyl action: $w \cdot n= w (n+\rho) - \rho$ where $\rho$ is the Weyl vector.} and the affine SU(3) lattice at level $k=2$ are displayed in figure~\ref{fondamentalrhombuslevel2}, left.

More generally, for a simple Lie group $G$ taken\footnote{meaning that we consider the fusion category ${\mathcal A}_k(G)$ or one of its module-categories} at level $k$, the obtained periodicity cell, 
once $F$ matrices have been appropriately extended to the whole weight lattice, 
 is a parallelotope $D$ with $\eta N^{r_G}$ vertices where
$N=k+g$, $g$ being the Coxeter number of $G$ and $r_G$ its rank, and where the value of $\eta$, a small integer, depends on the symmetries of the Dynkin diagram of $G$. We saw that $\eta=2$ for SU(2) but  $\eta=1$ for SU(3). 
The integer $N$ is sometimes called altitude of the module  \cite{DiFrancescoZuber}, or generalized Coxeter number (it coincides with the Coxeter number of the group SU(N)).

\subsection{On periodic essential matrices and the ribbon of hyper-roots}

\subsubsection{The ribbon ${\mathcal R}^\vee$}
\label{sec:ribbon}
Given a module-category ${\mathcal E}_k(G)$ over ${\mathcal A}_k(G)$ (the former, that we shall just call ${\mathcal E}$, if no confusion arises, using a notation that will also denote the set of isomorphisms classes of its simple objects,
can be chosen equal to the latter, but it is good to keep the distinction in mind), we defined, for each simple object $a$ of ${\mathcal E}_k(G)$, an essential matrix  $\tau_a$, with elements $(\tau_a)_{n,b}=F_{n,a,b}$, where $n,b$ refer respectively to the simple objects of ${\mathcal A}_k(G)$ and of ${\mathcal E}$. Since we have extended the definition of the fusion matrices $F_n$ to allow arguments $n$ belonging to the weight lattice of $G$ by using recursion relations, symmetries, and periodicity, we can do the same for the $\tau$'s, keeping the same notation: the indices $a,b$ of $(\tau_a)_{nb}$ still refer to simple objects of ${\mathcal E}$ but the index $n$ labels weights of $G$. 
The infinite matrices $\tau_a$ can be thought as a rectangular, with columns indexed by the elements of ${\mathcal E}$ (a finite number) and lines indexed by the weights of $\Lambda$, the weight lattice of $G$.

For every choice of $a \in {\mathcal E}$, $\tau_a$ is therefore a periodic, integer-valued function, on $\Lambda \times {\mathcal E} $.
Actually, in many cases the definition domain of $\tau_a$ can be further restricted. 
Indeed, there are many modules ${\mathcal E}_k(G)$ that have a non trivial grading with respect to the center ${\mathcal Z}$ of $G$;  in those cases, not only the weights of $G$, its irreducible representations, the simple objects of ${\mathcal A}_k(G)$, but also the simple objects of  ${\mathcal E}_k(G)$, have a well defined grading (denoted $\partial$) with respect to  ${\mathcal Z}$, and the module structure is compatible with this grading: matrix elements of $\tau_a$ in position $(n,b)$ will  automatically vanish if $\partial n + \partial a \neq \partial b$.  Unless stated otherwise we shall assume in the rest of the paper that we are in this situation.
The function $\tau_a$ is  periodic, integer-valued on $\Lambda \times_{\mathcal Z} {\mathcal E}$, and it is specified (see figure~\ref{fondamentalrhombuslevel2}, right) by the values that it takes on the finite set ${\mathcal R}^\vee = D \times_{\mathcal Z}  {\mathcal E}$ where $D$ is the period parallelotope.\
 The set ${\mathcal R}^\vee$,  a finite rectangular table made periodic, may be thought as a closed ribbon\footnote{The terminology ``ribbon'' comes from A. Ocneanu.}. 
 For most choices of ${\mathcal E}$, in particular if one takes ${\mathcal E} =  {\mathcal A}_k(G)$, 
  the group ${\mathcal Z}$ acts non trivially and ${\mathcal R}^\vee$ has $r_{\mathcal E}  \; \vert{D}\vert / \vert {\mathcal Z} \vert $ elements, where the rank $r_{\mathcal E}$ is the number of simple objects of ${\mathcal E}$,
  and $\vert{D}\vert = \eta N^{r_G}$.
The elements of ${\mathcal R}^\vee$ will be called (restricted\footnote{for reasons explained in section~\ref{sec:unrestricted}.}) hyper-roots of type $G$ defined by the module ${\mathcal E}_k(G)$.

\begin{figure}[ht]
\centering
\includegraphics[width=20pc]{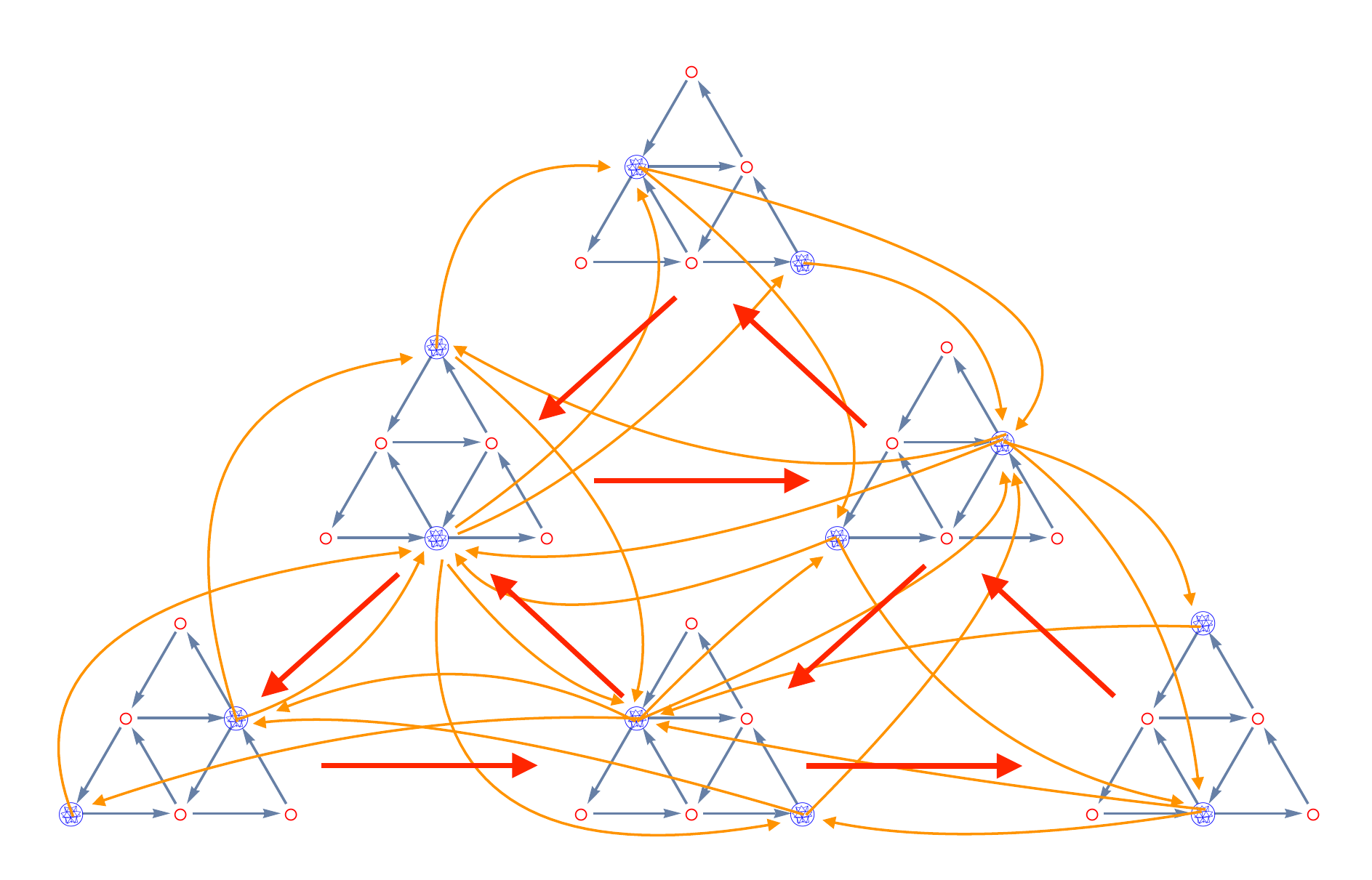}
\caption{\label{A2SU3QuiverStart} 
Edges (orange arrows) of the quiver ${\mathcal R}^\vee$ for ${\mathcal A}_2(\SU{3})$, $\pi$ being one of the two fundamental representations: 
start from an allowed vertex (a blue star) of a subgraph ${\mathcal E}$ (a blue triangle),  go  to the vertex, or vertices, located in the same position in the neighboring subgraphs by using the reversed edges of $\Lambda$, \ie
following the red arrows backwards, then follow the blue arrows on the latter subgraphs.}
\end{figure}

The choice of a fundamental irrep $\pi$ of $G$, with the constraint that it should exist at level $k$ (so that $\pi$ defines a particular non-trivial simple object of ${\mathcal A}_k(G)$)
allows one to consider ${\mathcal E}$ as a graph or, rather, to associate with ${\mathcal E}$, a graph denoted by the same symbol, once $\pi$ is chosen once and for all :  it is the graph of multiplication by $\pi$, sometimes called fusion graph, representation graph, nimrep graph, or McKay graph associated with $\pi$. If $\pi$ is complex, like the fundamental representation(s) of $\SU{3}$,  edges of ${\mathcal E}$ are oriented. If $\pi$ is self-conjugate, like the fundamental representation of $\SU{2}$, or like the antisymmetric square of the vector representation of SU(4), edges carry both orientations and can be considered as non-oriented. In general ${\mathcal E}$ is actually a quiver since it is a directed graph where loops and multiple arrows between two vertices are allowed. For any choice of a fundamental irrep $\pi$ of $G$ existing at level $k$, the set ${\mathcal E}$, and therefore the ribbon ${\mathcal R}^\vee$ as well, become quivers (there is an edge from one vertex of ${\mathcal R}^\vee$ to another if there are edges between their respective two projections  in  ${\mathcal E}$ and in $\Lambda$).  The definition of ${\mathcal R}^\vee$ as a set of vertices does not depend on the choice of $\pi$.
Figure~\ref{A2SU3QuiverStart} displays the first few edges (orange arrows) between the vertices of the bottom left corner of  the ${\mathcal R}^\vee$ quiver for ${\mathcal E}={\mathcal A}_2(\SU{3})$; since this is essentially a cartesian product\footnote{actually we follow the reversed red arrows in figure~\ref{A2SU3QuiverStart}, which means that we use the opposite $\Lambda$, but this choice is purely conventional and plays no role in the sequel.} of two multigraphs (red and blue arrows in the same picture) we shall no longer displays the edges of ${\mathcal R}^\vee$ in subsequent illustrations. 
If $\SU{3}$ is replaced by $\SU{2}$, one obtains in the same way the quivers of roots for all simple Lie groups; several examples of this construction in the case of usual roots, for instance the quiver of roots of $E_6$, can be found in \cite{RC:periodicquivers}. 

\subsubsection{The case $G=\SU{2}$}
When $G$ is $\SU{2}$, the module-categories ${\mathcal E}_k(\SU{2})$ are classified by ADE Dynkin diagrams. If one chooses for instance $k=10$, there are three of them, described by the Dynkin Diagrams $A_{11}, D_7, E_6$, the first being the modular fusion category ${\mathcal A}_k(\SU{2})$ itself, and $N=k+2$ is the Coxeter number. 
At the level of sets it is clear that the ribbon ${\mathcal R}^\vee$ is in bijection with the root system defined by the chosen Dynkin diagram since the number of roots is indeed equal to $r_{\mathcal E} \times N$. Let us briefly mention why  ${\mathcal R}^\vee =  \ZZ_{2N} \times_{\ZZ_2}{\mathcal E}$ can  be identified with the periodic quiver of roots.
As graphs, the Dynkin diagram ${\mathcal E}$ encodes the multiplication of the simple objects of ${\mathcal E}_k(\SU{2})$  by the fundamental irrep of $\SU{2}$.
The Chebyshev recursion relation \ref{SU2recursion} now reads 
$ \; F_{\{2\}} \, F_{\{p\}} = F_{\{p-1\}} + F_{\{p+1\}} $, 
but $F_{\{2\}}$ is also the adjacency matrix of the graph ${\mathcal E}$ since $\{2\}=(1)$ denotes the fundamental irrep of $\SU{2}$.
 In other words, the functions $\tau_a$ are such that the sum of neighbors taken vertically (\ie along  $\ZZ_{2N}$) equals the sum of neighbors taken horizontally (\ie along ${\mathcal E}$)  on the graph ${\mathcal R}^\vee$.
This can be written as a harmonicity property: define the laplacian  $\Delta_{\mathcal E}$ on ${\mathcal E}$ as the sum of neighbors, and similarly for the laplacian $\Delta_\Lambda$ on $\ZZ_{2N}$; a function $f$ such that $\Delta_\Lambda f = \Delta_{\mathcal E} f$ is called harmonic. The functions $\tau_a$ are therefore $\ZZ$-valued and harmonic. As the vertices of Dynkin diagrams also label fundamental weights, one has a function $\tau_a$ on ${\mathcal R}^\vee$ for every fundamental weight, therefore weights define functions that are harmonic.
If one thinks of a root as a point (a Dirac measure) of the ribbon, one can show (see section 3 of \cite{RC:periodicquivers}, in relation with \cite{Ocneanu:Bariloche}) that $(\tau_a)_{n,b}$ is the inner product  between the fundamental weight $a$ and the root localized at the point $(n,b)$.
This discussion justifies the terminology  ``ribbon of hyper-roots'' since when $G=\SU{2}$, the period is $2N$, the period parallelotope is the interval $D=\ZZ_{2N}$,  and the ribbon ${\mathcal R}^\vee = \ZZ_{2N} \times_{\ZZ_2} {\mathcal E}$ can be identified with the periodic quiver of roots. 
From another point of view, roots are weights, therefore roots also define (particular) $\ZZ$-valued harmonic functions on ${\mathcal R}^\vee$, and since $\langle \alpha, \alpha \rangle = 2$,
 the point where a root $\alpha$ is localized on the ribbon, as a Dirac measure, is obtained from the collection of inner products $\langle \alpha, \beta \rangle$ between $\alpha$ and all the roots as the (unique) point where this value is equal to $2$. 
 
Let $\alpha = (m,a)$ and $\beta=(n,b)$ two vertices of ${\mathcal R}^\vee$ \ie two roots.
One finds:
\be
\langle \alpha, \beta \rangle = F_{(m-n), a,b} +  F_{(n-m), a, b}  = F_{(m-n), a,b}  -  F_{(m-n-2), a, b} 
\label{scalarproductformula1}
\ee
where the fusion coefficients $F_{(n) a b}$ have been extended by periodicity as explained in section~\ref{fusioncatbasics}.
In terms of fusion matrices $F_{\{n\}}$ with shifted labels and matrix elements $({F_{\{n\}}})_{(a,b)}= F_{\{n\}, a,b}  = F_{(n-1),a,b}$, this relation reads
\be
  < (m,a),(n,b)> = ( F_{\{m-n+1\}} - F_{\{m-n -1\}})_{(a,b)}
  \label{scalarproductformula2}
\ee

Although expressed in terms of dimensions of spaces of essential paths on graphs (a concept that we shall not use in the present paper),  
equation \ref{scalarproductformula1} was explicitly written in section 1.7 of  \cite{Ocneanu:Bariloche}. What was then proposed, in this last reference, is to use this expression as a starting point in order to {\sl define} the inner product between {\sl all}  the vertices of the ribbon ${\mathcal R}^\vee$, \ie all the roots, without relying on the existence of a special basis (the simple roots) in which the inner products would be given by elements of the usual Cartan matrix, and finally to consider higher generalizations where $\SU{2}$ is replaced by another simple or semi-simple Lie group $G$.

\subsubsection{SU(3) and the general case}
\label{sec:harmonicity}

Again, up to notations and a different terminology, the content of the present section is already present in or can be inferred from reference \cite{Ocneanu:Bariloche}, pp 9-10.
The vector space of complex valued functions on the set of hyper-roots, the ribbon, is ${\CC^{\vert {\mathcal R}^\vee \vert}}$.  
It admits a canonical basis whose elements are identified with characteristic functions $\delta_\alpha$ (Dirac measures located at the points $\alpha \in {\mathcal R}^\vee$). An Euclidean structure is defined on this space by declaring that these Dirac masses are orthonormal.  
The elements $f$ of the subspace of harmonic functions are such\footnote{This harmonicity property is illustrated for ${\mathcal A}_2(SU3)$ in figure~\ref{A2roohexagon5harmonicity}.} that $\Delta_\Lambda f = \Delta_{\mathcal E} f$, where $\Delta_\Lambda$ and $\Delta_{\mathcal E}$ respectively denote the laplacian on the weight lattice $\Lambda$ of $G$ and the laplacian on the graph ${\mathcal E}$; in order to define the latter, one should also, in principle, select a fundamental irrep $\pi$ of $G$ existing at the chosen level, but this choice will be irrelevant for the study of the case $G=\SU{3}$ that we shall investigate later in more details.

In the classical situation weights are harmonic and should have integral inner products with roots, therefore one defines hyper-weights as $\ZZ$-valued functions that are harmonic on the ribbon.
Hyper-roots can be thought either as points $\alpha$ of the ribbon (Dirac masses), or as particular hyper-weights.
A point  $\alpha$ of ${\mathcal R}^\vee$ specifies a harmonic function, also denoted $\alpha$, defined as the orthonormal projection of the Dirac measure $\delta_\alpha \in {\CC^{\vert {\mathcal R}^\vee \vert}}$ on the subspace 
of harmonic functions. One finds  \cite{Ocneanu:Bariloche}
 that its value on $\beta \in {\mathcal R}^\vee$, denoted $<\alpha, \beta >$, is explicitly given,  if $\alpha = (m,a)$ and $\beta=(n,b)$,  by
\be
<\alpha, \beta > =  \sum_{w \in {\mathcal W}} \epsilon(w) F_{m-n + w\rho - \rho, a,b} 
\label{scalarproductofhyperroots1}
\ee
The above expression generalizes the one previously obtained for $G=\SU{2}$, where ${\mathcal W}=S_2$, see equation~\ref{scalarproductformula1}.
It is a finite signed sum with $\vert{\mathcal W}\vert$ terms, where ${\mathcal W}$ is the Weyl group of $G$,  $\epsilon$ is the determinant of Weyl reflections,  and $\rho$ the Weyl vector.
Using shifted labels, as in equation~\ref{scalarproductformula2}, the above reads
\be
<\alpha, \beta > =  \sum_{w \in {\mathcal W}} \epsilon(w) F_{\{m-n + w\rho\}, a,b} 
\label{scalarproductofhyperroots2}
\ee
One can take the latter expression as a definition of the inner product and extend $\langle  \; , \;  \rangle$ by linearity to the linear span of hyper-roots.  
One checks that it defines a positive definite\footnote{
Using eq.~\ref{scalarproductofhyperroots1}
one could define a periodic inner product  on  $\Lambda \times_{\mathcal Z}{\mathcal E}$ that would not be positive definite because of the periodicity, 
but we consider directly its non-degenerate quotient, naturally defined on the ribbon $D \times_{\mathcal Z}{\mathcal E}$.}
 inner product and therefore an Euclidean structure on the space of hyper-roots. This euclidean space will be denoted~${\mathfrak C}$.
Notice that\footnote{The group $G$  may be simply-laced or not, but for the modules considered in this paper (choices of ${\mathcal E}$), all hyper-roots have only one possible length.}
 $<\alpha, \alpha > =  \vert{\mathcal W}\vert$ for all hyper-roots $\alpha$. \\
Terminological conventions: elements of the hyper-root lattice $L$, the $\ZZ$ span of hyper-roots, are ``hyper-root vectors'' and the
elements of the dual lattice $L^\star$, are ``hyper-weight vectors''.

If $G=\SU{3}$, the Weyl group is ${\mathcal S}_3$ and  the inner product of two hyper-roots is obtained from equation~\ref{scalarproductofhyperroots1}, or \ref{scalarproductofhyperroots2},
as the sum of six fusion coefficients. 
Writing $\alpha=(m,a) = ((m_1,m_2),a)$, $\beta=(n,b) = ((n_1,n_2),b)$, setting $\lambda_1=m_1-n_1$, $\lambda_2=m_2-n_2$, and using shifted labels, we obtain
\begin{eqnarray}
& & <\alpha, \beta > \quad =\\
\nonumber
& & \left(F_{\{\lambda_1+1,\lambda_2+1\}}  +  F_{\{\lambda_1-2,\lambda_2+1\}}+F_{\{\lambda_1+1,\lambda_2-2)\}}-F_{\{\lambda_1-1,\lambda_2-1\}}-F_{\{\lambda_1-1,\lambda_2+2\}}-F_{\{\lambda_1+2,\lambda_2-1\}} \right)_{(a,b)}
\label{scalarproductofsu3roots}
\end{eqnarray}

\subsubsection{From the set ${\mathcal R}^\vee$ to the set ${\mathcal R}$}
\label{sec:unrestricted}
Let $\alpha$ be a point of ${\mathcal R}^\vee$, this defines a hyper-root, also called $\alpha$, as a vector of the euclidean space~${\mathfrak C}$. Obviously, its negative $-\alpha$ is another vector of ${\mathfrak C}$. 
For usual roots, \ie SU(2) hyper-roots, the opposite of a root is a root, as it is well known. However, for SU(3) hyper-roots, we can see, using the definition of the period parallelotope $D$, that $-\alpha$ does not correspond to any vertex of ${\mathcal R}^\vee$. This feature is not convenient. For all purposes it is useful to  generalize the previous definitions, keeping ${\mathcal R}={\mathcal R}^\vee$ for SU(2) but  setting ${\mathcal R}= {\mathcal R}^\vee  \cup - {\mathcal R}^\vee$
for SU(3), then  $\vert {\mathcal R} \vert =2  \vert  {\mathcal R}^\vee  \vert$.
The opposite of a hyper-root (an element of ${\mathcal R}$) is then always a hyper-root. 

For modules ${\mathcal E}$ endowed with a non trivial grading by the center ${\mathcal Z}$ of $G$, we already know that:
\be
\vert {\mathcal R}^\vee \vert =  \frac{  r_{\mathcal E}  \;   \vert {D} \vert }{  \vert {\mathcal Z} \vert } =   \frac{ r_{\mathcal E} \, N^{r_G}}  {\vert {\mathcal Z} \vert} \, \eta
\ee

Since we only study $G=$SU(3) hyper-roots in this paper, we shall not give more details about  ${\mathcal R}$ versus ${\mathcal R}^\vee$ for a general $G$, nevertheless we notice, 
in view of the last comment of section \ref{fusioncatbasics} on the period parallelotope $D$, its size being $\eta N^{r_G}$ with $\eta =1$ for SU(2) and $\eta =2$ for SU(3), and from the above definition of ${\mathcal R}$,
that we obtain in both cases a number of hyper-roots given by the same general expression: 
 \be
 \vert {\mathcal R} \vert = 2 \frac{ r_{\mathcal E} \, N^{r_G}}  {\vert {\mathcal Z} \vert}
 \ee  
which, more specifically, reads $ \vert {\mathcal R} \vert = r_{\mathcal E} \, N$ for $G=$SU(2), as it should, and reads for $G=\SU{3}$,
\be  \vert {\mathcal R}\vert = 2 \vert {\mathcal R}^\vee \vert = \frac{2}{3}  r_{\mathcal E} \, N^2   \ee
since $r_G=2$ and $\vert {\mathcal Z} \vert = 3$; we have also $N=k+3$.
Moreover if one chooses ${\mathcal E} = {\mathcal A}_k$ one has $r_{\mathcal E} = (N-2)(N-1)/2$, then
\be  \vert {\mathcal R}\vert =  2\vert {\mathcal R}^\vee \vert = (N-2)(N-1) N^2 /3   \ee
From now on, we shall usually not mention ${\mathcal R}^\vee$, the set of restricted hyper-roots  that we had to introduce in the first place, since ${\mathcal R}$ will be used most of the time.

\subsection{Dimension of the space of hyper-roots}

The dimension $\mathfrak{r} = \text{dim} \; {\mathfrak C}$ of the space of hyper-roots associated with ${\mathcal E}_k(G)$, in those cases
 where the center ${\mathcal Z}$ acts non-trivially on the set of vertices of ${\mathcal E}_k(G)$, is\footnote{This general result was claimed in the last two slides of \cite{Ocneanu:MSRI} and it can be explicitly checked in all the cases that we consider below.}
\be
\label{dimCgeneral}
 \mathfrak{r} = \text{dim} \; {\mathfrak C} = \frac{ r_{\mathcal E} \, \vert {\mathcal W} \vert }{\vert {\mathcal Z} \vert} 
 \ee
 where ${\mathcal W}$ is the Weyl group associated with the  simple Lie group $G$.
The term ${\vert {\mathcal W} \vert }/{\vert {{\mathcal Z}} \vert }$ cancels out for $G=\SU{2}$ since $ {\mathcal W}$ and $ {\mathcal Z}$ are both isomorphic with $\ZZ_2$. One then recovers the rank $\mathfrak{r}  = r_{\mathcal E}$ given by the number of vertices of the chosen Dynkin diagram.

\smallskip
For $G=\SU{3}$, ${\mathcal W} ={\mathcal S}_3$, $\vert {\mathcal W} \vert=3!$, and for modules ${\mathcal E}$ with non trivial triality (a property that holds for all the examples that we shall considered below) we have
${\mathcal Z}=Z_3$ therefore \be \textrm{For SU(3),} \qquad \mathfrak{r}=  2\, r_{\mathcal E}\ee\\
Moreover, If one chooses ${\mathcal E} = {\mathcal A}_k$,  then $\mathfrak{r} = (N-2)(N-1)$. 

\smallskip
More generally for $G=\SU{n}$, ${\mathcal W} ={\mathcal S}_n$ and for modules ${\mathcal E}$ endowed with a non-trivial grading by the center
 ${\mathcal Z}=Z_n$,  one has  $\mathfrak{r} =  r_{\mathcal E} \, n!/n = (n-1)! \,  r_{\mathcal E}$.

\bigskip
Remark : The Weyl group has order $\vert {\mathcal W} \vert = r_G ! \,  \Pi(\theta_\alpha) \, \vert {\mathcal Z} \vert$ where   $\Pi(\theta_\alpha)$ is the product of the components of the highest root of $G$ in the basis of simple roots\footnote{If $G$ is not simply-laced, one should be careful not to use here the basis of simple coroots.}.  One can therefore rewrite equation \ref{dimCgeneral} as
\be 
\mathfrak{r} =  r_{\mathcal E} \,  r_G ! \;  \Pi(\theta_\alpha)
\ee

\subsection{Example:   ${\mathcal E} = {\mathcal A}_k(\SU{3})$}

\begin{figure}[ht]
\centering
\includegraphics[width=32pc]{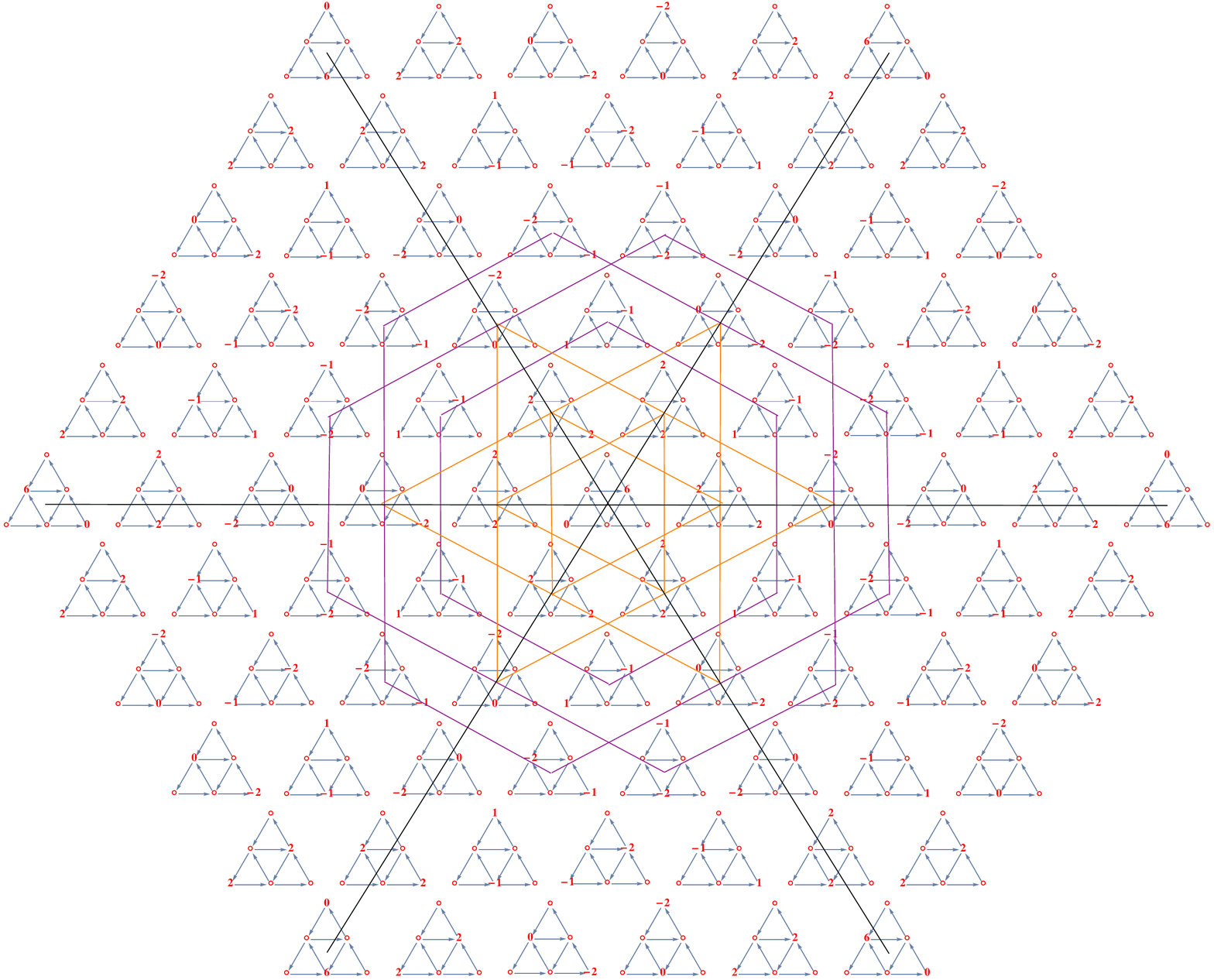}
\caption{\label{A2roohexagon5symmetries} $\SU{3}$ at level $2$: one of the hyper-root hexagons (symmetry properties).}
\end{figure}

At level $k$  we have $N=k+3$, and the values of $\vert {\mathcal R}^\vee \vert$, $\vert {\mathcal R}\vert$, $r_{\mathcal E}$ and $\mathfrak{r}$, in terms of $N$, were given before.
 Taking for instance $k=2$ there are $100$ hyper-roots,  $50$ strict hyper-roots, $r_{\mathcal E} =6$, $\mathfrak{r} =12$,
and we display in figure~\ref{A2roohexagon5symmetries}  one of the (hyper) root hexagons associated with some chosen hyper-root that sits where the integer $6$ (inner product with itself) appears\footnote{It may be useful to enlarge these pictures, using an online version of the present paper.}, in the central fusion diagram; the inner products of that hyper-root with all the others are given by the integers that appear in the figure.
The small circles with no integer marks are forbidden vertices (points where the constraint $\partial n + \partial a = \partial b$, see section~\ref{sec:ribbon},  is not obeyed).
Notice that $N=5$, so that drawing hexagons with edges of length $5$ would be sufficient to display all the inner products and their Weyl symmetries, but we drew an hexagon with edges of  length $N+1=6$ in order to illustrate the periodicity properties of the chosen hyper-root.

\begin{figure}[htbp]
\centering
\includegraphics[width=32pc]{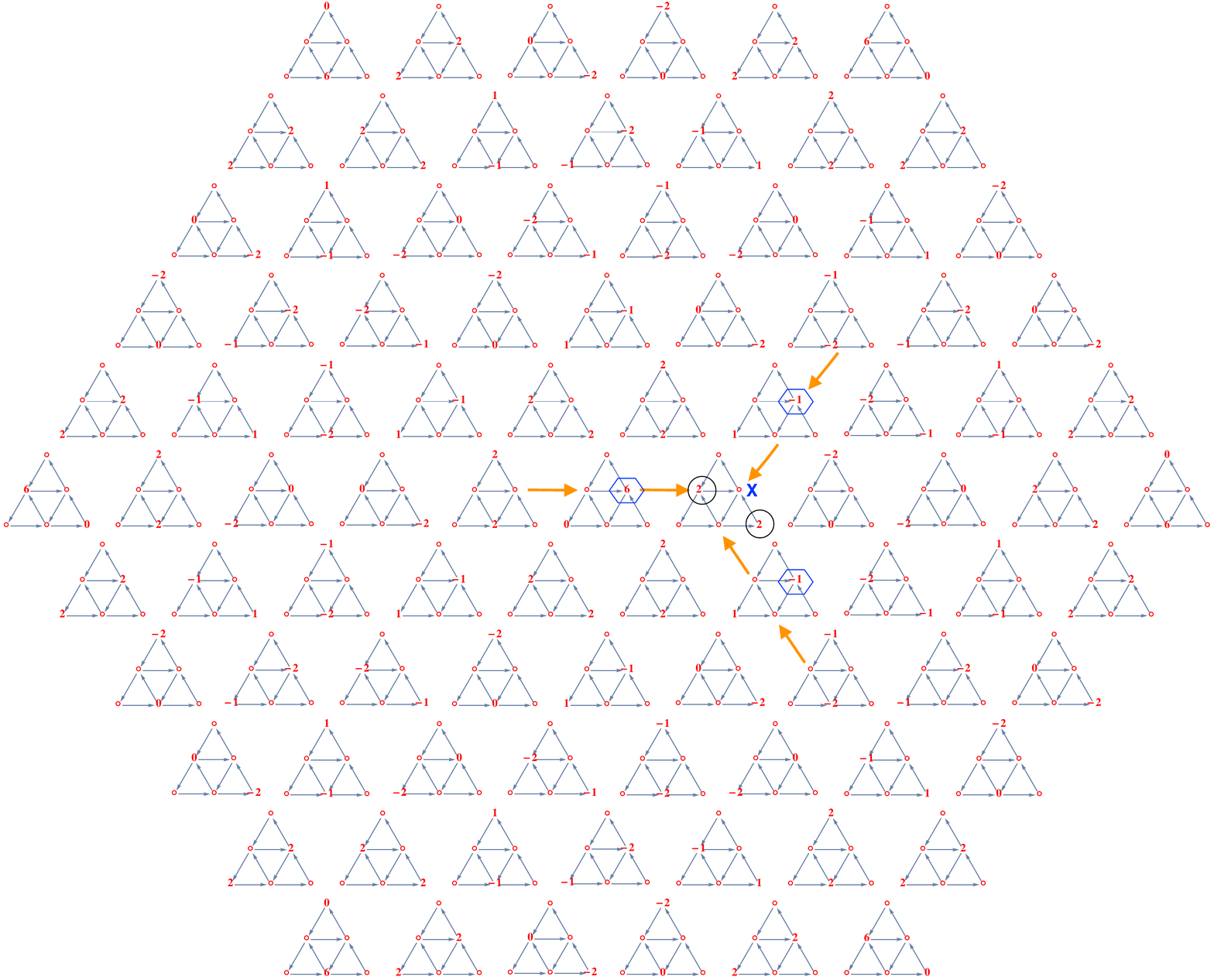}
\caption{\label{A2roohexagon5harmonicity}
Harmonicity property for one of the hyper-root hexagons of ${\mathcal A}_2(\SU{3})$, \ie $\SU{3}$ at level $2$.  The chosen root is located where the ``6'' stands.}
\end{figure}

\begin{figure}[htbp]
\centering
\includegraphics[width=12pc]{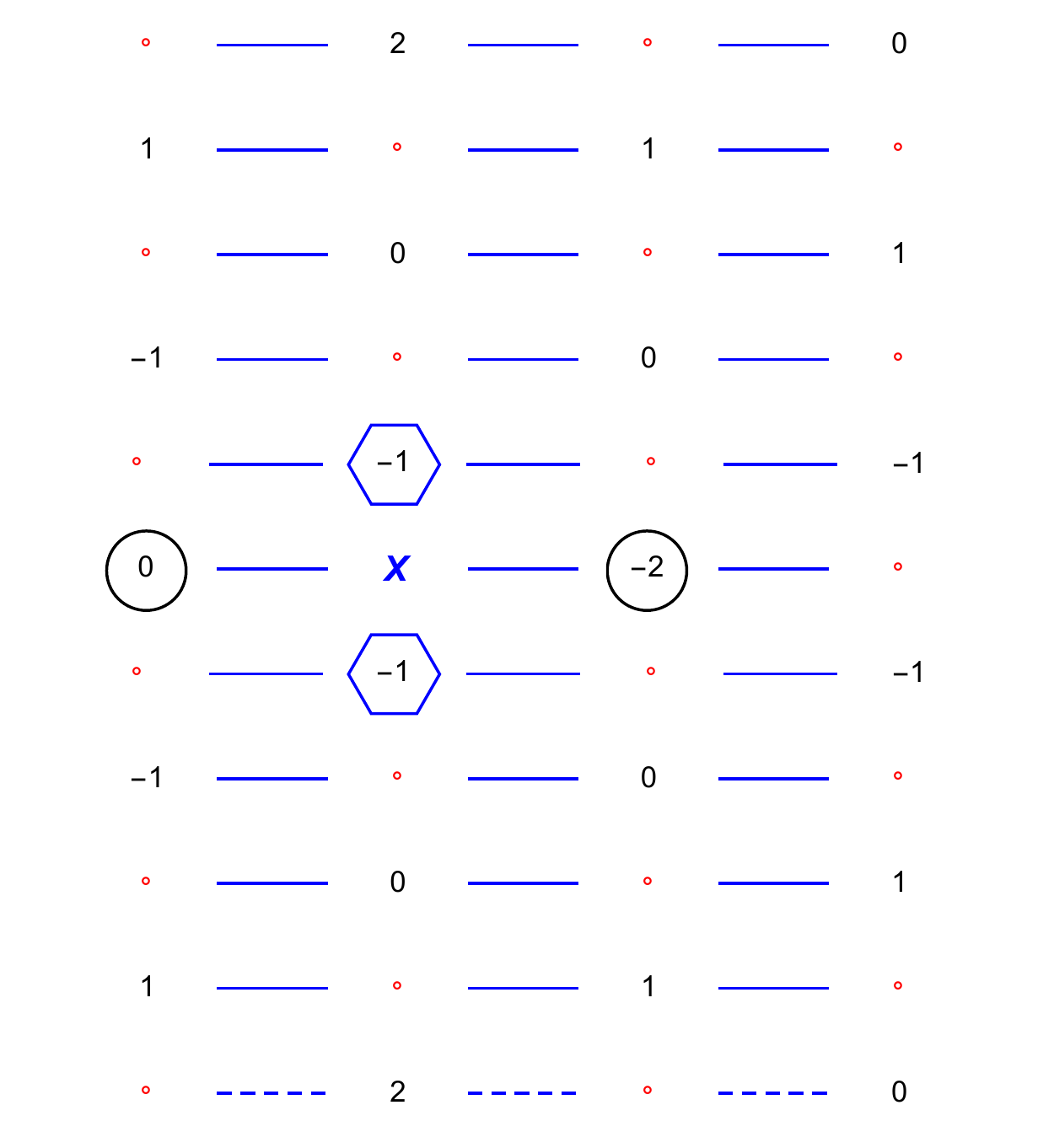}
\caption{\label{A1rootperiod5harmonicity} 
Harmonicity property for one of the $12$ roots of ${\mathcal A}_3(\SU{2})$, \ie $\SU{2}$ at level $3$. The period (read vertically) is $2 \times 5$ since $N=g+k=2+3=5$. Horizontally one recognizes the Dynkin diagram $A_4$ of SU(5).
 The chosen root is located where the ``2'' stands.}
\end{figure}

The harmonicity property  ($\Delta_\Lambda f = \Delta_{\mathcal E} f$, see section~\ref{sec:harmonicity}) of the chosen hyper-root, call it $f$, can be easily checked in this hexagon, see  figure~\ref{A2roohexagon5harmonicity}: consider for instance the point marked $X$, it belongs to a particular fusion diagram and there are two oriented edges (arrows) ending on $X$ in the same diagram, they start from vertices where $f$ have values $2$ and $2$, so that $(\Delta_{\mathcal E} f) (X) = 2+2 = 4$;
there are also oriented edges in the hexagon (or in the $\SU{3}$ weight lattice) that connect the fusion diagrams themselves, they follow the three orientations given in the figure, and those ending in the fusion diagram where $X$ is located therefore define three arrows with head $X$, starting from three neighbouring fusion diagrams, from vertices where $f$ has values $6$, $-1$, $-1$, so that $(\Delta_\Lambda f) (X) = 6-1-1=4$.\\
For illustration and comparison with the above $\SU{3}$ harmonicity, we display on figure~\ref{A1rootperiod5harmonicity} the corresponding  property for the simpler case ${\mathcal A}_k(\SU{2})$ with $k=3$;  the fusion graph coincides therefore with the Dynkin diagram $A_4$, and the $\SU{2}$ harmonicity property for this quiver of roots, which is the SU(5) quiver,  can be checked for instance at the point marked $X$ on figure~\ref{A1rootperiod5harmonicity}. This $\SU{2}$ harmonicity property, which holds for all simple Lie groups, is also illustrated (in particular for the roots of several exceptional Lie groups) in one of the sections of ref \cite{RC:periodicquivers}.

\begin{figure}[ht]
\centering{\includegraphics[width=16cm]{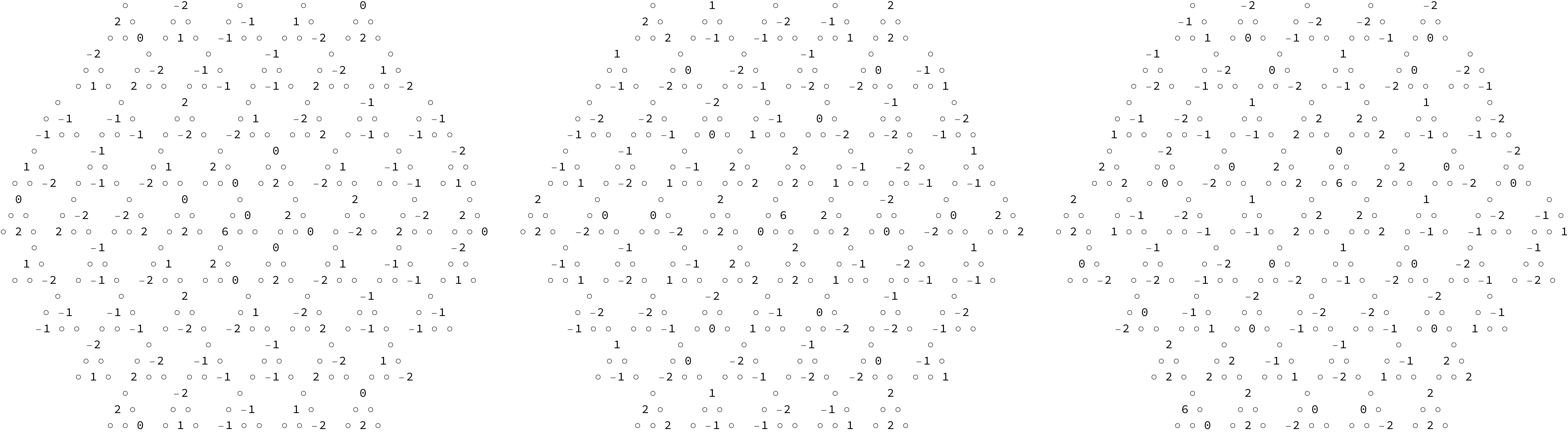}}\\
 \vspace*{1cm}
\centering{\includegraphics[width=16cm]{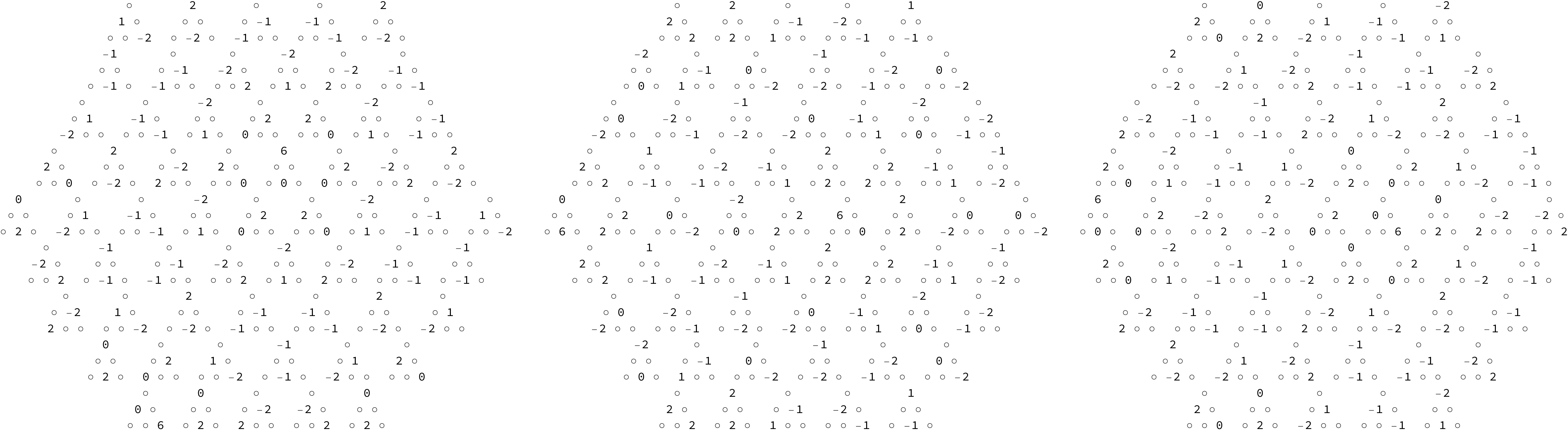}}\\
 \vspace*{1cm}
\centering{\includegraphics[width=16cm]{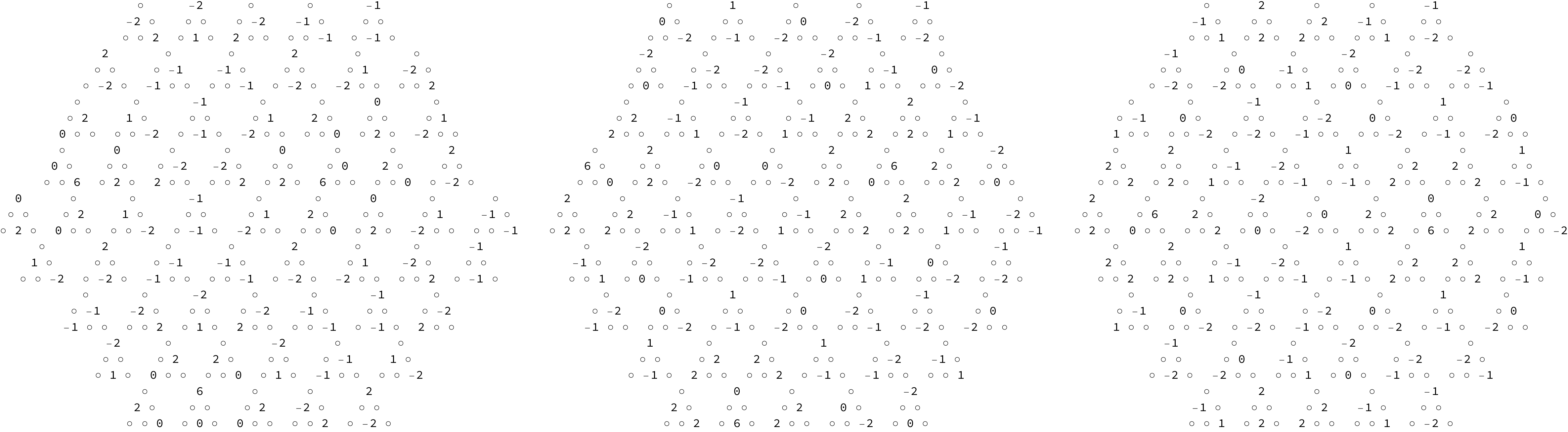}}\\
 \vspace*{1cm}
\centering{\includegraphics[width=16cm]{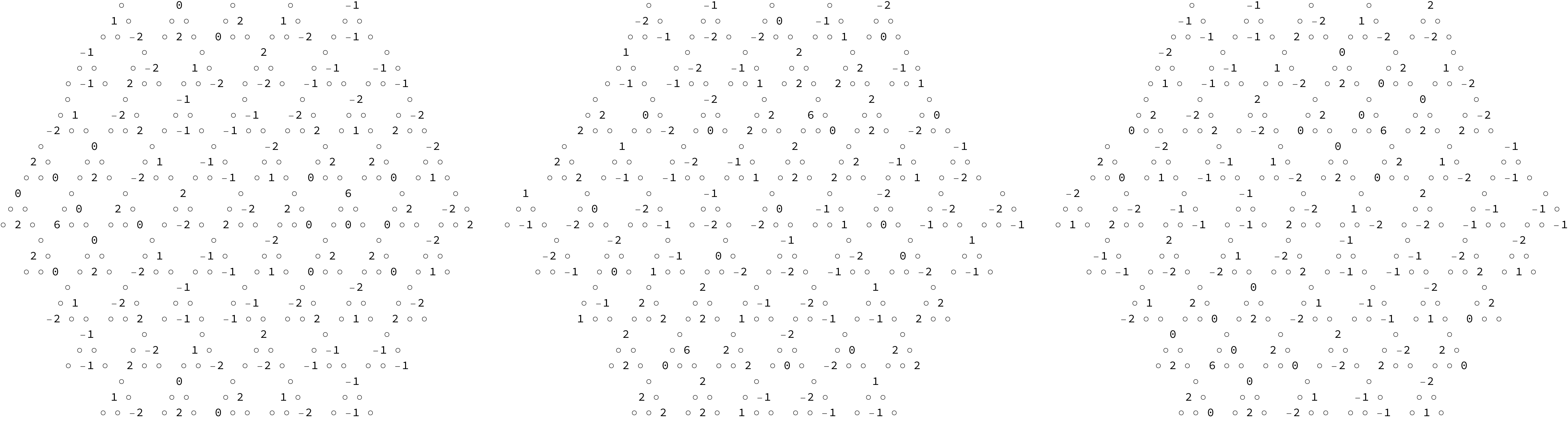}}
\caption{\label{simpleA2relativeHexagons} Twelve relative hexagons associated with a basis of hyper-roots for $\mathcal{A}_2(\SU{3})$. 
There are 50 positive (\ie restricted) hyper-roots and therefore also 50 relative hexagons.
Those displayed here correspond to the  basis ${\mathcal B}_2$.}
\end{figure}

{\bf Absolute and relative hexagons}.\\
The {\sl absolute hexagon} associated with $a$, a vertex of the fusion diagram,  is a hexagonal window, with edges of length $N$  (or $N+1$ if the hexagon is extended, as in figures \ref{A2roohexagon5symmetries}, \ref{A2roohexagon5harmonicity}), displaying the inner product between some hyper-root $\alpha = (m,a)$, where $m$ is an element of the weight lattice (or of the periodicity rhombus), and all others, using periodicity, the hexagon being chosen in such a way that the fusion diagram to which $a$ belongs is located at the center of the hexagon (also the origin of $\SU{3}$ weight lattice). This is the case in figures~\ref{A2roohexagon5symmetries}, \ref{A2roohexagon5harmonicity}.
There are of course $r_{\mathcal E} = 6$ such absolute hexagons if $k=2$.\\
 {\sl Relative hexagons} are also hexagonal windows displaying the inner products between one chosen hyper-root $\alpha = (m,a)$ and all others, they are not usually centered (in the sense that $\alpha$ does not always belong to the fusion diagram located at the center of the hexagon), but they are in good relative positions:  the vertex $a$ belongs to a fusion diagram itself located at the weight $m$ of the $\SU{3}$ weight lattice. By definition there are as many relative hexagons as there are hyper-roots, and if one makes the choice of a basis one can in particular consider the relative hexagons associated with these basis elements.
Relative hexagons are still symmetric (Weyl axes) with respect to the position of the chosen hyper-root, the position of the fusion diagram to which the vertex labelled $6$ belongs, since $<\alpha, \alpha >=6$, but this diagram is not necessarily located at the center of the hexagon. Since they are in good relative positions, relative hexagons can be added (pointwise) or multiplied by scalars; the resulting hexagons display arbitrary (non-necessarily integral) hyper-weights since they automatically obey the required harmonicity properties. Still with the same example $k=2$, we give in figure~\ref{simpleA2relativeHexagons} the relative root hexagons associated with the choice of a particular basis of hyper-roots ($12$ of them in this particular case).

\begin{figure}[ht]
\centering
\includegraphics[width=32pc]{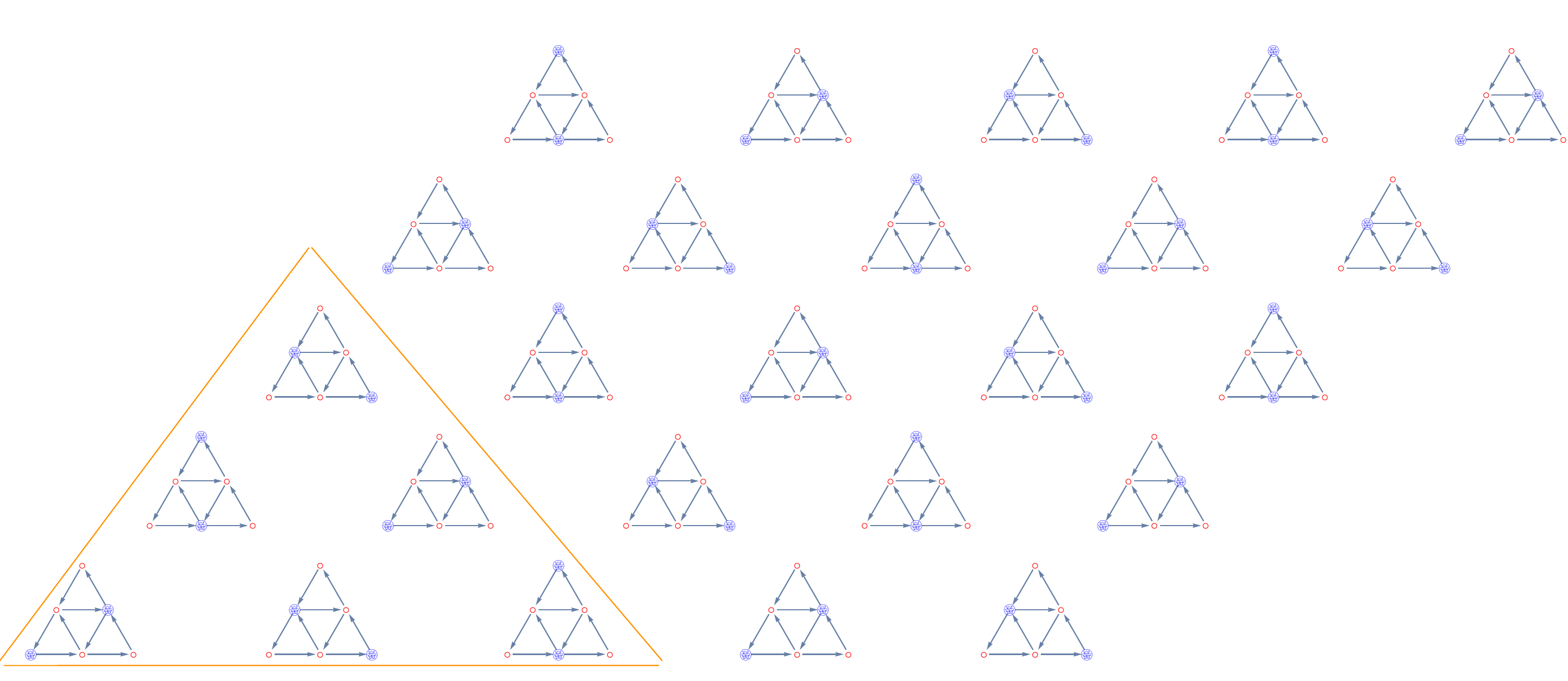}
\caption{\label{allrootpositionA2} $\SU{3}$ at level $2$: position of all hyper-roots in a periodicity rhombus).
The $12$ elements belonging to the triangle traced out in the left bottom can be used, after choosing some ordering, to build a  basis ${\mathcal B}_1$ of the space of hyper-roots.
This basis is used to write the Gram matrix given in equation~\ref{GramA2}.
}
\end{figure}

\begin{figure}[ht]
\centering
\includegraphics[width=32pc]{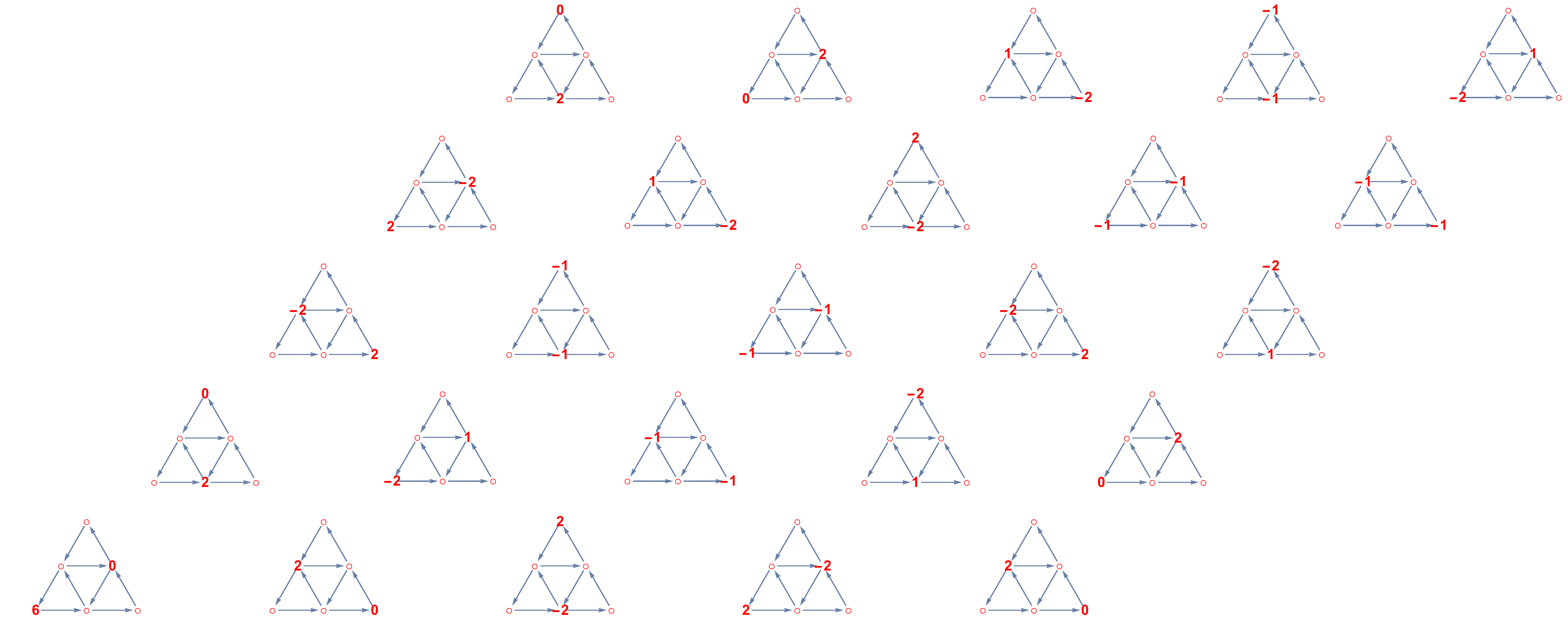}
\caption{\label{rootpositionA2alpha01} $\SU{3}$ at level $2$: the periodicity rhombus associated with one hyper-root (marked 6).}
\end{figure}

{\bf Position of hyper-roots and periodicity rhombus}.\\
Rather than displaying root hexagons, it is enough, albeit slightly less convenient, to display the periodicity rhombus. To each hyper-root one can associate such a rhombus. We still consider the example $k=2$ and display on figure~\ref{allrootpositionA2} the {\sl positions} of the $50$ restricted hyper-roots and, on figure~\ref{rootpositionA2alpha01}, one of the rhombuses associated with some chosen hyper-root (as usual the latter sits where the "6" is). 
The $12$ hyper-roots whose positions are included in the triangle located in the left corner (brown triangle in the picture) define, up to some chosen ordering, a basis of the space of the space of hyper-roots; 
there is nothing special about this basis (we did not introduce any notion of ``simple hyper-roots'')  but it will be used later to define a particular Gram matrix for the hyper-root lattice, and we give in figure~\ref{simpleA2relativeHexagons} the $12$ relative hexagons associated with this basis.

\section{Needed tools for lattices of hyper-roots of type $\SU{3}$}

Root systems defined by Dynkin diagrams (or simple Lie groups), are well known, and the corresponding lattices which are just, in our framework, hyper-root lattices of type $G=\SU{2}$, are described in many places.
Their associated lattices, as well as their theta functions, can be found in the literature, for instance in \cite{ConwaySloane}. 
In most cases they are explicitly given in terms of combinations of elliptic theta functions, but they could also be obtained from a method that uses the theory of modular forms twisted by appropriate Dirichlet characters (although this is not usually done).
This latter method, that we shall review in section \ref{ThetaFunctionsTheoryForLattices}, will be used to find explicit expressions for the lattices of hyper-roots of type $G=\SU{3}$.

\subsection{On the SU(3) classification (reminders)}
A given module ${\mathcal E}$ over ${\mathcal A}_k(\SU{3})$, and consequently a given lattice of hyper-roots, is fully specified by one of the two fundamental fusion matrices of $\SU{3}$
or equivalently by the fusion graph describing the action of one fundamental representation of $\SU{3}$ on the chosen module.
The matrices $F_{(1,0)}$ that we need are recalled in appendix~\ref{AppendixFusionMatrices}, their associated fusion graphs are displayed alongside the headings of section \ref{resultsforlattices}.

The classification of modules ${\mathcal E}$ is known, properties of the members of the different series and of the seven exceptional cases of the $\SU{3}$ family, together with their fusion graphs, 
are discussed in a number of places, see \cite{DiFrancescoZuber}, \cite{Ocneanu:Bariloche}, \cite{YellowBook}, \cite{CoquereauxSU3Maroc}, \cite{CoquereauxSchieberJMP},  \cite{EvansPughSU3}, see also \cite{RCsiteWebFusionGraphs}. 

In the following we shall mostly concentrate on the series ${\mathcal E} = {\mathcal A}_k(\SU{3})$ and call $L_k$ the corresponding lattices of hyper-roots; we consider explicitly the cases $k=1,2,3,4$.
We also give explicit results for theta functions associated with the modules ${\mathcal D}_3$, ${\mathcal D}_6$, 
and the three exceptional cases ${\mathcal E}_5$, ${\mathcal E}_9$ and  ${\mathcal E}_{21}$;
these modules have a non-trivial $Z_3$ grading (our previous general discussion should be slightly modified for modules that do not have a non-trivial $Z_3$ grading, this is why we do not give explicit results for such cases but the method would be identical).
Only  ${\mathcal A}_k$, ${\mathcal D}_k$ with $k=0\; \text{mod}\;  3$ and the three exceptional cases just mentioned  have ``self-fusion'',  or ``are flat'' (in operator algebra parlance), this notion will not be used in this article but these are the examples that we explicitly consider here. \\
Remember that, apart from the ${\mathcal A}_k$ themselves, the following modules have a non-trivial $Z_3$ grading: the ${\mathcal D}_k$ series with $k=0\; \text{mod}\;  3$, the ${\mathcal D}_k^\star$ series, the exceptional ${\mathcal E}_5$, ${\mathcal E}_9$ and  ${\mathcal E}_{21}$, the twisted ${\mathcal D}$ cases ${\mathcal D}_9^t$ and its own module ${\mathcal D}_9^{t \star}$ (these last two are often also flagged as ``exceptional'', the first being an analog of the $E_7$ of the $\SU{2}$ family, which  indeed appears as a twisted $D_{10}= {\mathcal D}_{16}(\SU{2})$), and the exceptional module ${\mathcal M}_9$  (which is also a module over ${\mathcal E}_9$). 
The other members of the $\SU{3}$ classification, namely the cases ${\mathcal D}_k$ with $k=1\; \text{or}\; 2 \; \text{mod}\;  3$,  ${\mathcal A}_k^\star$, and the exceptional module ${\mathcal M}_5$ over ${\mathcal E}_5$ have  trivial  $Z_3$ gradings.

As usual, the subindex $k$ in the above script capital letters indicates the existence of a module structure over ${\mathcal A}_k$ (for instance ${\mathcal D}_4^\star$ is a module over ${\mathcal A}_4$).
In the previous sections the notation ${\mathcal E}_k$ was used in a generic way,  but from now on we use specific notations to denote the module-categories of the $\SU{3}$ classification (all of them appear in the above lists), and we therefore keep the ``${\mathcal E}$'' notation to refer to the three exceptional cases ${\mathcal E}_5$, ${\mathcal E}_9$ and  ${\mathcal E}_{21}$. We hope that there should be no confusion.
Remember also that the Dynkin notation used for the $\SU{2}$ classification does not agree with the above convention: the subindex of a Dynkin diagrams refers to the number of simple objects, whereas the subindex of ${\mathcal A}$ (or ${\mathcal D}$, or  ${\mathcal E}$, \etc)  used in higher classifications usually refers to the level. For instance one has $A_{11} = {\mathcal A}_{10}(\SU{2})$,  $D_{7} = {\mathcal D}_{10}(\SU{2})$,  $E_{6} = {\mathcal E}_{10}(\SU{2})$.

\subsection{Choice of a basis}
\label{choicebasis}
There are many ways of choosing a basis for a lattice $L$. To every choice is associated a fundamental parallelotope and 
a Gram matrix $A$ (the matrix of inner products in this basis).
Two Gram matrices $A$ and $A^\prime$ give rise (or come from) congruent lattices when they are integrally equivalent, \ie when there exists a matrix $U$, with integer entries and determinant $\pm 1$, such that $A^\prime=U A U^{T}$. It is clear that the discriminant is an invariant for integral equivalence. It is also the square volume of a fundamental parallelotope and it is equal to the order of the dual quotient 
$L^\star/L$ where $L^\star$ is the dual lattice.
It is however sometimes useful to loosen a little bit the notion of equivalence and use rational equivalence rather than integral equivalence; this amounts to identify lattices associated with rationally equivalent Gram matrices (the matrix $U$ has rational elements, and its determinant is a non - zero rational number).

In the case of the $\SU{2}$ hyper-root systems (\ie root systems in the usual sense), one may choose for Gram matrix $A$ the Cartan matrix corresponding to a given Coxeter-Dynkin graph ${\mathcal E}$, namely 
$A=2 - F_{(1)}$ where $F_{(1)}$ is the fundamental fusion matrix of the module defined by ${\mathcal E}$.
For usual root lattices,  the notion of Cartan matrix (associated with a basis of simple positive roots) is unique, but one can find  for these lattices other integral basis and other Gram matrices than those associated with simple positive roots. 
For lattices of hyper-roots we did not define fundamental hyper-weights and did not define simple hyper-roots either: the notion of ``Cartan matrix'' is not available. 
However we can certainly consider several interesting choices for the Gram matrices. In what follows we shall usually present only one Gram matrix, called $A$, since it defines the lattice up to integral equivalence, but one should keep in mind that other choices are possible.

\paragraph{Remark.}
 Warning: a naive generalization of the equation $A = 2 - F_{(1)}$ that relates the adjacency matrix of Dynkin diagrams, and therefore also, in the simply-laced case, the fundamental fusion matrices of the $\SU{2}$ modules to the Cartan matrix $A$
 and to the usual Lie group root lattices,  suggests, in the case of $\SU{3}$,  to replace the Cartan matrix $A$ by $6 - (F_{(1,0)}+F_{(0,1)})$ \ie using the fundamental fusion matrices for $\SU{3}$ modules (nimreps). However the lattices obtained from this naive choice\footnote{Some properties of this matrix and of its inverse are investigated in one section of \cite{CoquereauxZuberNuclPhys}, see also \cite{CoquereauxSchieberJMP}.} of $A$ {\sl do not correspond} to the lattices of hyper-roots considered in the present paper; it is already clear that the dimensions do not match: the rank of the lattice of hyper-roots associated with ${\mathcal E}_k(\SU{3})$ is $2 r_{\mathcal E}$ whereas it would be only equal to $r_{\mathcal E}$ for the above naive choice.
 
 \paragraph{Basis ${\mathcal B}_1$.}
 
\begin{figure}[ht]
\centering
\includegraphics[width=28pc]{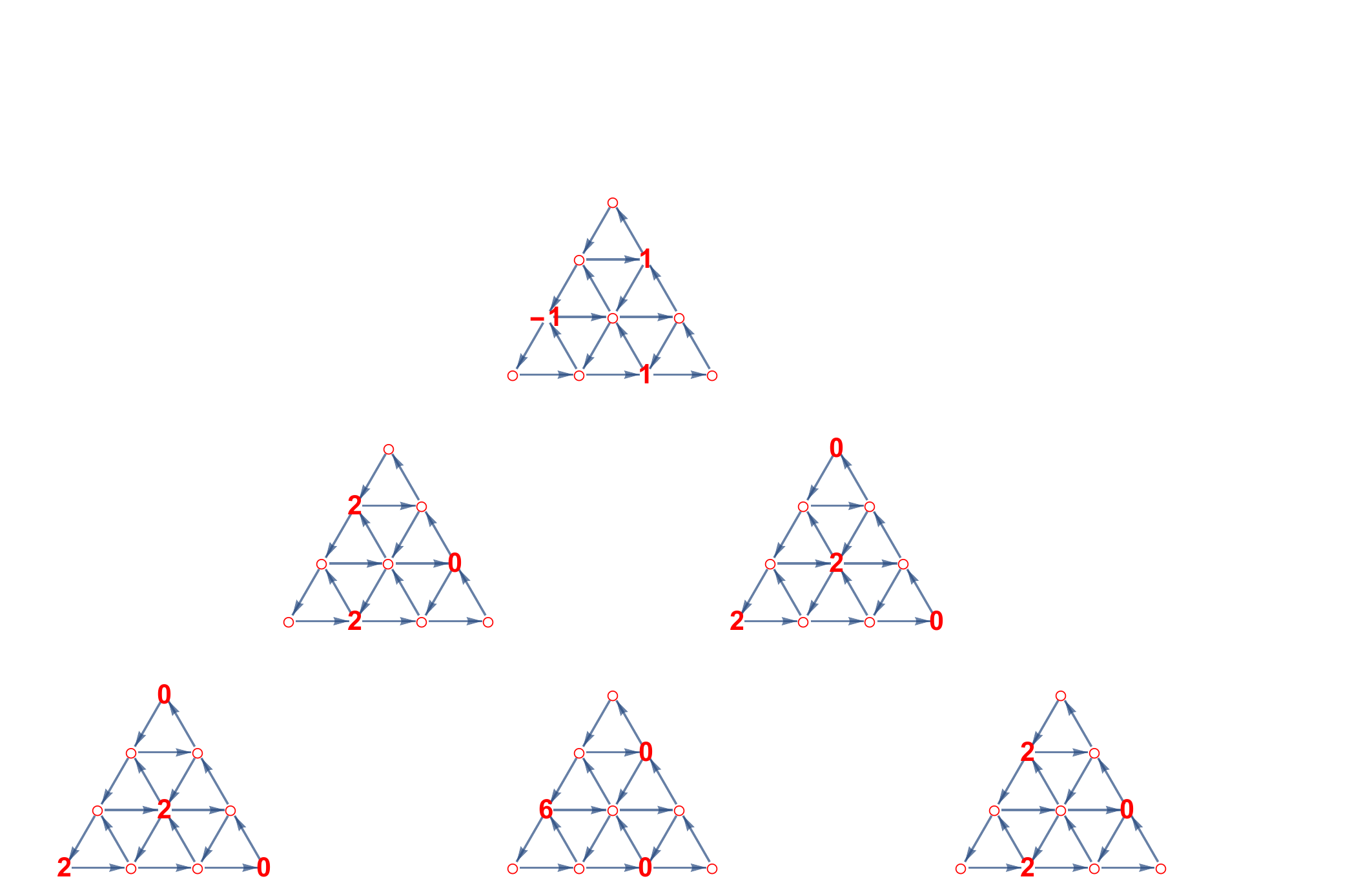}
\caption{\label{A3scalarproductssimplerootsA3} $\SU{3}$ at level $3$: the basis ${\mathcal B}_1$ in the periodicity rhombus ($20$ positions marked with integers). We display the scalar product between the hyper-root marked $6$ and all other basis elements (these values appear on line $5$ of the Gram matrix $A$ given in equation~\ref{GramA3}).}
\end{figure}

Its $2 r_{\mathcal E}$ elements (assuming $k>0$) belong to the bottom left corner of the $\SU{3}$ period parallelogram, more precisely,  we choose those hyper-roots located in the admissible vertices of the six fusion graphs sitting in positions 
 $\{\{0,0\},\{0,1\},\{1,0\},\{1,1\},\{2,0\},\{0,2\}\}$ of the weight lattice; one checks that this indeeds determines a basis which is fully specified once an ordering has been chosen. 
This basis corresponds to the highlighted triangle in figure~\ref{allrootpositionA2}, see also figure~\ref{A3scalarproductssimplerootsA3}.

\paragraph{Other basis.}
Many other convenient basis can be chosen, for instance  ${\mathcal B}_2$ and ${\mathcal B}_3$, respectively associated with the admissible vertices belonging to the fusion graphs located in positions
$\{ \{1, 1\}, \{2, 1\}, \{1, 2\}, \{3, 1\}, \{2, 2\}, \{1, 3\}\}$ and $\{\{0,0\},\{1,0\},\{0,1\},\{N-1,N-2\},\{N-2,N-1\},\{N-1,N-1\}\}$, with $N=k+3$.

\paragraph{Our choice.} With the exception of the lattice $L_1$ where, for illustration purposes, we shall present two equivalent Gram matrices respectively associated with the basis ${\mathcal B}_1$ and ${\mathcal B}_2$, 
the Gram matrices $A$ that will be given later are obtained from equation~\ref{scalarproductofsu3roots} using the basis ${\mathcal B}_1$.
 
\subsection{Summary of the procedure.}
\begin{itemize}
\item Choose a module ${\mathcal E}$ over  ${\mathcal A}_k(\SU3)$, for instance ${\mathcal A}_k$ itself.
\item From the fundamental fusion matrix $F_{(10)}$ of the chosen module given in Appendix, calculate the other fusion matrices $F_n$,  for instance using the SU(3) recurrence relation.
\item Extend the fusion matrices to the weight lattice of SU(3), using symmetries and periodicity.
\item It is useful to build the periodic essential matrices $\tau_a$, and not only the $F_n$, in particular if the module ${\mathcal E}$ is not ${\mathcal A}_k$ itself.
\item Using equation \ref{scalarproductofsu3roots} one can determine a table $A^{\mathrm{big}}$ of the $\vert {\mathcal R}\vert \times \vert {\mathcal R}\vert$ scalar products between all the hyper-roots (or only between those of ${\mathcal R}^\vee$).
The matrix $A^{\mathrm{big}}$ has rank $\mathfrak{r}  = 2 r_{\mathcal E}  <  \vert {\mathcal R}\vert =  (k+1)(k+2)(k+3)^2 /3$.
\item Select a family $(\alpha_i)$ of $\mathfrak{r}$ independent hyper-roots (\ie choose a basis) and call $A$ the $\mathfrak{r}\times \mathfrak{r}$ restriction of the previous table $A^{\mathrm{big}}$ to the chosen basis.
This $A$ will be a Gram matrix for the lattice of hyper-roots. However $A^{\mathrm{big}}$ can be huge.  It is shorter to proceed as follows:
use one of the hyper-root basis (for instance ${\mathcal B_1}$) described previously,  and determine the corresponding $A$ matrix by calculating only the $\mathfrak{r}\times \mathfrak{r}$ inner products between its basis elements.\\
One ends up with a Gram matrix $A$ that is, of course, basis dependent. 
The rest of the discussion is standard in the sense that it mimics what is done for usual roots and weights. Note: the above steps could {\sl also} be followed in that case, just replacing $\SU{3}$ by $\SU{2}$.\\
\item The choice of $A$ determines a basis $(\alpha_i)$ of hyper-roots that are such that $< \alpha_i, \alpha_j> = A_{ij}$.\\
Call $K=A^{-1}$ the inverse of $A$ and $(\omega_i)$ the dual basis of $(\alpha_i)$, then  $<\omega_i, \omega_j> = K_{ij}$ and $<\alpha_i, \omega_j> = \delta_{ij}$.
The family of vectors $(\omega_i)$ is, by definition, the hyper-weight basis associated with the hyper-root basis $(\alpha_i)$ --- there is no need to introduce co-roots or co-weights, since, for the systems considered here, all hyper-roots have the same norm, equal to $6$.\\
Warning: the indices $i,j\ldots$ of $(\alpha_i)$ or $(\omega_i)$ run from $1$ to ${\mathfrak r}=2 r_{\mathcal E}$ whereas the indices $a,b$ of $(\tau_a)$ refer to the irreps of ${\mathcal E}$ and therefore
run only from $1$ to $r_{\mathcal E}$. In particular essential matrices $\tau_a$ and elements $\omega_i$ of the  hyper-weight basis are distinct quantities.\\
Arbitrary linear combinations of the vectors  $(\omega_i)$ with integer coefficients are (integral) hyper-weights, by definition. They have integer scalar products with hyper-roots and they are harmonic functions on the ribbon. 
Hyper-roots are particular hyper-weights.
\item The last step is to study the lattice of hyper-roots and its theta function. How this is done will be described in the next section.
\end{itemize}

One can a posteriori check that the orthonormal projection of a Dirac measure on the ribbon on the subspace of harmonic functions (hyper-weights) is indeed a hyper-root.
This could have been used as a method to determine them. 
Take $\delta_u = (0,0,\ldots, 0, 1, 0, \ldots 0)$ with $\vert {\mathcal R}^\vee\vert$ components and a single $1$ in position $u$, the other components being $0$'s (so this is a Dirac measure on  $\CC^{{\mathcal R}^\vee}$);
its projection can be decomposed on the basis  $(\omega_i)$: $P_u= \sum_j P_{uj}\,  \omega_j$, where the $P_{uj}$  coefficients have to be determined.
Every $\omega_i$, with $i\in\{1\ldots \mathfrak{r}\}$, determines a vector $o_i$ of $\CC^{{\mathcal R}^\vee}$ with components  $< \omega_i, \beta>$ where $\beta$ runs in 
${{\mathcal R}^\vee}$. The projection
$P_u$ also determines a vector $p_u$ of $\CC^{{\mathcal R}^\vee}$ with components $< P_u, \beta>$.
The unknown $P_{uj}$ are obtained by imposing the orthogonality relations $(\delta_u  - p_u). o_i=0$ for all i.
Up to a rescaling factor $N^2$, one checks that the obtained result $P_u$ is indeed one of the hyper-roots, the one localized in position $u$ (where the coefficient $6$ stands).

\subsection{Lattices and theta functions (reminders)}
\label{ThetaFunctionsTheoryForLattices}
We remind the reader a few results about lattices and their theta functions. This material can be gathered from \cite{Zagier:modularforms}.

Consider a positive definite quadratic form $Q$ which takes integer values on $\ZZ^m$. We can write $Q = \frac{1}{2} x^T \, A \, x$, with $x \in \ZZ^m$ and  $A$ a symmetric $m \times m$ matrix. Integrality of $Q$ implies that $A$ is an even integral matrix (its matrix elements are integers and its diagonal elements are even). Therefore $A$ is a positive definite non singular matrix, and $det(A) >0$. So the inverse $A^{-1}$ exists, as a matrix with rational coefficients. 
The {\sl modular level} of $Q$, or of $A$, is the smallest integer $\ell$ such that  $\ell A^{-1}$ is again an even integral matrix -- this notion differs from the notion of conformal level $k$ used in the previous part of this article. $\Delta = (-1)^m \, det(A)$ is  the {\sl discriminant} of $A$.

Given $Q$, one defines the theta function  $\theta_Q (z) = \sum_{n=0}^\infty  \, p(n) \, q^n$ where\footnote{This parameter $q$ is not related to the root of unity, called ${\mathfrak q}$, that appears in section~\ref{fusioncatbasics}.} $q=exp(2 i \pi z)$  and  $p(n) \in \ZZ_{\geq 0}$  is the number of vectors $x\in Z^m$ that are such that $Q(x)=n$.
The function $\theta_Q$ is always a modular form of weight $m/2$.
In our framework $m$ will always be even (in particular $\Delta = det(A)$) so that we set $m=2s$ with $s$ an integer.

The following theorem  (Hecke-Schoenberg) is known \cite{Zagier:modularforms} and will be used: \\
{\sl
Let $Q : \ZZ^{2s} \mapsto \ZZ$ a positive definite quadratic form, integral, with $m=2s$ variables, of level $\ell$ and discriminant $\Delta$. Then the theta function $\theta_Q$ is a modular form  on the group $\Gamma_0(\ell)$, of weight $s$,  and character $\chi_\Delta$.
}
\\
In plain terms:  $\theta_Q(\frac{az+b}{cz+d}) = \chi_\Delta(a) \, (c z + d)^s \, \theta_Q(z)$ for all  $z\in {\mathfrak H}$ (upper half-plane) and 
$ \left( \begin{smallmatrix} a&b\\ c&d \end{smallmatrix} \right) \in \Gamma_0(\ell)$. Here $\Gamma_0(\ell)$ is the subgroup\footnote{As  $\Gamma_1(\ell) \subset \Gamma_0(\ell)$, one can sometimes use 
modular forms (and bases of spaces of modular forms) twisted by Dirichlet characters on the congruence subgroup $\Gamma_1(\ell)$.}
 of $SL(2,Z)$ defined by the condition $c \equiv 0 \, mod \, \ell$
and $\chi_\Delta$ is the unique Dirichlet character modulo $\ell$ which is such that $\chi_\Delta(p) = {\mathfrak L}(\Delta,p)$ for all odd primes $p$ that do not divide $\ell$, where 
${\mathfrak L}$ denotes the Legendre symbol.

Notice that $m$, as defined above, is also, in our framework, the dimension of the space ${\mathfrak C}$ of hyper-roots, which, for $G=\SU{3}$, is equal to $2 r_{\mathcal E}$. 
In that case,  the weight of the (twisted) modular form $\theta_Q$ is therefore equal to $r_{\mathcal E}$, the number of vertices of the fusion diagram, or the number of simple objects in ${\mathcal E} $.

\paragraph{About Dirichlet characters.} Dirichlet characters are particular functions from the integers to the complex numbers that arise as follows:  given a character on the group of invertible elements of the set of integers modulo $p$, one can lift it to a completely multiplicative function on integers relatively prime to $p$ and then extend this function to all integers by defining it to be $0$ on integers having a non-trivial factor in common with $p$. A Dirichlet character with modulus $p$ takes the same value on two integers that agree modulo $p$. The interested reader may consult the abundant literature on the subject but it is enough for us to remember that they are a particular kind of completely multiplicative complex valued functions on the set of integers,  that there are $\phi(p)$ characters modulo $p$, where $\phi$ is the Euler function, and that they are tabulated in many places ---there is even a command DirichletCharacter[p,j,n] in Mathematica \cite{Mathematica} that gives the Dirichlet character with modulus $p$ and index $j$ as a function of $n$ (the index $j$ running from $1$ to $\phi(p)$).

\section{Lattices of hyper-roots of type $\SU{3}$ and their theta functions}
\label{resultsforlattices}
\subsection{The hyper-root lattice $L_0$ of  ${\mathcal A}_0(\SU{3})$}

For $k=0$, the general formulae of section~\ref{sec:ribbon} give $N=3$, $r_{{\mathcal E}_0}=1$,  $\mathfrak{r} = 2$, $\vert {\mathcal R}^\vee \vert =3$ and $\vert {\mathcal R}\vert = 6$.
One expects\footnote{Warning: A simple counting argument shows that the lattice of hyper-roots obtained by taking $k=0$ for $G=\SU{n}$ and $n>3$ cannot be identified with the usual root lattice of SU(n).} that this lattice should coincide with the usual root lattice of $\SU{3}$, and ${\mathcal R}^\vee$ with the set of positive roots.
 The period is a rhombus $3\times 3$ but the small fusion graph has a single vertex, which $Z_3$ grading $0$, so, 
 in order to build a basis of the hyper-root lattice, only two of the six weights $\{\{1,1\},\{2,1\},\{1,2\},\{3,1\},\{2,2\},\{1,3\}\}$ contribute (see section~\ref{choicebasis}), those of grading $0$, namely $\{1,1\}$ and $\{2,2\}$;
 we therefore recover that the dimension is $\mathfrak{r} = 2$. The Gram matrix of the lattice of hyper-roots, in this basis, obtained from equation~(\ref{scalarproductofsu3roots}), is
 $\left(\begin{smallmatrix}6&-3\\-3&6\end{smallmatrix}\right)$, \ie three times the Cartan matrix of $\SU{3}$, this lattice can therefore be considered as a rescaled version of the  $\SU{3}$ root lattice (the hexagonal lattice).
 The theta function of the latter is well known and can be found in many textbooks, for instance in \cite{ConwaySloane}; its expression in terms of the elliptic theta function $\vartheta_3$ reads:
 \begin{eqnarray*}
 \theta(z) &=& \frac{\vartheta _3(0,q){}^3+\vartheta _3\left(\frac{\pi }{3},q\right){}^3+\vartheta
   _3\left(\frac{2 \pi }{3},q\right){}^3}{3 \vartheta _3\left(0,q^3\right)}\cr
   {}&=&1+6 \, q^2+6 \, q^6+6 \, q^8+12 \, q^{14}+6 \, q^{18}+6 \, q^{24}+12 \, q^{26}+6 \, q^{32}+O\left(q^{34}\right)
    \end{eqnarray*}
 The theta function of the hyper-root lattice $L_0$  is then simply $\theta(3z)$ ---replace $q$ by $q^3$.
 
 Although this special case ($k=0$) coincides, up to scale, with a well known lattice, it is instructive to look at its theta function by using the theorems recalled in the previous section.
The quadratic form defined by the Cartan matrix of $\SU{3}$ has level $3$, so that its theta function
is a modular form on the group $\Gamma_0(3)$, of weight $s=1$; it is twisted by a non-trivial 
Dirichlet character modulo $3$ (there are only two of them here, the first being trivial), and there is no constraint coming from the Legendre symbol condition since there are no odd primes that do not divide $3$. This vector space of modular forms is of dimension $1$, hence $\theta$ can be identified with its generator.  
As an application, here is a very short program, using the computer algebra package Magma \cite{Magma} that returns the above theta function, up to the same order $(q^{17})^2$ ---one has to rescale $q$--- and uses the above concepts:

\begin{verbatim}
H := DirichletGroup(3,CyclotomicField(EulerPhi(3)));
chars := Elements(H); eps := chars[2];
M := ModularForms([eps],1); Basis(M,17);
\end{verbatim}
 
\subsection{The hyper-root lattice $L_1$ of  ${\mathcal A}_1(\SU{3})$ \hspace{5.cm}\protect\icon{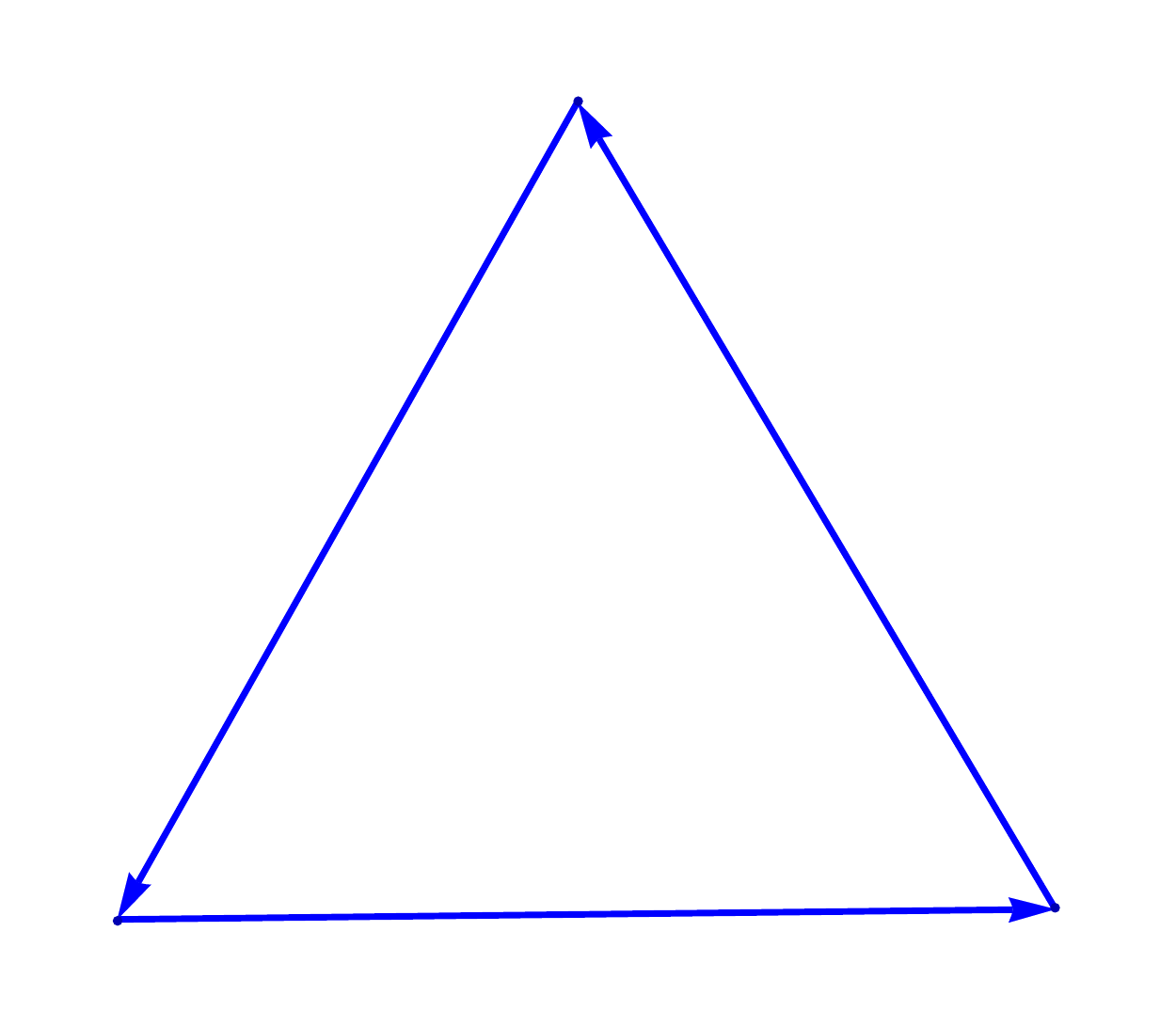}}
\label{secL1}
For $k=1$, $N=4$, $r_{{\mathcal A}_1} =3$. The rank of the lattice $L_1 $ is $\mathfrak{r} = 2\, r_{{\mathcal A}_1}  = 6$.
The period is a rhombus $4\times 4$ and  the number of hyper-roots is $\vert {\mathcal R}\vert = 32$, with $\vert {\mathcal R}^\vee\vert = 16$.
This last number being reasonably small,  we shall give more details for this lattice than for those that come after.

\subsubsection{Gram matrices}
\label{GramA1}

There are many possible Gram matrices for this lattice: they differ by the choice of the basis (integral equivalence).
They have a determinant equal to $4^{6}$.  The lattice is even, with minimal norm~$6$. 
For illustration, three possible Gram matrices denoted $A, A^\prime, A^{\prime \prime}$, are given below.
The first two are respectively associated with the basis choices ${\mathcal B}_1$ and ${\mathcal B}_2$ described in section~\ref{choicebasis}.
The third simply relates to a fundamental fusion matrix of ${\mathcal A}_1$ (see a comment in section~\ref{simpleGram}).
The matrix $K$ is the inverse of $A$. 
We denote $\{{\alpha_i}\}, i=1\ldots 6$, the elements of the hyper-root basis ${\mathcal B}_1$, \ie $<\alpha_i, \alpha_j> = A_{ij}$.
The members of the dual basis (the hyper-weight basis) are denoted $\{{\omega_j}\}$, so that  $<\omega_i,  \alpha_j> = \delta_{ij}$ and $<\omega_i, \omega_j> =K_{ij}$.
$$
\begin{array}{ccc}
A=
\left(
\begin{array}{cccccc}
 6 & 2 & 2 & -2 & -2 & -2 \\
 2 & 6 & 2 & 2 & -2 & 2 \\
 2 & 2 & 6 & 2 & 2 & -2 \\
 -2 & 2 & 2 & 6 & 2 & 2 \\
 -2 & -2 & 2 & 2 & 6 & -2 \\
 -2 & 2 & -2 & 2 & -2 & 6 \\
\end{array}
\right)
& &
A^\prime =
\left(
\begin{array}{cccccc}
 6 & 2 & 2 & -2 & -2 & -2 \\
 2 & 6 & 2 & 2 & 2 & -2 \\
 2 & 2 & 6 & -2 & 2 & 2 \\
 -2 & 2 & -2 & 6 & 2 & -2 \\
 -2 & 2 & 2 & 2 & 6 & 2 \\
 -2 & -2 & 2 & -2 & 2 & 6
\end{array}
\right)
\end{array}
$$

$$
\begin{array}{ccc}
A^{\prime\prime} = 
\left(
\begin{array}{cccccc}
 6 & 2 & 2 & 2 & -2 & -2 \\
 2 & 6 & 2 & -2 & 2 & -2 \\
 2 & 2 & 6 & -2 & -2 & 2 \\
 2 & -2 & -2 & 6 & -2 & -2 \\
 -2 & 2 & -2 & -2 & 6 & -2 \\
 -2 & -2 & 2 & -2 & -2 & 6 \\
\end{array}
\right)

&&

K= \frac{1}{8} \;
\left(
\begin{array}{cccccc}
 3 & -1 & -1 & 1 & 1 & 1 \\
 -1 & 3 & -1 & -1 & 1 & -1 \\
 -1 & -1 & 3 & -1 & -1 & 1 \\
 1 & -1 & -1 & 3 & -1 & -1 \\
 1 & 1 & -1 & -1 & 3 & 1 \\
 1 & -1 & 1 & -1 & 1 & 3 \\
\end{array}
\right)
\end{array}
$$

The $16$ elements of ${\mathcal R}^\vee$ may be called ``positive hyper-roots'' (their opposites, the elements of $-{\mathcal R}^\vee$, being ``negative'') and they can be expanded on the chosen root basis as follows:
\begin{equation*}
\begin{split}
\alpha _ 1 - \alpha _ 2 + \alpha _ 6,\,&
 -\alpha _ 2 + \alpha _ 3 - \alpha _ 5,\, -\alpha _ 1 + \alpha _ 3 - \alpha _ 4,\, -\alpha _ 4 + \alpha _ 5 + \alpha _ 6,\, \alpha_ 6,\, -\alpha _ 2 + \alpha _ 4 - \alpha _ 5,\, -\alpha _ 1 - \alpha _ 5 - \alpha _ 6,  \\
  -\alpha _ 1 + \alpha _ 2 - \alpha _ 4,&
   \, \alpha _ 2,\, \alpha _ 4,\, -\alpha _ 3 + \alpha_ 4 - \alpha _ 6,\, \alpha _ 2 - \alpha _ 3 - \alpha _ 6, \, \alpha _ 1,\, \alpha _ 3,\, \alpha _ 5,\, \alpha _ 1 - \alpha _ 3 + \alpha _ 5
\end{split}
\end{equation*}
With the same ordering,  the family of their mutual inner products build the following $16\times 16$ matrix  $A^{\mathrm{big}}$, which is of rank $6$, as expected:       
{\scriptsize
$$
A^{\mathrm{big}}=
\left(
\begin{array}{cccccccccccccccc}
 6 & 2 & -2 & 2 & 2 & 2 & -2 & -2 & -2 & -2 & -2 & -2 & 2 & -2 & -2 & 2 \\
 2 & 6 & 2 & -2 & -2 & 2 & 2 & -2 & -2 & -2 & -2 & -2 & 2 & 2 & -2 & -2 \\
 -2 & 2 & 6 & 2 & -2 & -2 & 2 & 2 & -2 & -2 & -2 & -2 & -2 & 2 & 2 & -2 \\
 2 & -2 & 2 & 6 & 2 & -2 & -2 & 2 & -2 & -2 & -2 & -2 & -2 & -2 & 2 & 2 \\
 2 & -2 & -2 & 2 & 6 & 2 & -2 & 2 & 2 & 2 & -2 & -2 & -2 & -2 & -2 & -2 \\
 2 & 2 & -2 & -2 & 2 & 6 & 2 & -2 & -2 & 2 & 2 & -2 & -2 & -2 & -2 & -2 \\
 -2 & 2 & 2 & -2 & -2 & 2 & 6 & 2 & -2 & -2 & 2 & 2 & -2 & -2 & -2 & -2 \\
 -2 & -2 & 2 & 2 & 2 & -2 & 2 & 6 & 2 & -2 & -2 & 2 & -2 & -2 & -2 & -2 \\
 -2 & -2 & -2 & -2 & 2 & -2 & -2 & 2 & 6 & 2 & -2 & 2 & 2 & 2 & -2 & -2 \\
 -2 & -2 & -2 & -2 & 2 & 2 & -2 & -2 & 2 & 6 & 2 & -2 & -2 & 2 & 2 & -2 \\
 -2 & -2 & -2 & -2 & -2 & 2 & 2 & -2 & -2 & 2 & 6 & 2 & -2 & -2 & 2 & 2 \\
 -2 & -2 & -2 & -2 & -2 & -2 & 2 & 2 & 2 & -2 & 2 & 6 & 2 & -2 & -2 & 2 \\
 2 & 2 & -2 & -2 & -2 & -2 & -2 & -2 & 2 & -2 & -2 & 2 & 6 & 2 & -2 & 2 \\
 -2 & 2 & 2 & -2 & -2 & -2 & -2 & -2 & 2 & 2 & -2 & -2 & 2 & 6 & 2 & -2 \\
 -2 & -2 & 2 & 2 & -2 & -2 & -2 & -2 & -2 & 2 & 2 & -2 & -2 & 2 & 6 & 2 \\
 2 & -2 & -2 & 2 & -2 & -2 & -2 & -2 & -2 & -2 & 2 & 2 & 2 & -2 & 2 & 6 \\
\end{array}
\right)
$$
}
The lattice is even, and if we rescale it, setting  $B = A/2$, where $A$ is one of the above Gram matrices, one finds $det (B)=64$, and the vectors of minimal norm belonging to  the lattice (no longer even) associated with the Gram matrix $B$ have norm $3$; however, the lattices $L_k$, for $k>1$ are usually not even.
In all coming sections we shall always choose for the lattices under consideration a basis made of hyper-roots, and the diagonal elements of the associated Gram matrix (keeping in mind an arbitrariness of choice) will therefore always be equal to $6$, that comes from the order of the Weyl group of SU(3).\\
The determinant of $L_1$, equal to $4096$, is sometimes called ``connection index'', it is also the order of the dual quotient $L_1^*/L_1$. The latter is an abelian group isomorphic with $Z_2\times (Z_4)^{ \times 4} \times Z_8$. The  lattice $L_1$ is obviously not self-dual. If we rescale $L_1$ as above in such a way that the minimal norm is $3$, and call $B$ this new lattice, we see that the connection index is $64$ and that the dual quotient $B^*/B$ is isomorphic with $(Z_2)^{ \times 4} \times Z_4$. Elements of the lattice $B^*$ belong to one and only one congruence class, an element of the dual quotient, they are  therefore be classified by 5-uplets $(c_{2_1},c_{2_2},c_{2_3},c_{2_4},c_4)$, with $c_{2_i} \in \{0,1\}$ and $c_4  \in \{0, 1, 2, 3\}$.

\subsubsection{Theta function}

A direct calculation leads to  $$ \theta(z) = 1 + 32 \, q^6 + 60 \, q^8 + 192 \, q^{14} + 252 \, q^{16} + 480 \, q^{22} + 544 \, q^{24} + 
 832 \, q^{30} + 1020 \, q^{32} + 1440 \, q^{38} + 1560 \, q^{40} + 2112 \, q^{46} + \ldots$$
 
 This is in agreement with the following theta series:
$$\frac{1}{2} \left(\vartheta _2\left(0,q^4\right){}^6+\vartheta _3\left(0,q^4\right){}^6+\vartheta _4\left(0,q^4\right){}^6\right)$$

The latter is recognized as the theta function for a (scaled version of) the shifted  $D_6$ lattice, called $D_6^{+} = D_6 \cup ([1] + D_6)$, see \cite{ConwaySloane}.
Notice that two inequivalent lattices may have the same theta series, so the stated coincidence, by itself, is not sufficient to allow identification of $L_1$ and $D_6^{+}$ which, ultimately, relies on the fact, as we shall see below, that one can choose the same Gram matrix to define both lattices.
It is known \cite{ConwaySloane} that the $D_n^{+}$ packing is a lattice packing if and only if $n$ is even. In particular this is so for $n=6$ -- and we know {\sl a priori} that $L_1$  is a lattice and not only a packing. The fact that $D_n^{+}$ is not a lattice for $n$ odd excludes a possible systematic identification with the lattices $L_k$, $k>1$, that are associated with higher hyper-root systems of $\SU{3}$ type.

\smallskip

We may recover the previous theta function for this lattice by applying the Hecke-Schoenberg theorem.
From the Gram matrix one finds that the discriminant  is $4^{6}$ and that the  (modular) level of the quadratic form is $16$. 
The odd primes not dividing $16$ are $3, 5, 7, 11, 13$ and their Legendre symbols are all equal to $1$.
From the $8 \times 16$ table of Dirichlet characters of modulus $16$ over the cyclotomic field of order $\varphi_{Euler}(16)=8$
restricted to odd primes not dividing the level,
one selects the unique character whose values 
coincide with the list obtained for the Legendre symbols.
The space of modular forms on $\Gamma_ 1(16)$ of weight $3$,  twisted by this Dirichlet character,
namely the Kronecker  character -4,  has dimension $7$. It is spanned by the following forms (in the remaining part of this section we set $q_2 = q^2$): 
{\small
\begin{eqnarray*}
  b_1&=&  1 + 12 \, q_2^8 + 64 \, q_2^{12} + 60 \, q_2^{16} + O(q_2^{24}), \qquad
  b_2 =  \, q_2 + 21 \, q_2^9 + 40 \, q_2^{13} + 30 \, q_2^{17} + 72 \, q_2^{21} + O(q_2^{24}),\\
  b_3&=&  \, q_2^2 + 26 \, q_2^{10} + 73 \, q_2^{18} + O(q_2^{24}), \qquad
  b_4 =  \, q_2^3 + 6 \, q_2^7 + 15 \, q_2^{11} + 26 \, q_2^{15} + 45 \, q_2^{19} + 66 \, q_2^{23} + O(q_2^{24}),\\
  b_5&=&  \, q_2^4 + 4 \, q_2^8 + 8 \, q_2^{12} + 16 \, q_2^{16} + 26 \, q_2^{20} + O(q_2^{24}),\quad
  b_6 =   \, q_2^5 + 2 \, q_2^9 + 5 \, q_2^{13} + 10 \, q_2^{17} + 12 \, q_2^{21} + O(q_2^{24}),\\
  b_7&=&  \, q_2^6 + 6 \, q_2^{14} + 15 \, q_2^{22} + O(q_2^{24})
\end{eqnarray*}
}
An explicit determination of the vectors (and their norms) belonging to the first shells of the hyper-roots lattice of ${\mathcal A}_1(\SU{3})$ shows that the theta function starts as $1+32 q_2^3+60 q_2^4 + O(q^{14})$.
The components of this modular form on the previous basis are therefore ${1, 0, 0, 32, 60, 0, 0}$. 
In other words, $$\theta = b_1 + 32 \, b_4 + 60 \, b_5$$
Using a computer package, one can quickly obtain the $q$-expansion of the functions $b_n$ to very large orders and recover or extend the result that was  given for $\theta$.
 As an alternative to the expression of  $\theta$  previously given  in terms of elliptic theta functions, here is a Magma program that returns its series expansion up to order 24 in $q_2$ and uses the above ideas: 
\begin{verbatim}
H := DirichletGroup(16,CyclotomicField(EulerPhi(16))); 
chars := Elements(H); eps := chars[2];
M := ModularForms([eps],3); order:=24;
PowerSeries(M![1,0,0,32,60,0,0],order);
\end{verbatim}

\subsubsection{The automorphism group of the lattice $L_1$}

For $\SU{2}$ hyper-roots (\ie usual roots) the Weyl group is a subgroup of the automorphism group of the lattice. In the case of $\SU{3}$ hyper-roots latices, one can also, in each case, consider the automorphism group $aut$ of the lattice.
Using Magma, we find that the automorphism group of $L_1$ is of order $23040$ and that it is  isomorphic with the semi-direct product of $A_6$ (the alternated group of order $6!/2= 360$) times an abelian group of order  $64$, actually with
$((C_2)^{\times 5} \rtimes A_6)\rtimes C_2$. Orbits of the basis vectors under the $aut$ action coincide and contain the $32$ hyper-roots (the $16$ positive and the $16$ negative ones). 
The group $aut$ is generated by the following matrices 
{\footnotesize
$$
\begin{array}{ccc}
\left(
\begin{array}{cccccc}
 0 & 0 & 0 & -1 & 1 & -1 \\
 -1 & 0 & 0 & -1 & 0 & -1 \\
 -1 & 0 & 1 & 0 & 0 & -1 \\
 0 & 0 & 0 & 0 & -1 & 0 \\
 -1 & 1 & 0 & 0 & -1 & 0 \\
 0 & 0 & 0 & 1 & 0 & 0
\end{array}
\right)
, &
\left(
\begin{array}{cccccc}
 1 & -1 & 0 & 0 & 1 & 0 \\
 0 & -1 & 0 & 0 & 1 & -1 \\
 0 & -1 & 0 & 0 & 0 & 0 \\
 -1 & 0 & 1 & 0 & 0 & -1 \\
 -1 & 0 & 0 & -1 & 0 & -1 \\
 0 & 0 & 0 & 0 & -1 & 0
\end{array}
\right)
, &
\left(
\begin{array}{cccccc}
 1 & 0 & 0 & 0 & 0 & 0 \\
 0 & 1 & 0 & 0 & 0 & 0 \\
 0 & 0 & 0 & -1 & 1 & -1 \\
 0 & 0 & 0 & 1 & 0 & 0 \\
 0 & 0 & 0 & 0 & 1 & 0 \\
 0 & 0 & -1 & -1 & 1 & 0
\end{array}
\right)
\end{array}
$$
}

\subsubsection{Other avatars of the lattice $L_1$}

We already identified the lattice $L_1$ with a scaled version of the shifted  $D_6$ lattice, called $D_6^{+}$. Here are a few others.

\paragraph{The generalized laminated lattice $\Lambda_6[3]$ with minimal norm $3$.} 

It belongs to a family of lattices $\Lambda_n[3]$ that was studied and classified in \cite{PleskenPohst}. 
These lattices have a kind of periodicity, and the relevant information is encoded in a tree of inclusions up to some maximal object that appears to be $\Lambda_{23}[3]$, isomorphic to the so-called ``shorter Leech lattice'' (a nice sublattice of the Leech lattice). The authors construct the tree of inclusions and give the Gram matrices for all  the $\Lambda_n[3]$ of the sequence --- actually they give the Gram matrix for $n=23$ and a few others but this information is sufficient to reconstruct Gram matrices for all of them.  In particular, for $n=6$, one recovers half the matrix $A_2$ already obtained in this section. The fact that $\Lambda_6[3]$ can be identified with  $D_6^{+}$ (and in particular with ``our'' $L_1$) is not mentioned in \cite{PleskenPohst} but the previous observation shows that it is so.

One could again be tempted to identity the lattices $L_k$, defined by the fusion graphs associated with ${\mathcal E} = {\mathcal A}_k(\SU{3})$, for $k>1$, 
with other $\Lambda_n$'s. This is however not the case because the minimal norm of the (unrescaled) $L_k$, for $k > 1$, is $6$, not $3$. Moreover the lattices 
$\Lambda_n[3]$ in dimensions $n=12$ and $n=20$ have kissing numbers respectively equal to $136$ and $1280$ whereas the lattices associated with fusion graphs of $\SU{3}$ at levels $k=2,3$, \ie also in dimensions $12$ and $20$, have kissing numbers $100$ and $240$.

\paragraph{The lattice ${\mathcal L}_4$ generated by cuts of the complete graph on a set of $4$ vertices.}
In \cite{DezaGrishukhin}, the authors study the ``Delaunay Polytopes of Cut Lattices'', \ie the real span of the lattice ${\mathcal L}_n$ generated by cuts of the complete graph on a set of $n$ vertices, which  is a vector space of dimension $({}^ n_2)$. In particular the dimension is $6$ when $n=4$. The authors are interested in the Delaunay polytopes for the lattices ${\mathcal L}_n$. The point is that when $n=4$, this lattice is isomorphic with $D_6^{+}$. Actually, the precise relation is ${\mathcal L}_4 = {\sqrt 2} \, D_6^{+}$.
Let $e_{ij}$,   $1 \leq i < j \leq n$  be an orthonormal basis of $R^{({}^ n_2)}$.  In this basis, 
a vector $d$ of the lattice has coordinates $d_{ij}$.
An integral vector $d \in {\mathcal L}_n$ if and only if  $d_{ij} + d_{jk} + d_{ki} = 0 \,  (mod \,  2) $,  
for all triples  $\{i,j,k\}$. 

Again, one can see that lattices associated with fusion graphs of $\SU{3}$ cannot, in general, be identified with the lattices ${\mathcal L}_n$, unless $n=4$. 
So the above properties (characterization of lattice vectors) hold only for the lattice $L_1$.

\subsubsection{A short description of the Voronoi cells}
The Voronoi cells of the lattice $L_1$ have $588$ vertices ($576$ correspond to shallow holes and $12$ to deep holes).
The deep holes are of maximal norm, equal to $4$ (the covering radius).
The Voronoi polytope has $92$ $5$-dimensional facets, $32$ of them are orthogonal to the vectors of norm $6$ and $60$ are orthogonal to the vectors of norm $8$ (these norms are $3$ and $4$ if one uses the previously mentioned rescaled version of this lattice), $4896$ edges and $588$ vertices. These results can be found from the Gram matrix, for example using Magma \cite{Magma}, and agree with those of \cite{DezaGrishukhin} who, in another framework already studied the Voronoi and Delaunay dual tesselations of the lattice $L_1 \sim {\mathcal L}_4 \sim D_6^{+}$ (see previous paragraph).

\subsection{The hyper-root lattice $L_2$ of  ${\mathcal A}_2(\SU{3})$\hspace{5.cm}\protect\icon{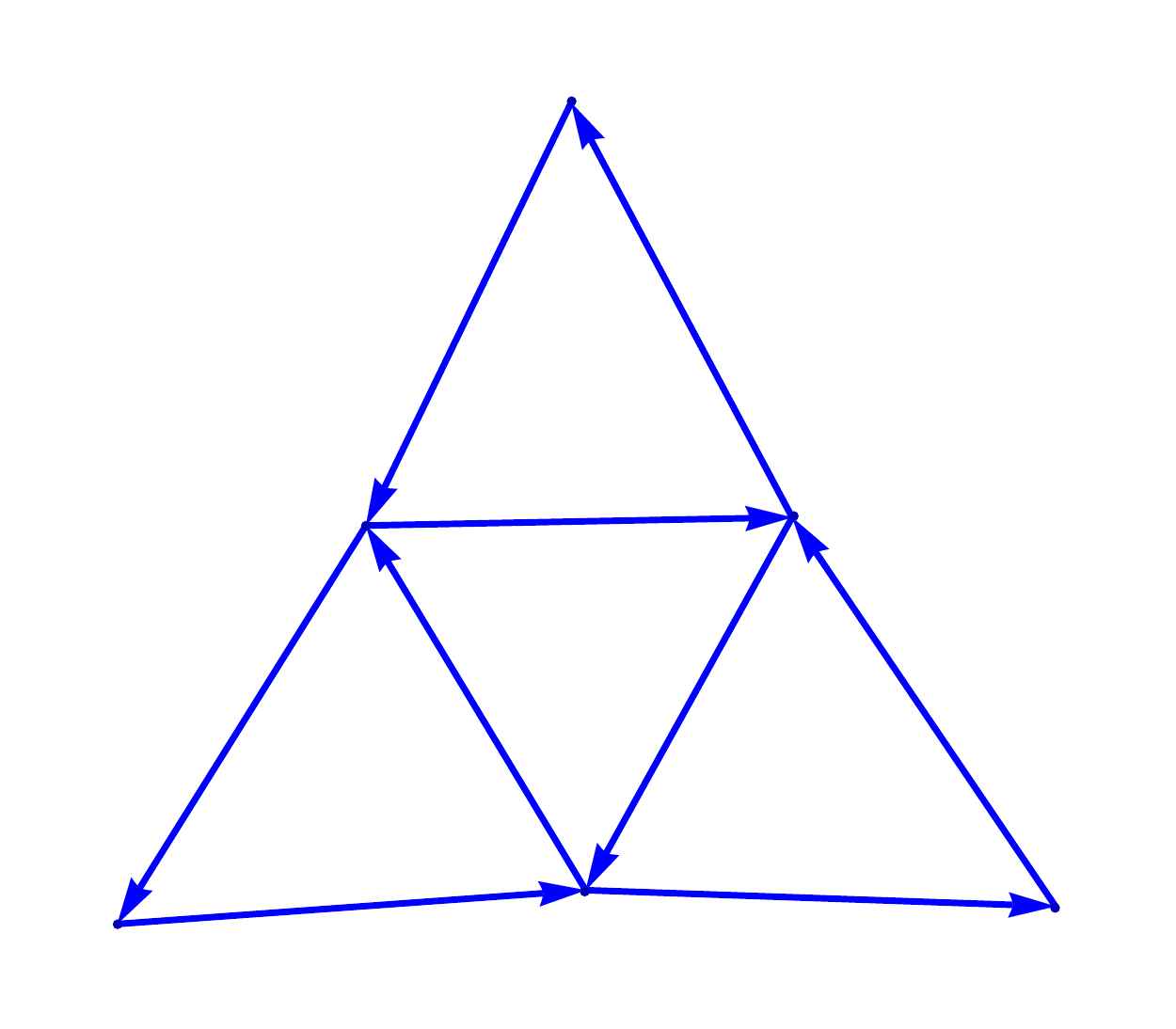}}

For $k=2$, $N=5$, $r_{{\mathcal A}_2} =6$,  the rank of the lattice $L_2 $ is  $\mathfrak{r} = 2\, r_{{\mathcal A}_2}  = 12$.\\
The period is a rhombus $5\times 5$ and  $\vert {\mathcal R}\vert = 100$.
A Gram matrix $A$ for the lattice is given below
\be
A = \left(
\begin{array}{cccccccccccc}
 6 & 0 & 2 & 0 & 2 & 0 & -2 & 1 & -2 & 2 & -2 & 2 \\
 0 & 6 & 2 & 2 & 2 & 2 & 1 & -1 & 0 & -2 & 0 & -2 \\
 2 & 2 & 6 & 0 & 2 & 2 & 2 & 2 & -1 & 1 & 2 & 2 \\
 0 & 2 & 0 & 6 & 2 & 0 & 0 & 2 & 1 & -2 & 2 & 0 \\
 2 & 2 & 2 & 2 & 6 & 0 & 2 & 2 & 2 & 2 & -1 & 1 \\
 0 & 2 & 2 & 0 & 0 & 6 & 0 & 2 & 2 & 0 & 1 & -2 \\
 -2 & 1 & 2 & 0 & 2 & 0 & 6 & 0 & 2 & 0 & 2 & 0 \\
 1 & -1 & 2 & 2 & 2 & 2 & 0 & 6 & 2 & 2 & 2 & 2 \\
 -2 & 0 & -1 & 1 & 2 & 2 & 2 & 2 & 6 & 0 & 0 & -2 \\
 2 & -2 & 1 & -2 & 2 & 0 & 0 & 2 & 0 & 6 & -2 & 2 \\
 -2 & 0 & 2 & 2 & -1 & 1 & 2 & 2 & 0 & -2 & 6 & 0 \\
 2 & -2 & 2 & 0 & 1 & -2 & 0 & 2 & -2 & 2 & 0 & 6 \\
\end{array}
\right)
\label{GramA2}
\ee
The discriminant is readily calculated:  $\Delta = 5^9$. The modular level is $\ell = N^2 = 25$. \\

\subsubsection{Theta function}

Applying the Hecke-Schoenberg theorem leads to the following result:
the theta function of this lattice of hyper-roots of type $\SU{3}$ at conformal level  $k=2$ is of weight $6$, modular level $\ell = 5^2 = 25$ (the square of the altitude) and Dirichlet character $\chi(11)$ for the characters modulo 25 on a cyclotomic field of order 20. It is the only character (namely the Kronecker character 5), the eleventh on a collection of $20 = \Phi_{Euler}(25)$) that coincides with the value of the Legendre symbol ${\mathfrak L}(\Delta,p)$ for all odd primes $p$ that do not divide $25$. This space of modular forms has dimension $16$. The theta function, in the variable $q_2=q^2$, is therefore fully determined by its $16$ first Fourier coefficients (the first being $1$). The coefficients of $q_2^a$ with $a>15$ are then predicted. The series starts as $\theta(z) = 1 + 100\,q_2^{3 } + 450\,q_2^{4 } + 960\,q_2^{5 } + 2800\,q_2^{6 } + 6600\,q_2^{7 } + 12300\,q_2^{8} + \ldots$. Here are the first coefficients, up to $q_2^{48}=q^{96}$: 
\\
{1, 100, 450, 960, 2800, 6600, 12300, 22400, 30690, 63000, 93150, \
144000, 203100, 236080, 392850, 550800, 708350, 961800, 972780, \
1581600, 1937250, 2495400, 2977400, 3063360, 4469400, 5547700, \
6477600, 7963200, 7344920, 11094000, 12627000, 15127200, 17091900, \
16459440, 22670850, 26899200, 29779950, 34869600, 31131750, 44964000, \
48927900, 57061200, 62034900, 57598720, 77425500, 89018400, 95469650,$\ldots$ }
Here is the Magma code leading to this result:
\begin{verbatim}
H := DirichletGroup(25,CyclotomicField(EulerPhi(25)));
chars := Elements(H); eps := chars[11];
M := ModularForms([eps],6); order:=48;
PowerSeries(M![1, 0, 0,100, 450, 960, 2800, 6600, 12300, 22400, 30690, 63000, 
93150, 144000, 203100, 236080],order);
\end{verbatim}
The first Fourier coefficients have to be computed by a brute force approach that relies, ultimately, on the explicitly given Gram matrix.

\subsubsection{Other properties of this lattice}
\paragraph{The automorphism group of the lattice $L_2$.}
The automorphism group $aut$ of $L_2$ is of order $1200$ and its structure, in terms of direct and semi-direct products, is $C_2 \times ((((C_5 \times C_5) \rtimes C_4) \rtimes C_3) \rtimes C_2)$.
Orbits of the basis vectors under the $aut$ action coincide and contain the $100$ hyper-roots (the $50$ positive and the $50$ negative ones). 
Since all the hyper-roots belong to a single orbit of $aut$, all their stabilizers are conjugated in $aut$, and found to be isomorphic with the group $D_{12}$ (which is itself isomorphic with $S_3 \times C_2$).
\paragraph{A short description of the Voronoi cells.}
We  only mention that the Voronoi polytope has $5410$ $11$-dimensional facets.

\subsection{The hyper-root lattice $L_3$ of  ${\mathcal A}_3 (\SU{3})$\hspace{5.cm}\protect\icon{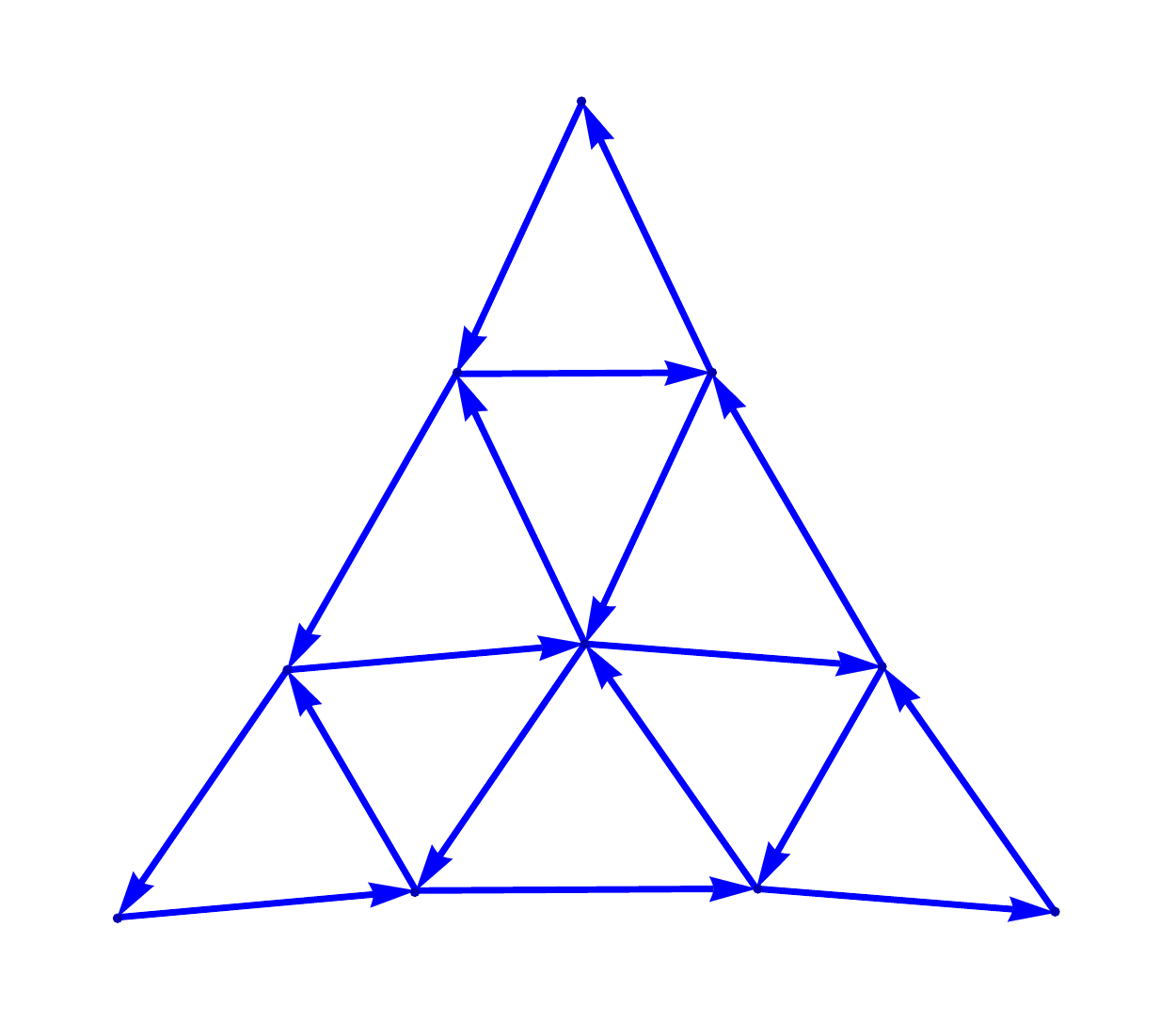}}
A Gram matrix is given in appendix \ref{AppendixGram}.
For $k=3$, $N=6$, $r_{{\mathcal A}_3} =10$, the rank of the lattice $L_3 $ is   $\mathfrak{r} = 2\, r_{{\mathcal A}_3} = 20$.
The period is a rhombus $6\times 6$ and  $\vert {\mathcal R}\vert = 240$.

{\bf Theta function}

The discriminant is readily calculated:  $\Delta = 6^{12}$. The modular level is $\ell = 18$. \\
The theta function belongs to  a space of modular forms on $\Gamma_0(18)$, of weight $10$, twisted by an appropriate character of modulus $18$ on a cyclotomic field of order $6 = \Phi_{Euler}(18)$.
The corresponding space of modular forms has dimension $31$ and the theta function of the lattice, determined by its first Fourier coefficients starts as 
\begin{eqnarray*}
 \theta(z) &=& 1+240 \, q^6+1782 \, q^8+9072 \, q^{10}+59328 \, q^{12}+216432 \, q^{14}+810000 \, q^{16}+2059152 \, q^{18}+ \\ 
 {} &{} &6080832 \, q^{20}+ 12349584 \, q^{22} + 31045596 \, q^{24} +  O\left(q^{25}\right)
\end{eqnarray*}   

Here is the list of its coefficients, up to order $60$ in $q$:
\footnotesize
\begin{equation*}
\begin{split}
& 1,0,0,240,1782,9072,59328,216432,810000,2059152,6080832,12349584,31045596,57036960,\\ &
122715648,204193872,418822650,622067040,1193611392,1734272208,3043596384,4217152080,\\ &
7354100160,9446435136,15901091892,20507712192,32268036096,40493364288,64454759856,\\ &
76079125584,118436670720,142127536464,209154411792,246451249296,369868125312,\\ &
413358056928,611268619740,698624989632,981886883328,1108342458624,1597262339340,\\ &
1716946287264,2447106074496,2701744008624,3674391470784,4018040848656,5617678157568,\\ &
5869298618208,8140982862948,8753718885120,11607623460864,12394567905984,16938128525364,\\ &
17305593381648,23493640620096,24756714700128,32196165379200,33726641096496,45246801175488,45433065648240
\end{split}
   \end{equation*}
   \normalsize
   
\bigskip

\paragraph{The automorphism group of the lattice $L_3$.}
The automorphism group $aut$ of $L_3$ is of order $864=2^5 3^3$. Its structure is $C_2 \times ((((C_6 \times C_6) \rtimes C_3) \rtimes C_2) \rtimes C_2)$.
\paragraph{Voronoi cells.}
The Voronoi polytope has $539214$ $19$-dimensional facets. 

\paragraph{Note.}
 {\sl In what follows we shall only provide basic information about the lattices. 
All lattice related properties ultimately rely on the explicit expression of Gram matrices ---that will be displayed in the coming sections or in appendix \ref{AppendixGram}.
In particular, for most following examples, we only mention the order of the automorphism group, as determined by the computer algebra system Magma, without discussing its structure or the way it is generated.}

\subsection{The hyper-root lattice $L_4$ of  ${\mathcal A}_4(\SU{3})$\hspace{5.cm}\protect\icon{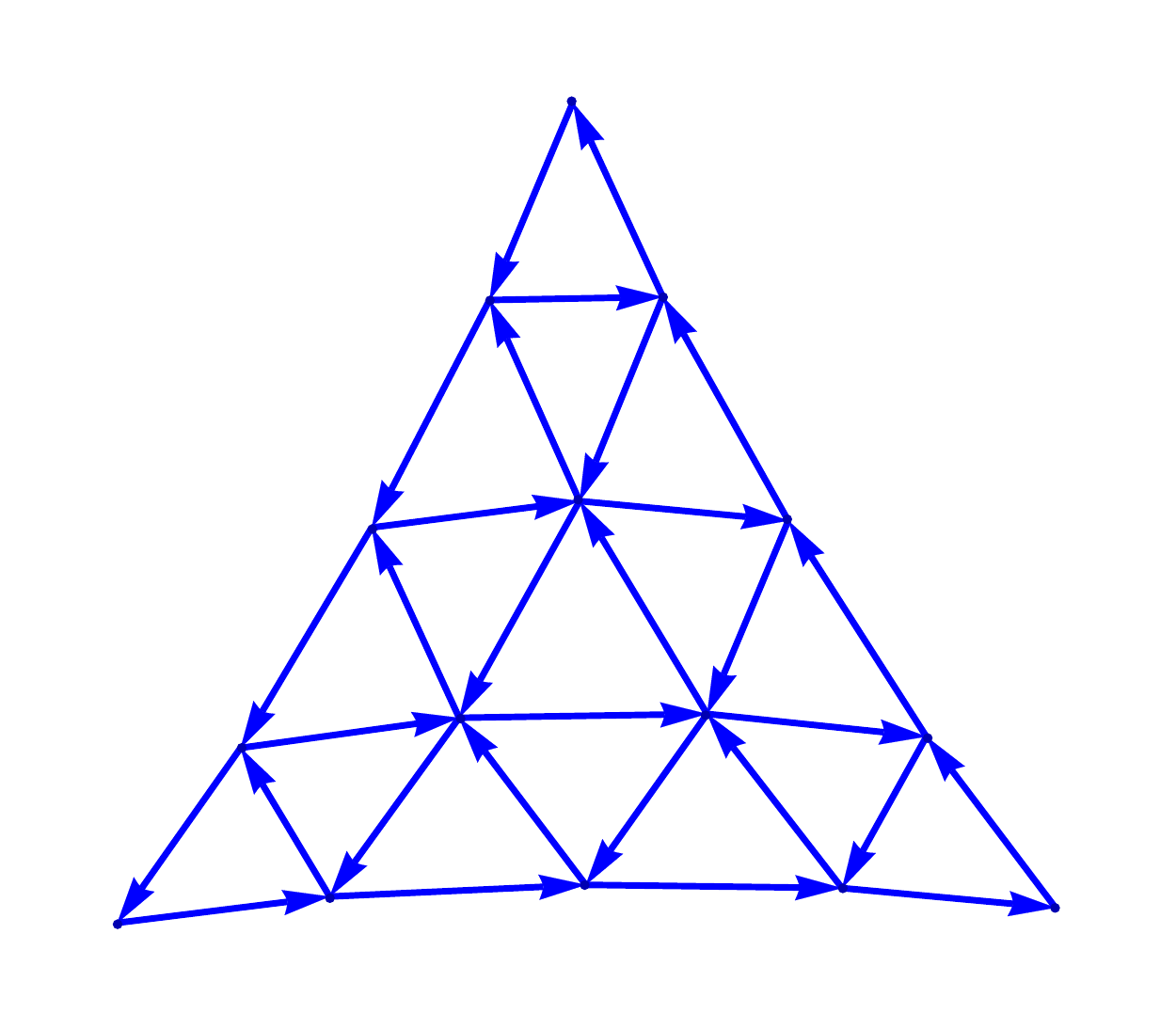}}

For $k=4$, $N=7$, $r_{{\mathcal A}_4} =15$, the rank of the lattice $L_4 $ is  $\mathfrak{r} = 2\, r_{{\mathcal A}_4} = 30$.
The period is a rhombus $7\times 7$ and  $\vert {\mathcal R}\vert = 490 $.
The discriminant is $7^{15}$, and the level is 49.  The automorphism group has order  $2^3 \, 3^2 \, 7^2 = 3528$.
 The dimension of the appropriate space of modular forms (modular forms on $\Gamma_1(49)$ with character Kronecker character $-7$ and
weight 15), is quite large: it has dimension 70 over the ring of integers.  The theta series starts as follows:
\begin{eqnarray*}
 \theta(z) &=& 1 + 490 \, q^6 + 4998 \, q^8 + 45864 \, q^{10} + 464422 \, q^{12} + 3429426 \, q^{14} + 21668094 \, q^{16} + 111678742 \, q^{18} + \\ 
  {} &{} &492567012 \, q^{20} + 1876801038 \, q^{22} + 6352945942 \, q^{24} + 19484903508 \, q^{26} + 54935857326 \, q^{28} + \\
  {} & {} & 144330551050 \, q^{30} + O(q^{31})
\end{eqnarray*}  

\subsection{Higher $L_k$'s}

We only mention that the theta series of $L_5$ (of rank $42$) and $L_6$ (of rank $56$) start as follows:
$${\mathcal A}_5: \quad \theta(z)  = 1+896 \, q^6+11856 \, q^8+154368 \, q^{10}+2331648 \, q^{12}+ 27065088 \, q^{14}  + 281311128 \, q^{16} +O\left(q^{17}\right)$$
 $${\mathcal A}_6: \quad \theta(z) = 1+1512 \, q^6+24300 \, q^8+425736 \, q^{10}+8530758 \, q^{12}+O\left(q^{13}\right)$$

\subsection{Some remarks about the lattices $L_k$}
\label{simpleGram}
\paragraph{About the determination of a Gram matrix, for general $L_k$.}
We already described in section~\ref{choicebasis} one way to select a basis made of hyper-roots. However, for each particular value of $k$, the determination of a Gram matrix, using equation~\ref{scalarproductofhyperroots1}
and a chosen basis, is not a computationally totally trivial task, and it would be nice to have a way to deduce such a matrix from the fusion coefficients of the module by a simpler algorithm.
As already commented in section~\ref{choicebasis}, we do not have any canonical choice here for the Gram matrix (no available Cartan matrix) and the naive generalization of the $\SU{2}$ algorithm to the $\SU{3}$ family fails. 
Let us nevertheless mention that in the cases $k = 1$ (see matrix $A^{\prime\prime}$ in section~\ref{GramA1}) and $k=2$, the following simple expressions, written in terms of fusion matrices, are Gram matrices for the lattices $L_1$, $L_2$ and are equivalent to those given previously:
{\small
$$
6 \one_{6} \,  + \, 
\left(
\begin{array}{cc}
2 (F_{\{1,2\}}+F_{\{2,1\}}) & 2 \one_3 -2
   (F_{\{1,2\}}+F_{\{2,1\}}) \\
 2 \one_3 -2 (F_{\{1,2\}}+F_{\{2,1\}}) & -2
   (F_{\{1,2\}}+F_{\{2,1\}}) \\
\end{array}
\right)
$$
$$
6 \one_{12} \, + \,
\left(
\begin{array}{cc}
 2 (F_{\{1,2\}}+F_{\{2,1\}}) &
 2  \one_6+(F_{\{1,2\}}+F_{\{2,1\}})+(F_{\{1,3\}}+F_{\{3,1\}})-F_{\{2,2\}} \\
 2 \one_6+(F_{\{1,2\}}+F_{\{2,1\}})+(F_{\{1,3\}}+F_{\{3,1\}})-F_{\{2,2\}} & 
 2 (F_{\{1,2\}}+F_{\{2,1\}}) \\
\end{array}
\right)
$$
}

\paragraph{About the determination of $\theta(z)$, for general $k$.}
The theta function of $L_k$, as a modular form twisted by a character, can, in principle, be obtained by following the method explained in the previous sections and illustrated in the case of the first few members of the $L_k$ series.
In this respect we observed that the (quadratic form) level of $L_k$ is often equal to  $\ell = (k+3)^2$ but it is not so for $L_3$ where the level is $18$ and not $36$.
Notice that for $L_1$, the matrix $8 A^{-1}$ is integral but its diagonal elements are odd, so the level is indeed 16.\\
The discriminant of the lattice $L_k$ is $(k + 3)^{3 (k + 1)}$,
the weight is $r_{{\mathcal A}_k} = (k+1)(k+2)/2$, the quadratic form level is readily obtained from the Gram matrix, and the determination of the appropriate character requires a discussion relying on the arithmetic properties of the discriminant and of the level. 
However, the first coefficients of the Fourier series expansion have to be found, and the number of needed coefficients depends on the properties of an appropriate space of modular forms. The determination of the needed coefficients is done by brute force, namely by computing the norm of the vectors belonging to the first shells, using the Gram matrix as an input.  Moreover, the explicit determination of a Gram matrix for $L_k$ ($k$ being given) also becomes a non-trivial exercise when $k$ is large (see the previous comment).
The present method may therefore become rapidly intractable if we increase $k$ too much. 
Admittedly it would be nice to have a general formula, like the one that we have for the root lattices of type $A_{n-1}$, that would be valid for all $k$'s, and would express the theta function of $L_k$ in terms of known functions (for instance elliptic theta's). This was not done but we hope that our results will trigger new developments in that direction.

\paragraph{About the vectors of smallest norm.}  For all the lattices $L_k$ that we considered explicitly, 
 the lattice vectors of shortest length are precisely the hyper-roots ($100$ of them for $L_3$, for instance), the kissing number of those lattices are then given by the number of hyper-roots. 
 As it is well known, this property holds for all usual root lattices, \ie hyper-root lattices of the $\SU{2}$ family.
However, as we shall see below, this property does not always hold for those lattices associated with modules  of the $\SU{3}$ family that are not of type ${\mathcal A}_k(\SU{3})$.

\subsection{Theta function for ${\mathcal D}_3(\SU{3})$\hspace{5.cm}\protect\icon{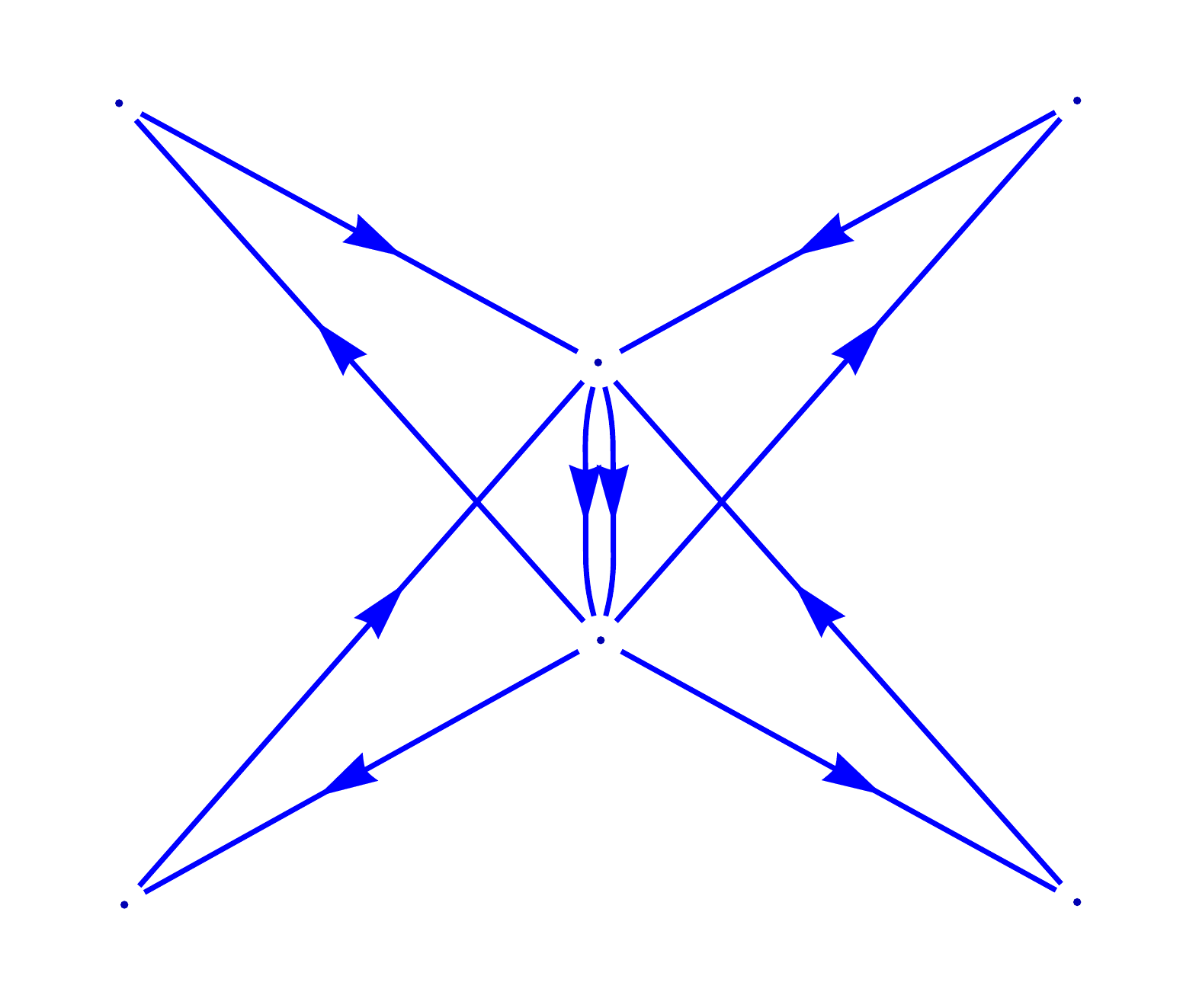}}
A Gram matrix is given in appendix \ref{AppendixGram}.
The fusion graphs of the ${\mathcal D}_k(\SU{3})$ series are $\ZZ_3$ orbifolds of the ${\mathcal A}_k(\SU{3})$, and, when $k=0\; \text{mod}\;  3$, their number of vertices (simple objects of the category) is $\tfrac{1}{3}(r_{{\mathcal A}_k} -1) + 3$, \ie $ \tfrac{1}{3} (\tfrac{(k+1)(k+2)}{2} -1)+3$. So, for $k=3$ we have $r_{ {\mathcal D}_3}=6$, and the rank of the quadratic form is $\mathfrak{r} = 2\, r_{ {\mathcal D}_3} = 12$. 
The discriminant of the quadratic form is $3^{12}$ and the (modular) level is $2\times 9$. The order of the automorphism group is $2^8 3^3$.
The reader can check that $\theta$, given below, belongs to a  space of modular forms on $\Gamma_0(18)$, of weight $6$, twisted by an appropriate character.

\small
 $$ \theta(z) = 1+36 \, q^4+144 \, q^6+486 \, q^8+2880 \, q^{10}+5724 \, q^{12}+
  7776 \, q^{14} + 31068\, q^{16} + 
    40320 \, q^{18} + 47628 \, q^{20} +
 O\left(q^{21}\right)$$
 \normalsize

The number of hyper-roots is $\vert {\mathcal R} \vert = 2 N^2  r_{{\mathcal D}_3} /3 = 2(3+3)^2 6/3 = 144$, whereas the number of vectors of smallest norm is $36$.
This is the first manifestation of a phenomenon that we mentioned in the previous paragraph and that never occurs for usual root lattices. In the present case, the first shell is made of vectors of norm $4$, that are not hyper-roots, and the only vectors of the lattice that belong to the second shell, or norm $6$, are precisely the hyper-roots.  Vectors of smallest norm can of course be expanded on a chosen basis of hyper-roots; here are, for instance, the components of one of them, on the basis that is chosen to write the Gram matrix $A$ of ${\mathcal D}_3$ in appendix~\ref{AppendixGram}~: taking $v = \{1, 1, 1, 1, -2, -1, 0, 0, 0, 0, 1, 1\}$, one can check that $<v,v>=4$. 
One finds that the vector space spanned by the $36$ vectors of shortest length (it is enough to choose $18$ of them from the pairs $(v, -v)$) is of dimension $6$.

\subsection{Theta function for ${\mathcal D}_6(\SU{3})$\hspace{5.cm}\protect\icon{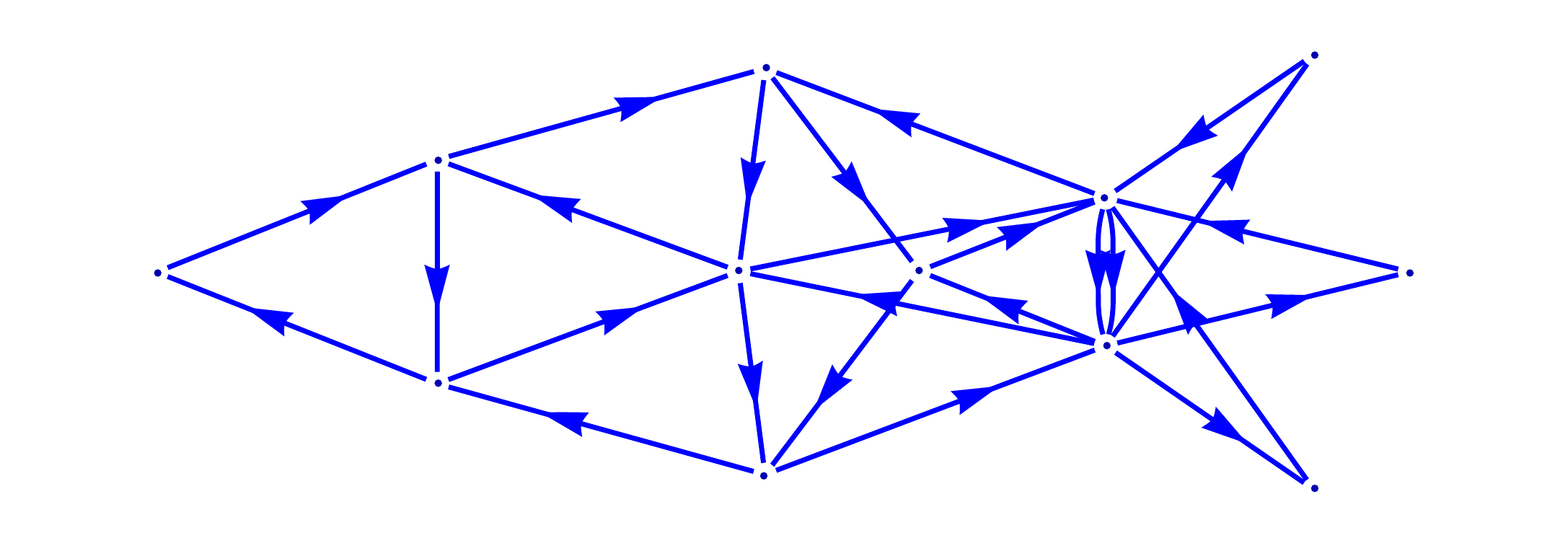}}
A Gram matrix is given in appendix \ref{AppendixGram}.
The number of simple objects is  $r_{{\mathcal D}_6}=12$ and the rank of the quadratic form is $\mathfrak{r} = 2\, r_{ {\mathcal D}_6} = 24$.
The discriminant of the quadratic form is $3^{18}$ and the modular level is $2\times 27$. The order of the automorphism group is $2^6 3^{11}$. The theta function reads:
{\small 
\begin{equation*}
\begin{split}
 \theta(z) &= 1+162 \, q^4+2322 \, q^6+35478 \, q^8+273942 \, q^{10}+1771326 \, q^{12}+9680148 \, q^{14}+40813632
   \, q^{16}+ \\
   &150043014 \, q^{18}+484705782 \, q^{20}+O\left(q^{21}\right)
   \end{split}
   \end{equation*}}
The number of hyper-roots is $\vert {\mathcal R} \vert = 2 N^2 r_{{\mathcal D}_6} /3 = 2(6+3)^2 12/3 = 648$.
The first shell is made of $162$ vectors of norm $4$, that are not hyper-roots, and the second shell, of norm $6$, contains not only the hyper-roots themselves, but $3522-648=2874$ other vectors.

\ommit{
\subsection{Theta function for ${\mathcal D}_4^\star (\SU{3})$ \hspace{5.cm}\protect\icon{adjgraphD4star.pdf}}
The rank is $\mathfrak{r} = 2\times 9 = 18$.
\blue This one is special (no triality). What ! Check this ! Results to be checked again. Or Remove it (also in the text). \normalcolor
 $$\theta(z) =1+336 \, q^6+4410 \, q^8+22176 \, q^{10}+100464 \, q^{12}+346752 \, q^{14}+993132 \, q^{16}+2562896
   \, q^{18}+5974416 \, q^{20}+O\left(q^{21}\right)$$
   \blue but R is 294 \normalcolor
   \\
The order of the automorphism group is $2^9 3^2 7^2$.}

\subsection{Theta function for ${\mathcal E}_5 (\SU{3})$\hspace{5.cm}\protect\icon{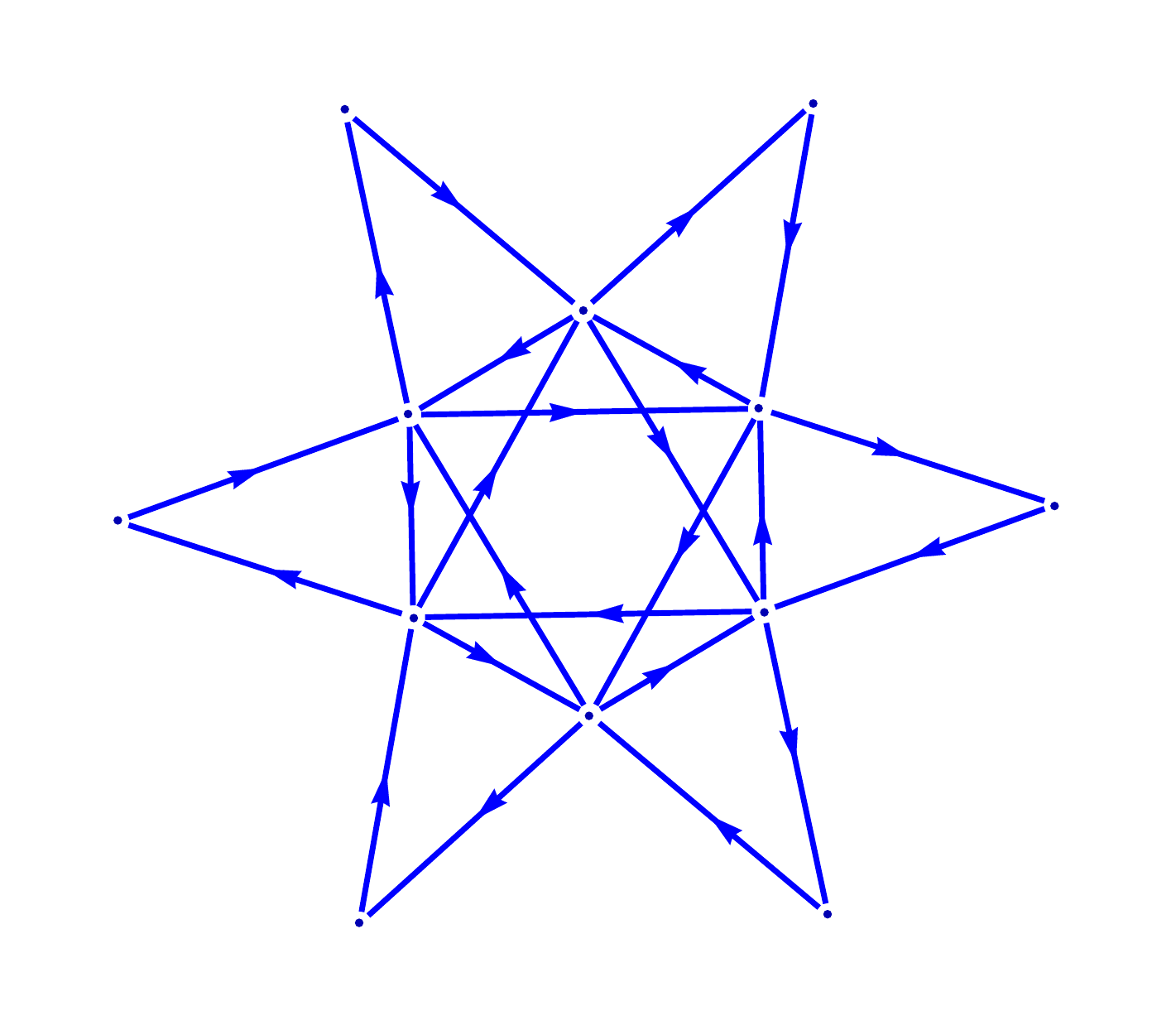}}
A Gram matrix is given in appendix \ref{AppendixGram}.
The rank of the category is $r_{{\mathcal E}_5} = 12$ and the rank of the quadratic form is $\mathfrak{r} = 2\, r_{{\mathcal E}_5} =24$.
Its discriminant is $2^{30}$ and the modular level is $2\times 8$.  The order of the automorphism group is $2^{11} 3$. The theta function reads:
{\small 
\begin{equation*}
\begin{split}
\theta(z) &= 1+512 \, q^6+11232 \, q^8+145920 \, q^{10}+1055616 \, q^{12}+5618688 \, q^{14}+25330128 \, q^{16}+ \\
& 89127936 \, q^{18}+295067136 \, q^{20}+O\left(q^{21}\right)
   \end{split}
   \end{equation*}}
Here, the hyper-roots ($2 N^2  r_{{\mathcal E}_5}  /3 = 2 8^2 12/3 = 512$ of them), like for the  ${\mathcal A}_k$ series, coincide with the vectors of smallest length.

\subsection{Theta function for ${\mathcal E}_9(\SU{3})$ \hspace{5.cm}\protect\icon{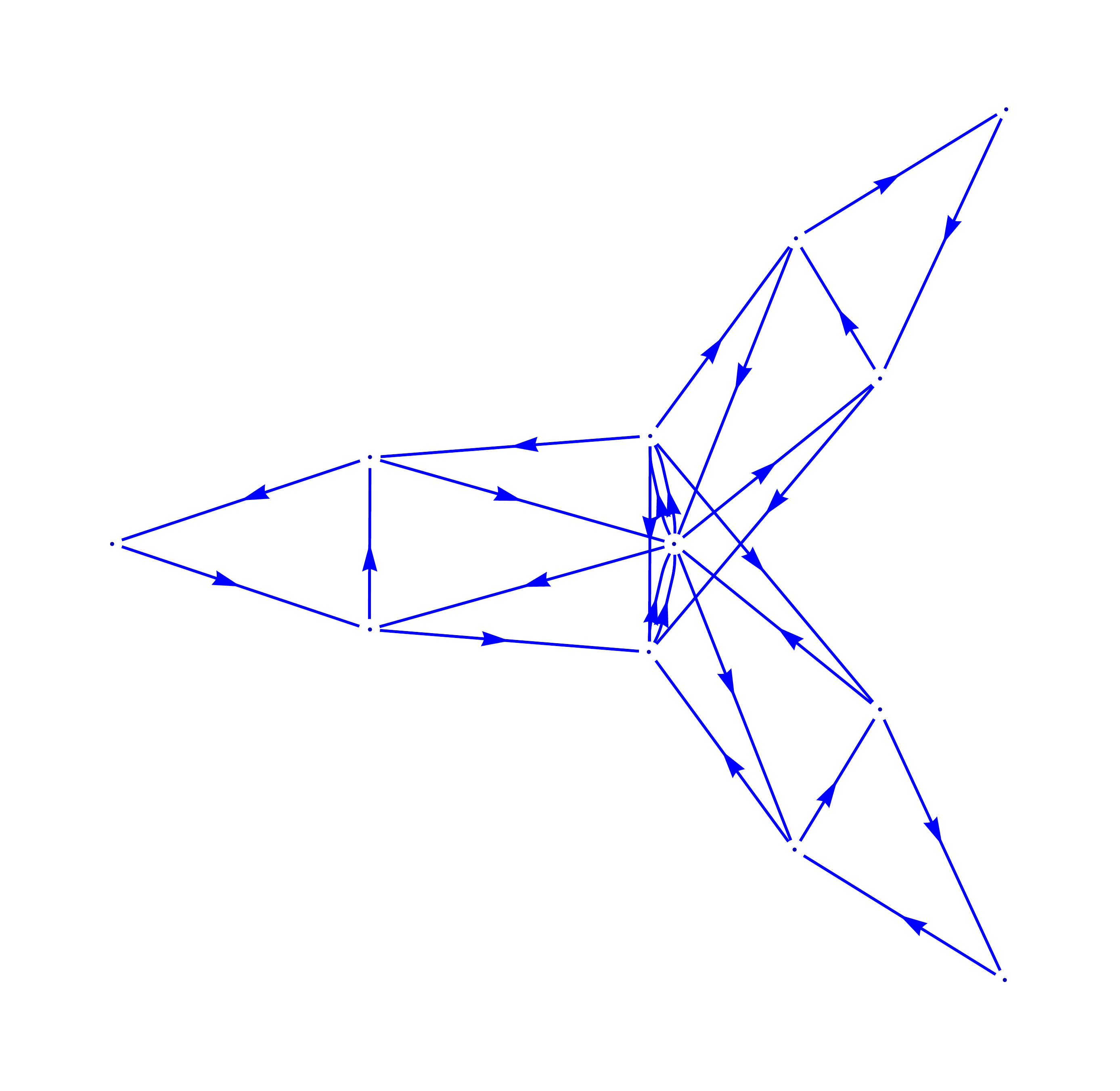}}
A Gram matrix is given in appendix \ref{AppendixGram}.
The rank of the category is $ r_{{\mathcal E}_9}  = 12$ and the rank of the quadratic form is $\mathfrak{r} = 2\,  r_{{\mathcal E}_9} =24$.
Its discriminant is $2^{24}$ and the modular level is $2\times 8$. One finds:
{\small 
\begin{equation*}
\begin{split}
\theta(z) &= 
1 + 756\, q^4 + 5760\, q^6 + 98928\, q^8 + 1092096\, q^{10} + 8435760\, q^{12} + 45142272\, q^{14} + 202712400 \, q^{16} +\\
& 715373568 \, q^{18} + 2350118808\, q^{20} + O\left(q^{21}\right)
   \end{split}
   \end{equation*}}

The number of hyper-roots is $\vert {\mathcal R} \vert = 2 (9+3)^2 12 / 3 = 1152$ and we observe that there are $756$ vectors of smaller norm ($4$) that build the first shell, and that the second shell contains not only the hyper-roots, but other vectors as well.

\ommit{
> ThetaSeriesModularFormSpace (E9);
Space of modular forms on Gamma_0(8) of weight 12 and dimension 13 over Integer Ring.
}

\subsection{Theta function for ${\mathcal E}_{21}(\SU{3})$ \hspace{5.cm}\protect\icon{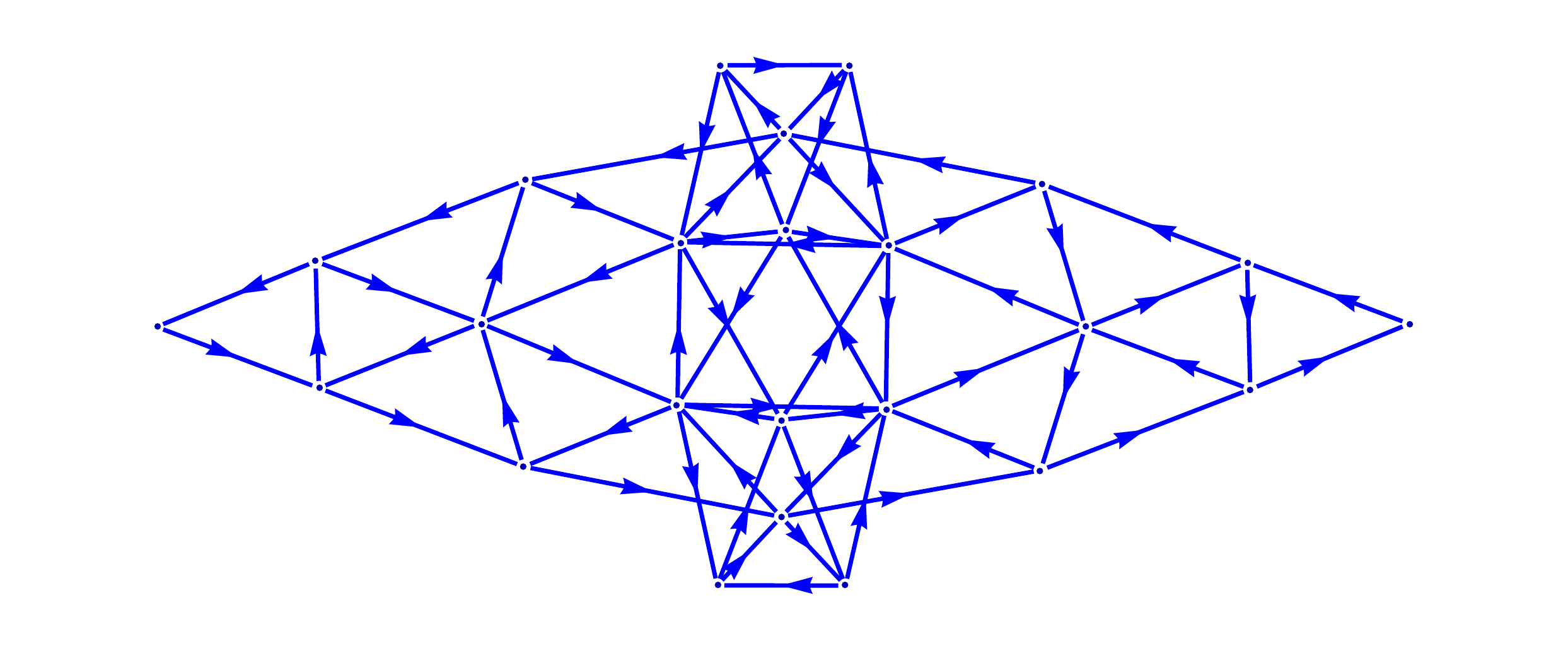}}
A Gram matrix is given in appendix \ref{AppendixGram}.
Here we have $r_{{\mathcal E}_{21}}  = 24$, so that the rank of the quadratic form is  $\mathfrak{r} = 2\, r_{{\mathcal E}_{21}}=48$.
Its discriminant is $3^{12}$ and the modular level is $2\times 3$. One finds:
{\small 
$$ \theta(z) =1+144 \, q^4+64512 \, q^6+ 54181224\, q^8 + O\left(q^9\right)
$$
}

The number of hyper-roots is $\vert {\mathcal R} \vert = 2 (21+3)^2 24 / 3 = 9216$ but the kissing number is only $144$. 
Therefore, here again, the vectors of smallest norm are not hyper-roots, and not all the lattice vectors of the second shell are hyper-roots.

\section{Appendices}
\subsection{Fundamental fusion matrices}
\label{AppendixFusionMatrices}
For completeness sake we give the fundamental fusion matrices $F_{(1,0)} = F_{\{2,1\}}$ for the cases considered in the text (we have replaced zeroes by dots). 
They are also adjacency matrices for the associated fusion graphs.  
These expressions are needed to determine, first, the other fundamental fusion matrices, using the recurrence relation~\ref{recursion}, then the inner product $<\alpha, \beta>$ of hyper-roots, using equations~\ref{scalarproductofhyperroots1} or \ref{scalarproductofhyperroots2}.
Similar expressions for the other $\SU{3}$ cases can be gathered from the available literature, and also from the website \cite{RCsiteWebFusionGraphs}.
{\scriptsize
$$ {\mathcal A}_1 : \quad \left(
\begin{array}{ccc}
 .  & 1 & .  \\
 .  & .  & 1 \\
 1 & .  & .  \\
\end{array}
\right)
\qquad 
{\mathcal A}_2 : \quad \left(
\begin{array}{cccccc}
 .  & 1 & .  & .  & .  & .  \\
 .  & .  & 1 & 1 & .  & .  \\
 1 & .  & .  & .  & 1 & .  \\
 .  & .  & .  & .  & 1 & .  \\
 .  & 1 & .  & .  & .  & 1 \\
 .  & .  & 1 & .  & .  & .  \\
\end{array}
\right)
\qquad
 {\mathcal A}_3 : \quad 
\left(
\begin{array}{cccccccccc}
 .  & 1 & .  & .  & .  & .  & .  & .  & .  & .  \\
 .  & .  & 1 & 1 & .  & .  & .  & .  & .  & .  \\
 1 & .  & .  & .  & 1 & .  & .  & .  & .  & .  \\
 .  & .  & .  & .  & 1 & .  & 1 & .  & .  & .  \\
 .  & 1 & .  & .  & .  & 1 & .  & 1 & .  & .  \\
 .  & .  & 1 & .  & .  & .  & .  & .  & 1 & .  \\
 .  & .  & .  & .  & .  & .  & .  & 1 & .  & .  \\
 .  & .  & .  & 1 & .  & .  & .  & .  & 1 & .  \\
 .  & .  & .  & .  & 1 & .  & .  & .  & .  & 1 \\
 .  & .  & .  & .  & .  & 1 & .  & .  & .  & .  \\
\end{array}
\right)$$
}
{\scriptsize
$$
 {\mathcal A}_4 : \quad
\left(
\begin{array}{ccccccccccccccc}
 .  & 1 & .  & .  & .  & .  & .  & .  & .  & .  & .  & .  & .  & .  & .  \\
 .  & .  & 1 & 1 & .  & .  & .  & .  & .  & .  & .  & .  & .  & .  & .  \\
 1 & .  & .  & .  & 1 & .  & .  & .  & .  & .  & .  & .  & .  & .  & .  \\
 .  & .  & .  & .  & 1 & .  & 1 & .  & .  & .  & .  & .  & .  & .  & .  \\
 .  & 1 & .  & .  & .  & 1 & .  & 1 & .  & .  & .  & .  & .  & .  & .  \\
 .  & .  & 1 & .  & .  & .  & .  & .  & 1 & .  & .  & .  & .  & .  & .  \\
 .  & .  & .  & .  & .  & .  & .  & 1 & .  & .  & 1 & .  & .  & .  & .  \\
 .  & .  & .  & 1 & .  & .  & .  & .  & 1 & .  & .  & 1 & .  & .  & .  \\
 .  & .  & .  & .  & 1 & .  & .  & .  & .  & 1 & .  & .  & 1 & .  & .  \\
 .  & .  & .  & .  & .  & 1 & .  & .  & .  & .  & .  & .  & .  & 1 & .  \\
 .  & .  & .  & .  & .  & .  & .  & .  & .  & .  & .  & 1 & .  & .  & .  \\
 .  & .  & .  & .  & .  & .  & 1 & .  & .  & .  & .  & .  & 1 & .  & .  \\
 .  & .  & .  & .  & .  & .  & .  & 1 & .  & .  & .  & .  & .  & 1 & .  \\
 .  & .  & .  & .  & .  & .  & .  & .  & 1 & .  & .  & .  & .  & .  & 1 \\
 .  & .  & .  & .  & .  & .  & .  & .  & .  & 1 & .  & .  & .  & .  & .  \\
\end{array}
\right)
\qquad
 {\mathcal D}_6 : \quad 
\left(
\begin{array}{cccccccccccc}
 .  & .  & .  & .  & .  & .  & .  & .  & 1 & .  & .  & .  \\
 .  & .  & .  & .  & .  & .  & .  & .  & 1 & .  & .  & .  \\
 .  & .  & .  & .  & .  & .  & .  & 1 & .  & .  & .  & .  \\
 .  & .  & .  & .  & .  & .  & 1 & 1 & 1 & .  & .  & .  \\
 .  & .  & .  & .  & .  & .  & .  & .  & 1 & .  & .  & .  \\
 .  & .  & .  & .  & .  & .  & 1 & .  & 1 & .  & .  & .  \\
 .  & .  & .  & .  & .  & .  & .  & .  & .  & 1 & 1 & .  \\
 .  & .  & .  & .  & .  & .  & .  & .  & .  & 1 & .  & 1 \\
 .  & .  & .  & .  & .  & .  & .  & .  & .  & .  & 2 & 1 \\
 .  & .  & 1 & 1 & .  & .  & .  & .  & .  & .  & .  & .  \\
 1 & 1 & .  & 1 & 1 & 1 & .  & .  & .  & .  & .  & .  \\
 .  & .  & .  & 1 & .  & 1 & .  & .  & .  & .  & .  & .  \\
\end{array}
\right)$$
}
{\scriptsize
$$
{\mathcal A}_5 : \quad \left(
\begin{array}{ccccccccccccccccccccc}
 .  & 1 & .  & .  & .  & .  & .  & .  & .  & .  & .  & .  & .  & .  & .  & .  & .  & .  & .  & .  & .  \\
 .  & .  & 1 & 1 & .  & .  & .  & .  & .  & .  & .  & .  & .  & .  & .  & .  & .  & .  & .  & .  & .  \\
 1 & .  & .  & .  & 1 & .  & .  & .  & .  & .  & .  & .  & .  & .  & .  & .  & .  & .  & .  & .  & .  \\
 .  & .  & .  & .  & 1 & .  & 1 & .  & .  & .  & .  & .  & .  & .  & .  & .  & .  & .  & .  & .  & .  \\
 .  & 1 & .  & .  & .  & 1 & .  & 1 & .  & .  & .  & .  & .  & .  & .  & .  & .  & .  & .  & .  & .  \\
 .  & .  & 1 & .  & .  & .  & .  & .  & 1 & .  & .  & .  & .  & .  & .  & .  & .  & .  & .  & .  & .  \\
 .  & .  & .  & .  & .  & .  & .  & 1 & .  & .  & 1 & .  & .  & .  & .  & .  & .  & .  & .  & .  & .  \\
 .  & .  & .  & 1 & .  & .  & .  & .  & 1 & .  & .  & 1 & .  & .  & .  & .  & .  & .  & .  & .  & .  \\
 .  & .  & .  & .  & 1 & .  & .  & .  & .  & 1 & .  & .  & 1 & .  & .  & .  & .  & .  & .  & .  & .  \\
 .  & .  & .  & .  & .  & 1 & .  & .  & .  & .  & .  & .  & .  & 1 & .  & .  & .  & .  & .  & .  & .  \\
 .  & .  & .  & .  & .  & .  & .  & .  & .  & .  & .  & 1 & .  & .  & .  & 1 & .  & .  & .  & .  & .  \\
 .  & .  & .  & .  & .  & .  & 1 & .  & .  & .  & .  & .  & 1 & .  & .  & .  & 1 & .  & .  & .  & .  \\
 .  & .  & .  & .  & .  & .  & .  & 1 & .  & .  & .  & .  & .  & 1 & .  & .  & .  & 1 & .  & .  & .  \\
 .  & .  & .  & .  & .  & .  & .  & .  & 1 & .  & .  & .  & .  & .  & 1 & .  & .  & .  & 1 & .  & .  \\
 .  & .  & .  & .  & .  & .  & .  & .  & .  & 1 & .  & .  & .  & .  & .  & .  & .  & .  & .  & 1 & .  \\
 .  & .  & .  & .  & .  & .  & .  & .  & .  & .  & .  & .  & .  & .  & .  & .  & 1 & .  & .  & .  & .  \\
 .  & .  & .  & .  & .  & .  & .  & .  & .  & .  & 1 & .  & .  & .  & .  & .  & .  & 1 & .  & .  & .  \\
 .  & .  & .  & .  & .  & .  & .  & .  & .  & .  & .  & 1 & .  & .  & .  & .  & .  & .  & 1 & .  & .  \\
 .  & .  & .  & .  & .  & .  & .  & .  & .  & .  & .  & .  & 1 & .  & .  & .  & .  & .  & .  & 1 & .  \\
 .  & .  & .  & .  & .  & .  & .  & .  & .  & .  & .  & .  & .  & 1 & .  & .  & .  & .  & .  & .  & 1 \\
 .  & .  & .  & .  & .  & .  & .  & .  & .  & .  & .  & .  & .  & .  & 1 & .  & .  & .  & .  & .  & .  \\
\end{array}
\right)
\qquad
 {\mathcal D}_3 : \quad \left(
\begin{array}{cccccc}
 .  & .  & .  & .  & 1 & .  \\
 .  & .  & .  & .  & 1 & .  \\
 .  & .  & .  & .  & 1 & .  \\
 .  & .  & .  & .  & 1 & .  \\
 .  & .  & .  & .  & .  & 2 \\
 1 & 1 & 1 & 1 & .  & .  \\
\end{array}
\right)$$
}
{\scriptsize
$$ {\mathcal E}_5 : \quad \left(
\begin{array}{cccccccccccc}
 .  & .  & .  & .  & .  & .  & 1 & .  & .  & .  & .  & .  \\
 .  & .  & .  & .  & .  & .  & .  & 1 & .  & .  & .  & .  \\
 .  & .  & .  & .  & 1 & .  & 1 & 1 & .  & .  & .  & .  \\
 .  & .  & .  & .  & .  & 1 & 1 & 1 & .  & .  & .  & .  \\
 .  & .  & .  & .  & .  & .  & .  & .  & .  & .  & .  & 1 \\
 .  & .  & .  & .  & .  & .  & .  & .  & .  & .  & 1 & .  \\
 .  & .  & .  & .  & .  & .  & .  & .  & 1 & .  & 1 & 1 \\
 .  & .  & .  & .  & .  & .  & .  & .  & .  & 1 & 1 & 1 \\
 .  & .  & 1 & .  & .  & .  & .  & .  & .  & .  & .  & .  \\
 .  & .  & .  & 1 & .  & .  & .  & .  & .  & .  & .  & .  \\
 1 & .  & 1 & 1 & .  & .  & .  & .  & .  & .  & .  & .  \\
 .  & 1 & 1 & 1 & .  & .  & .  & .  & .  & .  & .  & .  \\
\end{array}
\right)
\qquad
 {\mathcal E}_9 : \quad \left(
\begin{array}{cccccccccccc}
 .  & .  & .  & .  & 1 & .  & .  & .  & .  & .  & .  & .  \\
 .  & .  & .  & .  & .  & 1 & .  & .  & .  & .  & .  & .  \\
 .  & .  & .  & .  & .  & .  & 1 & .  & .  & .  & .  & .  \\
 .  & .  & .  & .  & 1 & 1 & 1 & 2 & .  & .  & .  & .  \\
 .  & .  & .  & .  & .  & .  & .  & .  & 1 & .  & .  & 1 \\
 .  & .  & .  & .  & .  & .  & .  & .  & .  & 1 & .  & 1 \\
 .  & .  & .  & .  & .  & .  & .  & .  & .  & .  & 1 & 1 \\
 .  & .  & .  & .  & .  & .  & .  & .  & 1 & 1 & 1 & 1 \\
 1 & .  & .  & 1 & .  & .  & .  & .  & .  & .  & .  & .  \\
 .  & 1 & .  & 1 & .  & .  & .  & .  & .  & .  & .  & .  \\
 .  & .  & 1 & 1 & .  & .  & .  & .  & .  & .  & .  & .  \\
 .  & .  & .  & 2 & .  & .  & .  & .  & .  & .  & .  & .  \\
\end{array}
\right)$$
}
{\scriptsize
$$ {\mathcal E}_{21} : \quad \left(
\begin{array}{cccccccccccccccccccccccc}
 .  & .  & .  & .  & .  & .  & .  & .  & 1 & .  & .  & .  & .  & .  & .  & .  & .  & .  & .  & .  & .  & .  &
   .  & .  \\
 .  & .  & .  & .  & .  & .  & .  & .  & 1 & 1 & 1 & .  & .  & .  & .  & .  & .  & .  & .  & .  & .  & .  &
   .  & .  \\
 .  & .  & .  & .  & .  & .  & .  & .  & .  & .  & 1 & .  & 1 & .  & 1 & .  & .  & .  & .  & .  & .  & .  &
   .  & .  \\
 .  & .  & .  & .  & .  & .  & .  & .  & .  & 1 & .  & 1 & .  & 1 & .  & .  & .  & .  & .  & .  & .  & .  &
   .  & .  \\
 .  & .  & .  & .  & .  & .  & .  & .  & .  & .  & 1 & .  & 1 & 1 & .  & .  & .  & .  & .  & .  & .  & .  &
   .  & .  \\
 .  & .  & .  & .  & .  & .  & .  & .  & .  & .  & 1 & 1 & .  & 1 & .  & .  & .  & .  & .  & .  & .  & .  &
   .  & .  \\
 .  & .  & .  & .  & .  & .  & .  & .  & .  & .  & .  & .  & .  & 1 & 1 & 1 & .  & .  & .  & .  & .  & .  &
   .  & .  \\
 .  & .  & .  & .  & .  & .  & .  & .  & .  & .  & .  & .  & .  & .  & .  & 1 & .  & .  & .  & .  & .  & .  &
   .  & .  \\
 .  & .  & .  & .  & .  & .  & .  & .  & .  & .  & .  & .  & .  & .  & .  & .  & 1 & 1 & .  & .  & .  & .  &
   .  & .  \\
 .  & .  & .  & .  & .  & .  & .  & .  & .  & .  & .  & .  & .  & .  & .  & .  & 1 & .  & 1 & .  & .  & .  &
   .  & .  \\
 .  & .  & .  & .  & .  & .  & .  & .  & .  & .  & .  & .  & .  & .  & .  & .  & .  & 1 & 1 & 1 & .  & 1 &
   .  & .  \\
 .  & .  & .  & .  & .  & .  & .  & .  & .  & .  & .  & .  & .  & .  & .  & .  & .  & .  & 1 & .  & 1 & .  &
   .  & .  \\
 .  & .  & .  & .  & .  & .  & .  & .  & .  & .  & .  & .  & .  & .  & .  & .  & .  & .  & .  & 1 & .  & 1 &
   .  & .  \\
 .  & .  & .  & .  & .  & .  & .  & .  & .  & .  & .  & .  & .  & .  & .  & .  & .  & .  & 1 & .  & 1 & 1 &
   1 & .  \\
 .  & .  & .  & .  & .  & .  & .  & .  & .  & .  & .  & .  & .  & .  & .  & .  & .  & .  & .  & .  & .  & 1 &
   .  & 1 \\
 .  & .  & .  & .  & .  & .  & .  & .  & .  & .  & .  & .  & .  & .  & .  & .  & .  & .  & .  & .  & .  & .  &
   1 & 1 \\
 1 & 1 & .  & .  & .  & .  & .  & .  & .  & .  & .  & .  & .  & .  & .  & .  & .  & .  & .  & .  & .  & .  &
   .  & .  \\
 .  & 1 & 1 & .  & .  & .  & .  & .  & .  & .  & .  & .  & .  & .  & .  & .  & .  & .  & .  & .  & .  & .  &
   .  & .  \\
 .  & 1 & .  & 1 & 1 & 1 & .  & .  & .  & .  & .  & .  & .  & .  & .  & .  & .  & .  & .  & .  & .  & .  &
   .  & .  \\
 .  & .  & 1 & .  & 1 & .  & .  & .  & .  & .  & .  & .  & .  & .  & .  & .  & .  & .  & .  & .  & .  & .  &
   .  & .  \\
 .  & .  & .  & 1 & .  & 1 & .  & .  & .  & .  & .  & .  & .  & .  & .  & .  & .  & .  & .  & .  & .  & .  &
   .  & .  \\
 .  & .  & 1 & .  & 1 & 1 & 1 & .  & .  & .  & .  & .  & .  & .  & .  & .  & .  & .  & .  & .  & .  & .  &
   .  & .  \\
 .  & .  & .  & 1 & .  & .  & 1 & .  & .  & .  & .  & .  & .  & .  & .  & .  & .  & .  & .  & .  & .  & .  &
   .  & .  \\
 .  & .  & .  & .  & .  & .  & 1 & 1 & .  & .  & .  & .  & .  & .  & .  & .  & .  & .  & .  & .  & .  & .  &
   .  & .  \\
\end{array}
\right)$$
}

\subsection{Gram matrices for lattices considered in the text}
\label{AppendixGram}
Gram matrices for lattices $L_1$ and $L_2$, associated with the ${\mathcal A}_k$ series, were given before. 
Here we display $L_3$ (it is small enough to fit in a single page), and we also give those associated with
${\mathcal D}_3$, ${\mathcal D}_6$, and with the three exceptional  ${\mathcal E}_5$, ${\mathcal E}_9$, ${\mathcal E}_{21}$. 

\smallskip

Lattice $L_3$ associated with ${\mathcal A}_3$.
{\scriptsize
$$ A=
\left(
\begin{array}{cccccccccccccccccccc}
 6 & 0 & 0 & 0 & 2 & 0 & 0 & 2 & 0 & 0 & -2 & 1 & 0 & 0 & -2 & 2 & 0 & -2 & 2 & 0 \\
 0 & 6 & 0 & 0 & 2 & 2 & 2 & 2 & 2 & 2 & 1 & 0 & 1 & 1 & 0 & 0 & 0 & 0 & 0 & 0 \\
 0 & 0 & 6 & 0 & 0 & 0 & 2 & 0 & 2 & 0 & 0 & 1 & -2 & 0 & 2 & 0 & -2 & 0 & -2 & 2 \\
 0 & 0 & 0 & 6 & 0 & 2 & 0 & 0 & 0 & 2 & 0 & 1 & 0 & -2 & 0 & -2 & 2 & 2 & 0 & -2 \\
 2 & 2 & 0 & 0 & 6 & 0 & 0 & 2 & 2 & 0 & 2 & 2 & 0 & 0 & -1 & 1 & 1 & 2 & 2 & 0 \\
 0 & 2 & 0 & 2 & 0 & 6 & 0 & 2 & 0 & 2 & 0 & 2 & 0 & 2 & 1 & -1 & 1 & 2 & 0 & 2 \\
 0 & 2 & 2 & 0 & 0 & 0 & 6 & 0 & 2 & 2 & 0 & 2 & 2 & 0 & 1 & 1 & -1 & 0 & 2 & 2 \\
 2 & 2 & 0 & 0 & 2 & 2 & 0 & 6 & 0 & 0 & 2 & 2 & 0 & 0 & 2 & 2 & 0 & -1 & 1 & 1 \\
 0 & 2 & 2 & 0 & 2 & 0 & 2 & 0 & 6 & 0 & 0 & 2 & 2 & 0 & 2 & 0 & 2 & 1 & -1 & 1 \\
 0 & 2 & 0 & 2 & 0 & 2 & 2 & 0 & 0 & 6 & 0 & 2 & 0 & 2 & 0 & 2 & 2 & 1 & 1 & -1 \\
 -2 & 1 & 0 & 0 & 2 & 0 & 0 & 2 & 0 & 0 & 6 & 0 & 0 & 0 & 2 & 0 & 0 & 2 & 0 & 0 \\
 1 & 0 & 1 & 1 & 2 & 2 & 2 & 2 & 2 & 2 & 0 & 6 & 0 & 0 & 2 & 2 & 2 & 2 & 2 & 2 \\
 0 & 1 & -2 & 0 & 0 & 0 & 2 & 0 & 2 & 0 & 0 & 0 & 6 & 0 & 0 & 0 & 2 & 0 & 2 & 0 \\
 0 & 1 & 0 & -2 & 0 & 2 & 0 & 0 & 0 & 2 & 0 & 0 & 0 & 6 & 0 & 2 & 0 & 0 & 0 & 2 \\
 -2 & 0 & 2 & 0 & -1 & 1 & 1 & 2 & 2 & 0 & 2 & 2 & 0 & 0 & 6 & 0 & 0 & 0 & -2 & 2 \\
 2 & 0 & 0 & -2 & 1 & -1 & 1 & 2 & 0 & 2 & 0 & 2 & 0 & 2 & 0 & 6 & 0 & -2 & 2 & 0 \\
 0 & 0 & -2 & 2 & 1 & 1 & -1 & 0 & 2 & 2 & 0 & 2 & 2 & 0 & 0 & 0 & 6 & 2 & 0 & -2 \\
 -2 & 0 & 0 & 2 & 2 & 2 & 0 & -1 & 1 & 1 & 2 & 2 & 0 & 0 & 0 & -2 & 2 & 6 & 0 & 0 \\
 2 & 0 & -2 & 0 & 2 & 0 & 2 & 1 & -1 & 1 & 0 & 2 & 2 & 0 & -2 & 2 & 0 & 0 & 6 & 0 \\
 0 & 0 & 2 & -2 & 0 & 2 & 2 & 1 & 1 & -1 & 0 & 2 & 0 & 2 & 2 & 0 & -2 & 0 & 0 & 6 \\
\end{array}
\right)
\label{GramA3}
$$
}

Lattice ${\mathcal D}_3$
{\scriptsize
$$ A=
\left(
\begin{array}{cccccccccccc}
 6 & 0 & 0 & 0 & 2 & 2 & -2 & 1 & 1 & 1 & 0 & 0 \\
 0 & 6 & 0 & 0 & 2 & 2 & 1 & -2 & 1 & 1 & 0 & 0 \\
 0 & 0 & 6 & 0 & 2 & 2 & 1 & 1 & -2 & 1 & 0 & 0 \\
 0 & 0 & 0 & 6 & 2 & 2 & 1 & 1 & 1 & -2 & 0 & 0 \\
 2 & 2 & 2 & 2 & 6 & 4 & 2 & 2 & 2 & 2 & 1 & 4 \\
 2 & 2 & 2 & 2 & 4 & 6 & 2 & 2 & 2 & 2 & 4 & 1 \\
 -2 & 1 & 1 & 1 & 2 & 2 & 6 & 0 & 0 & 0 & 2 & 2 \\
 1 & -2 & 1 & 1 & 2 & 2 & 0 & 6 & 0 & 0 & 2 & 2 \\
 1 & 1 & -2 & 1 & 2 & 2 & 0 & 0 & 6 & 0 & 2 & 2 \\
 1 & 1 & 1 & -2 & 2 & 2 & 0 & 0 & 0 & 6 & 2 & 2 \\
 0 & 0 & 0 & 0 & 1 & 4 & 2 & 2 & 2 & 2 & 6 & 0 \\
 0 & 0 & 0 & 0 & 4 & 1 & 2 & 2 & 2 & 2 & 0 & 6 \\
\end{array}
\right)
$$
}

Lattice ${\mathcal D}_6$
{\scriptsize
$$ A=
\left(
\begin{array}{cccccccccccccccccccccccc}
 6 & 0 & 0 & 0 & 0 & 0 & 0 & 0 & 2 & 0 & 2 & 0 & -2 & 1 & 0 & 1 & 1 & 1 & 2 & 0 & 0 & 0 & 0 & 2 \\
 0 & 6 & 0 & 0 & 0 & 0 & 0 & 0 & 2 & 0 & 2 & 0 & 1 & -2 & 0 & 1 & 1 & 1 & 2 & 0 & 0 & 0 & 0 & 2 \\
 0 & 0 & 6 & 0 & 0 & 0 & 0 & 2 & 0 & 2 & 0 & 0 & 0 & 0 & -2 & 1 & 0 & 0 & 2 & -2 & 0 & -2 & 0 & 2 \\
 0 & 0 & 0 & 6 & 0 & 0 & 2 & 2 & 2 & 2 & 2 & 2 & 1 & 1 & 1 & 0 & 1 & 2 & 0 & 0 & 2 & 0 & 2 & 0 \\
 0 & 0 & 0 & 0 & 6 & 0 & 0 & 0 & 2 & 0 & 2 & 0 & 1 & 1 & 0 & 1 & -2 & 1 & 2 & 0 & 0 & 0 & 0 & 2 \\
 0 & 0 & 0 & 0 & 0 & 6 & 2 & 0 & 2 & 0 & 2 & 2 & 1 & 1 & 0 & 2 & 1 & -1 & -2 & 2 & 2 & 2 & 2 & -2 \\
 0 & 0 & 0 & 2 & 0 & 2 & 6 & 0 & 0 & 2 & 2 & 0 & 0 & 0 & 0 & 2 & 0 & 2 & -1 & 1 & 2 & 2 & 2 & 0 \\
 0 & 0 & 2 & 2 & 0 & 0 & 0 & 6 & 0 & 2 & 0 & 2 & 0 & 0 & 2 & 2 & 0 & 0 & 1 & -1 & 1 & 2 & 0 & 2 \\
 2 & 2 & 0 & 2 & 2 & 2 & 0 & 0 & 6 & 0 & 4 & 2 & 2 & 2 & 0 & 2 & 2 & 2 & 2 & 1 & 2 & 0 & 4 & 2 \\
 0 & 0 & 2 & 2 & 0 & 0 & 2 & 2 & 0 & 6 & 0 & 0 & 0 & 0 & 2 & 2 & 0 & 0 & 2 & 2 & 0 & -1 & 1 & 1 \\
 2 & 2 & 0 & 2 & 2 & 2 & 2 & 0 & 4 & 0 & 6 & 0 & 2 & 2 & 0 & 2 & 2 & 2 & 2 & 0 & 4 & 1 & 2 & 2 \\
 0 & 0 & 0 & 2 & 0 & 2 & 0 & 2 & 2 & 0 & 0 & 6 & 0 & 0 & 0 & 2 & 0 & 2 & 0 & 2 & 2 & 1 & 2 & -1 \\
 -2 & 1 & 0 & 1 & 1 & 1 & 0 & 0 & 2 & 0 & 2 & 0 & 6 & 0 & 0 & 0 & 0 & 0 & 0 & 0 & 2 & 0 & 2 & 0 \\
 1 & -2 & 0 & 1 & 1 & 1 & 0 & 0 & 2 & 0 & 2 & 0 & 0 & 6 & 0 & 0 & 0 & 0 & 0 & 0 & 2 & 0 & 2 & 0 \\
 0 & 0 & -2 & 1 & 0 & 0 & 0 & 2 & 0 & 2 & 0 & 0 & 0 & 0 & 6 & 0 & 0 & 0 & 0 & 2 & 0 & 2 & 0 & 0 \\
 1 & 1 & 1 & 0 & 1 & 2 & 2 & 2 & 2 & 2 & 2 & 2 & 0 & 0 & 0 & 6 & 0 & 0 & 2 & 2 & 2 & 2 & 2 & 2 \\
 1 & 1 & 0 & 1 & -2 & 1 & 0 & 0 & 2 & 0 & 2 & 0 & 0 & 0 & 0 & 0 & 6 & 0 & 0 & 0 & 2 & 0 & 2 & 0 \\
 1 & 1 & 0 & 2 & 1 & -1 & 2 & 0 & 2 & 0 & 2 & 2 & 0 & 0 & 0 & 0 & 0 & 6 & 2 & 0 & 2 & 0 & 2 & 2 \\
 2 & 2 & 2 & 0 & 2 & -2 & -1 & 1 & 2 & 2 & 2 & 0 & 0 & 0 & 0 & 2 & 0 & 2 & 6 & 0 & 0 & -2 & 0 & 4 \\
 0 & 0 & -2 & 0 & 0 & 2 & 1 & -1 & 1 & 2 & 0 & 2 & 0 & 0 & 2 & 2 & 0 & 0 & 0 & 6 & 0 & 0 & 2 & -2 \\
 0 & 0 & 0 & 2 & 0 & 2 & 2 & 1 & 2 & 0 & 4 & 2 & 2 & 2 & 0 & 2 & 2 & 2 & 0 & 0 & 6 & 2 & 2 & 0 \\
 0 & 0 & -2 & 0 & 0 & 2 & 2 & 2 & 0 & -1 & 1 & 1 & 0 & 0 & 2 & 2 & 0 & 0 & -2 & 0 & 2 & 6 & 0 & 0 \\
 0 & 0 & 0 & 2 & 0 & 2 & 2 & 0 & 4 & 1 & 2 & 2 & 2 & 2 & 0 & 2 & 2 & 2 & 0 & 2 & 2 & 0 & 6 & 0 \\
 2 & 2 & 2 & 0 & 2 & -2 & 0 & 2 & 2 & 1 & 2 & -1 & 0 & 0 & 0 & 2 & 0 & 2 & 4 & -2 & 0 & 0 & 0 & 6 \\
\end{array}
\right)
$$
}

Lattice ${\mathcal E}_5$
{\scriptsize
$$ A=
\left(
\begin{array}{cccccccccccccccccccccccc}
 6 & 0 & 0 & 0 & 0 & 0 & 2 & 0 & 0 & 0 & 2 & 0 & -2 & 0 & 1 & 1 & 0 & 2 & -2 & 2 & 2 & 0 & -2 & 2 \\
 0 & 6 & 0 & 0 & 0 & 0 & 0 & 2 & 0 & 0 & 0 & 2 & 0 & -2 & 1 & 1 & 2 & 0 & 2 & -2 & 0 & 2 & 2 & -2 \\
 0 & 0 & 6 & 0 & 2 & 0 & 2 & 2 & 2 & 0 & 2 & 2 & 1 & 1 & 0 & 2 & -2 & 2 & 2 & 0 & -2 & 2 & 0 & 2 \\
 0 & 0 & 0 & 6 & 0 & 2 & 2 & 2 & 0 & 2 & 2 & 2 & 1 & 1 & 2 & 0 & 2 & -2 & 0 & 2 & 2 & -2 & 2 & 0 \\
 0 & 0 & 2 & 0 & 6 & 0 & 0 & 0 & 0 & 0 & 0 & 2 & 0 & 0 & 2 & 0 & -2 & 0 & 1 & 1 & 0 & 0 & 0 & 2 \\
 0 & 0 & 0 & 2 & 0 & 6 & 0 & 0 & 0 & 0 & 2 & 0 & 0 & 0 & 0 & 2 & 0 & -2 & 1 & 1 & 0 & 0 & 2 & 0 \\
 2 & 0 & 2 & 2 & 0 & 0 & 6 & 0 & 2 & 0 & 2 & 2 & 2 & 0 & 2 & 2 & 1 & 1 & 0 & 2 & 2 & 0 & 2 & 2 \\
 0 & 2 & 2 & 2 & 0 & 0 & 0 & 6 & 0 & 2 & 2 & 2 & 0 & 2 & 2 & 2 & 1 & 1 & 2 & 0 & 0 & 2 & 2 & 2 \\
 0 & 0 & 2 & 0 & 0 & 0 & 2 & 0 & 6 & 0 & 0 & 0 & 0 & 0 & 2 & 0 & 0 & 0 & 2 & 0 & -2 & 0 & 1 & 1 \\
 0 & 0 & 0 & 2 & 0 & 0 & 0 & 2 & 0 & 6 & 0 & 0 & 0 & 0 & 0 & 2 & 0 & 0 & 0 & 2 & 0 & -2 & 1 & 1 \\
 2 & 0 & 2 & 2 & 0 & 2 & 2 & 2 & 0 & 0 & 6 & 0 & 2 & 0 & 2 & 2 & 0 & 2 & 2 & 2 & 1 & 1 & 0 & 2 \\
 0 & 2 & 2 & 2 & 2 & 0 & 2 & 2 & 0 & 0 & 0 & 6 & 0 & 2 & 2 & 2 & 2 & 0 & 2 & 2 & 1 & 1 & 2 & 0 \\
 -2 & 0 & 1 & 1 & 0 & 0 & 2 & 0 & 0 & 0 & 2 & 0 & 6 & 0 & 0 & 0 & 0 & 0 & 2 & 0 & 0 & 0 & 2 & 0 \\
 0 & -2 & 1 & 1 & 0 & 0 & 0 & 2 & 0 & 0 & 0 & 2 & 0 & 6 & 0 & 0 & 0 & 0 & 0 & 2 & 0 & 0 & 0 & 2 \\
 1 & 1 & 0 & 2 & 2 & 0 & 2 & 2 & 2 & 0 & 2 & 2 & 0 & 0 & 6 & 0 & 2 & 0 & 2 & 2 & 2 & 0 & 2 & 2 \\
 1 & 1 & 2 & 0 & 0 & 2 & 2 & 2 & 0 & 2 & 2 & 2 & 0 & 0 & 0 & 6 & 0 & 2 & 2 & 2 & 0 & 2 & 2 & 2 \\
 0 & 2 & -2 & 2 & -2 & 0 & 1 & 1 & 0 & 0 & 0 & 2 & 0 & 0 & 2 & 0 & 6 & 0 & 0 & 0 & 2 & 0 & 2 & -2 \\
 2 & 0 & 2 & -2 & 0 & -2 & 1 & 1 & 0 & 0 & 2 & 0 & 0 & 0 & 0 & 2 & 0 & 6 & 0 & 0 & 0 & 2 & -2 & 2 \\
 -2 & 2 & 2 & 0 & 1 & 1 & 0 & 2 & 2 & 0 & 2 & 2 & 2 & 0 & 2 & 2 & 0 & 0 & 6 & 0 & -2 & 2 & 2 & 0 \\
 2 & -2 & 0 & 2 & 1 & 1 & 2 & 0 & 0 & 2 & 2 & 2 & 0 & 2 & 2 & 2 & 0 & 0 & 0 & 6 & 2 & -2 & 0 & 2 \\
 2 & 0 & -2 & 2 & 0 & 0 & 2 & 0 & -2 & 0 & 1 & 1 & 0 & 0 & 2 & 0 & 2 & 0 & -2 & 2 & 6 & 0 & 0 & 0 \\
 0 & 2 & 2 & -2 & 0 & 0 & 0 & 2 & 0 & -2 & 1 & 1 & 0 & 0 & 0 & 2 & 0 & 2 & 2 & -2 & 0 & 6 & 0 & 0 \\
 -2 & 2 & 0 & 2 & 0 & 2 & 2 & 2 & 1 & 1 & 0 & 2 & 2 & 0 & 2 & 2 & 2 & -2 & 2 & 0 & 0 & 0 & 6 & 0 \\
 2 & -2 & 2 & 0 & 2 & 0 & 2 & 2 & 1 & 1 & 2 & 0 & 0 & 2 & 2 & 2 & -2 & 2 & 0 & 2 & 0 & 0 & 0 & 6 \\
\end{array}
\right)
$$
}

Lattice ${\mathcal E}_9$
{\scriptsize
$$ A=
\left(
\begin{array}{cccccccccccccccccccccccc}
 6 & 0 & 0 & 0 & 2 & 0 & 0 & 0 & 2 & 0 & 0 & 0 & -2 & 0 & 0 & 1 & -2 & 0 & 0 & 2 & -2 & 0 & 0 & 2 \\
 0 & 6 & 0 & 0 & 0 & 2 & 0 & 0 & 0 & 2 & 0 & 0 & 0 & -2 & 0 & 1 & 0 & -2 & 0 & 2 & 0 & -2 & 0 & 2 \\
 0 & 0 & 6 & 0 & 0 & 0 & 2 & 0 & 0 & 0 & 2 & 0 & 0 & 0 & -2 & 1 & 0 & 0 & -2 & 2 & 0 & 0 & -2 & 2 \\
 0 & 0 & 0 & 6 & 2 & 2 & 2 & 4 & 2 & 2 & 2 & 4 & 1 & 1 & 1 & 4 & 2 & 2 & 2 & 2 & 2 & 2 & 2 & 2 \\
 2 & 0 & 0 & 2 & 6 & 0 & 0 & 0 & 2 & 0 & 0 & 2 & 2 & 0 & 0 & 2 & -1 & 1 & 1 & 2 & 2 & 0 & 0 & 2 \\
 0 & 2 & 0 & 2 & 0 & 6 & 0 & 0 & 0 & 2 & 0 & 2 & 0 & 2 & 0 & 2 & 1 & -1 & 1 & 2 & 0 & 2 & 0 & 2 \\
 0 & 0 & 2 & 2 & 0 & 0 & 6 & 0 & 0 & 0 & 2 & 2 & 0 & 0 & 2 & 2 & 1 & 1 & -1 & 2 & 0 & 0 & 2 & 2 \\
 0 & 0 & 0 & 4 & 0 & 0 & 0 & 6 & 2 & 2 & 2 & 2 & 0 & 0 & 0 & 4 & 2 & 2 & 2 & 1 & 2 & 2 & 2 & 2 \\
 2 & 0 & 0 & 2 & 2 & 0 & 0 & 2 & 6 & 0 & 0 & 0 & 2 & 0 & 0 & 2 & 2 & 0 & 0 & 2 & -1 & 1 & 1 & 2 \\
 0 & 2 & 0 & 2 & 0 & 2 & 0 & 2 & 0 & 6 & 0 & 0 & 0 & 2 & 0 & 2 & 0 & 2 & 0 & 2 & 1 & -1 & 1 & 2 \\
 0 & 0 & 2 & 2 & 0 & 0 & 2 & 2 & 0 & 0 & 6 & 0 & 0 & 0 & 2 & 2 & 0 & 0 & 2 & 2 & 1 & 1 & -1 & 2 \\
 0 & 0 & 0 & 4 & 2 & 2 & 2 & 2 & 0 & 0 & 0 & 6 & 0 & 0 & 0 & 4 & 2 & 2 & 2 & 2 & 2 & 2 & 2 & 1 \\
 -2 & 0 & 0 & 1 & 2 & 0 & 0 & 0 & 2 & 0 & 0 & 0 & 6 & 0 & 0 & 0 & 2 & 0 & 0 & 0 & 2 & 0 & 0 & 0 \\
 0 & -2 & 0 & 1 & 0 & 2 & 0 & 0 & 0 & 2 & 0 & 0 & 0 & 6 & 0 & 0 & 0 & 2 & 0 & 0 & 0 & 2 & 0 & 0 \\
 0 & 0 & -2 & 1 & 0 & 0 & 2 & 0 & 0 & 0 & 2 & 0 & 0 & 0 & 6 & 0 & 0 & 0 & 2 & 0 & 0 & 0 & 2 & 0 \\
 1 & 1 & 1 & 4 & 2 & 2 & 2 & 4 & 2 & 2 & 2 & 4 & 0 & 0 & 0 & 6 & 2 & 2 & 2 & 4 & 2 & 2 & 2 & 4 \\
 -2 & 0 & 0 & 2 & -1 & 1 & 1 & 2 & 2 & 0 & 0 & 2 & 2 & 0 & 0 & 2 & 6 & 0 & 0 & 0 & 0 & 2 & 2 & 0 \\
 0 & -2 & 0 & 2 & 1 & -1 & 1 & 2 & 0 & 2 & 0 & 2 & 0 & 2 & 0 & 2 & 0 & 6 & 0 & 0 & 2 & 0 & 2 & 0 \\
 0 & 0 & -2 & 2 & 1 & 1 & -1 & 2 & 0 & 0 & 2 & 2 & 0 & 0 & 2 & 2 & 0 & 0 & 6 & 0 & 2 & 2 & 0 & 0 \\
 2 & 2 & 2 & 2 & 2 & 2 & 2 & 1 & 2 & 2 & 2 & 2 & 0 & 0 & 0 & 4 & 0 & 0 & 0 & 6 & 0 & 0 & 0 & 4 \\
 -2 & 0 & 0 & 2 & 2 & 0 & 0 & 2 & -1 & 1 & 1 & 2 & 2 & 0 & 0 & 2 & 0 & 2 & 2 & 0 & 6 & 0 & 0 & 0 \\
 0 & -2 & 0 & 2 & 0 & 2 & 0 & 2 & 1 & -1 & 1 & 2 & 0 & 2 & 0 & 2 & 2 & 0 & 2 & 0 & 0 & 6 & 0 & 0 \\
 0 & 0 & -2 & 2 & 0 & 0 & 2 & 2 & 1 & 1 & -1 & 2 & 0 & 0 & 2 & 2 & 2 & 2 & 0 & 0 & 0 & 0 & 6 & 0 \\
 2 & 2 & 2 & 2 & 2 & 2 & 2 & 2 & 2 & 2 & 2 & 1 & 0 & 0 & 0 & 4 & 0 & 0 & 0 & 4 & 0 & 0 & 0 & 6 \\
\end{array}
\right)
$$
}

\begin{landscape}
Lattice ${\mathcal E}_{21}$.  $A=$
{
\tiny
\arraycolsep=3pt
$$
\left(
\begin{array}{@{}l*{48}{c}@{}}
 6 & 0 & 0 & 0 & 0 & 0 & 0 & 0 & 2 & 0 & 0 & 0 & 0 & 0 & 0 & 0 & 2 & 0 & 0 & 0 & 0 & 0 & 0 & 0 & -2 & 1 & 0 & 0 & 0 & 0 & 0 & 0 & -2 & 2 & 0 & 0 & 0 & 0 & 0 & 0 & -2 & 2
   & 0 & 0 & 0 & 0 & 0 & 0 \\
 0 & 6 & 0 & 0 & 0 & 0 & 0 & 0 & 2 & 2 & 2 & 0 & 0 & 0 & 0 & 0 & 2 & 2 & 2 & 0 & 0 & 0 & 0 & 0 & 1 & 0 & 1 & 1 & 1 & 1 & 0 & 0 & 0 & 0 & 0 & 2 & 0 & 2 & 0 & 0 & 0 & 0 &
   0 & 2 & 0 & 2 & 0 & 0 \\
 0 & 0 & 6 & 0 & 0 & 0 & 0 & 0 & 0 & 0 & 2 & 0 & 2 & 0 & 2 & 0 & 0 & 2 & 0 & 2 & 0 & 2 & 0 & 0 & 0 & 1 & 0 & 0 & 2 & 1 & 1 & 0 & 2 & 0 & 2 & 0 & 0 & 2 & -2 & 0 & 0 & -2
   & 2 & 0 & 0 & 2 & 0 & 2 \\
 0 & 0 & 0 & 6 & 0 & 0 & 0 & 0 & 0 & 2 & 0 & 2 & 0 & 2 & 0 & 0 & 0 & 0 & 2 & 0 & 2 & 0 & 2 & 0 & 0 & 1 & 0 & 0 & 1 & 2 & 1 & 0 & 0 & -2 & 2 & 0 & 0 & 2 & 0 & 2 & 2 & 0 &
   2 & 0 & 0 & 2 & -2 & 0 \\
 0 & 0 & 0 & 0 & 6 & 0 & 0 & 0 & 0 & 0 & 2 & 0 & 2 & 2 & 0 & 0 & 0 & 0 & 2 & 2 & 0 & 2 & 0 & 0 & 0 & 1 & 2 & 1 & 0 & 2 & 1 & 0 & 0 & 2 & 2 & 2 & 0 & 0 & 2 & 0 & 0 & 2 &
   0 & 0 & 2 & 2 & 2 & 0 \\
 0 & 0 & 0 & 0 & 0 & 6 & 0 & 0 & 0 & 0 & 2 & 2 & 0 & 2 & 0 & 0 & 0 & 0 & 2 & 0 & 2 & 2 & 0 & 0 & 0 & 1 & 1 & 2 & 2 & 0 & 1 & 0 & 0 & 2 & 0 & 0 & 2 & 2 & 2 & 0 & 0 & 2 &
   2 & 2 & 0 & 0 & 2 & 0 \\
 0 & 0 & 0 & 0 & 0 & 0 & 6 & 0 & 0 & 0 & 0 & 0 & 0 & 2 & 2 & 2 & 0 & 0 & 0 & 0 & 0 & 2 & 2 & 2 & 0 & 0 & 1 & 1 & 1 & 1 & 0 & 1 & 0 & 0 & 2 & 0 & 2 & 0 & 0 & 0 & 0 & 0 &
   2 & 0 & 2 & 0 & 0 & 0 \\
 0 & 0 & 0 & 0 & 0 & 0 & 0 & 6 & 0 & 0 & 0 & 0 & 0 & 0 & 0 & 2 & 0 & 0 & 0 & 0 & 0 & 0 & 0 & 2 & 0 & 0 & 0 & 0 & 0 & 0 & 1 & -2 & 0 & 0 & 0 & 0 & 0 & 0 & 2 & -2 & 0 & 0
   & 0 & 0 & 0 & 0 & 2 & -2 \\
 2 & 2 & 0 & 0 & 0 & 0 & 0 & 0 & 6 & 0 & 0 & 0 & 0 & 0 & 0 & 0 & 2 & 2 & 0 & 0 & 0 & 0 & 0 & 0 & 2 & 2 & 0 & 0 & 0 & 0 & 0 & 0 & -1 & 1 & 1 & 0 & 0 & 0 & 0 & 0 & 2 & 2 &
   0 & 0 & 0 & 0 & 0 & 0 \\
 0 & 2 & 0 & 2 & 0 & 0 & 0 & 0 & 0 & 6 & 0 & 0 & 0 & 0 & 0 & 0 & 2 & 0 & 2 & 0 & 0 & 0 & 0 & 0 & 0 & 2 & 0 & 2 & 0 & 0 & 0 & 0 & 1 & -1 & 1 & 1 & 0 & 1 & 0 & 0 & 2 & 0 &
   2 & 0 & 0 & 0 & 0 & 0 \\
 0 & 2 & 2 & 0 & 2 & 2 & 0 & 0 & 0 & 0 & 6 & 0 & 0 & 0 & 0 & 0 & 0 & 2 & 2 & 2 & 0 & 2 & 0 & 0 & 0 & 2 & 2 & 0 & 2 & 2 & 0 & 0 & 1 & 1 & 1 & 1 & 2 & 2 & 1 & 0 & 0 & 2 &
   2 & 2 & 0 & 2 & 0 & 0 \\
 0 & 0 & 0 & 2 & 0 & 2 & 0 & 0 & 0 & 0 & 0 & 6 & 0 & 0 & 0 & 0 & 0 & 0 & 2 & 0 & 2 & 0 & 0 & 0 & 0 & 0 & 0 & 2 & 0 & 2 & 0 & 0 & 0 & 1 & 1 & -1 & 0 & 2 & 0 & 0 & 0 & 0 &
   2 & 0 & 2 & 0 & 0 & 0 \\
 0 & 0 & 2 & 0 & 2 & 0 & 0 & 0 & 0 & 0 & 0 & 0 & 6 & 0 & 0 & 0 & 0 & 0 & 0 & 2 & 0 & 2 & 0 & 0 & 0 & 0 & 2 & 0 & 2 & 0 & 0 & 0 & 0 & 0 & 2 & 0 & -1 & 1 & 1 & 0 & 0 & 0 &
   0 & 2 & 0 & 2 & 0 & 0 \\
 0 & 0 & 0 & 2 & 2 & 2 & 2 & 0 & 0 & 0 & 0 & 0 & 0 & 6 & 0 & 0 & 0 & 0 & 2 & 0 & 2 & 2 & 2 & 0 & 0 & 0 & 0 & 2 & 2 & 2 & 2 & 0 & 0 & 1 & 2 & 2 & 1 & 1 & 1 & 1 & 0 & 0 &
   2 & 0 & 2 & 2 & 2 & 0 \\
 0 & 0 & 2 & 0 & 0 & 0 & 2 & 0 & 0 & 0 & 0 & 0 & 0 & 0 & 6 & 0 & 0 & 0 & 0 & 0 & 0 & 2 & 0 & 2 & 0 & 0 & 2 & 0 & 0 & 0 & 2 & 0 & 0 & 0 & 1 & 0 & 1 & 1 & -1 & 1 & 0 & 0 &
   0 & 0 & 0 & 2 & 0 & 2 \\
 0 & 0 & 0 & 0 & 0 & 0 & 2 & 2 & 0 & 0 & 0 & 0 & 0 & 0 & 0 & 6 & 0 & 0 & 0 & 0 & 0 & 0 & 2 & 2 & 0 & 0 & 0 & 0 & 0 & 0 & 2 & 2 & 0 & 0 & 0 & 0 & 0 & 1 & 1 & -1 & 0 & 0 &
   0 & 0 & 0 & 0 & 2 & 2 \\
 2 & 2 & 0 & 0 & 0 & 0 & 0 & 0 & 2 & 2 & 0 & 0 & 0 & 0 & 0 & 0 & 6 & 0 & 0 & 0 & 0 & 0 & 0 & 0 & 2 & 2 & 0 & 0 & 0 & 0 & 0 & 0 & 2 & 2 & 0 & 0 & 0 & 0 & 0 & 0 & -1 & 1 &
   1 & 0 & 0 & 0 & 0 & 0 \\
 0 & 2 & 2 & 0 & 0 & 0 & 0 & 0 & 2 & 0 & 2 & 0 & 0 & 0 & 0 & 0 & 0 & 6 & 0 & 0 & 0 & 0 & 0 & 0 & 0 & 2 & 2 & 0 & 0 & 0 & 0 & 0 & 2 & 0 & 2 & 0 & 0 & 0 & 0 & 0 & 1 & -1 &
   1 & 1 & 0 & 1 & 0 & 0 \\
 0 & 2 & 0 & 2 & 2 & 2 & 0 & 0 & 0 & 2 & 2 & 2 & 0 & 2 & 0 & 0 & 0 & 0 & 6 & 0 & 0 & 0 & 0 & 0 & 0 & 2 & 0 & 2 & 2 & 2 & 0 & 0 & 0 & 2 & 2 & 2 & 0 & 2 & 0 & 0 & 1 & 1 &
   1 & 1 & 2 & 2 & 1 & 0 \\
 0 & 0 & 2 & 0 & 2 & 0 & 0 & 0 & 0 & 0 & 2 & 0 & 2 & 0 & 0 & 0 & 0 & 0 & 0 & 6 & 0 & 0 & 0 & 0 & 0 & 0 & 2 & 0 & 2 & 0 & 0 & 0 & 0 & 0 & 2 & 0 & 2 & 0 & 0 & 0 & 0 & 1 &
   1 & -1 & 0 & 2 & 0 & 0 \\
 0 & 0 & 0 & 2 & 0 & 2 & 0 & 0 & 0 & 0 & 0 & 2 & 0 & 2 & 0 & 0 & 0 & 0 & 0 & 0 & 6 & 0 & 0 & 0 & 0 & 0 & 0 & 2 & 0 & 2 & 0 & 0 & 0 & 0 & 0 & 2 & 0 & 2 & 0 & 0 & 0 & 0 &
   2 & 0 & -1 & 1 & 1 & 0 \\
 0 & 0 & 2 & 0 & 2 & 2 & 2 & 0 & 0 & 0 & 2 & 0 & 2 & 2 & 2 & 0 & 0 & 0 & 0 & 0 & 0 & 6 & 0 & 0 & 0 & 0 & 2 & 0 & 2 & 2 & 2 & 0 & 0 & 0 & 2 & 0 & 2 & 2 & 2 & 0 & 0 & 1 &
   2 & 2 & 1 & 1 & 1 & 1 \\
 0 & 0 & 0 & 2 & 0 & 0 & 2 & 0 & 0 & 0 & 0 & 0 & 0 & 2 & 0 & 2 & 0 & 0 & 0 & 0 & 0 & 0 & 6 & 0 & 0 & 0 & 0 & 2 & 0 & 0 & 2 & 0 & 0 & 0 & 0 & 0 & 0 & 2 & 0 & 2 & 0 & 0 &
   1 & 0 & 1 & 1 & -1 & 1 \\
 0 & 0 & 0 & 0 & 0 & 0 & 2 & 2 & 0 & 0 & 0 & 0 & 0 & 0 & 2 & 2 & 0 & 0 & 0 & 0 & 0 & 0 & 0 & 6 & 0 & 0 & 0 & 0 & 0 & 0 & 2 & 2 & 0 & 0 & 0 & 0 & 0 & 0 & 2 & 2 & 0 & 0 &
   0 & 0 & 0 & 1 & 1 & -1 \\
 -2 & 1 & 0 & 0 & 0 & 0 & 0 & 0 & 2 & 0 & 0 & 0 & 0 & 0 & 0 & 0 & 2 & 0 & 0 & 0 & 0 & 0 & 0 & 0 & 6 & 0 & 0 & 0 & 0 & 0 & 0 & 0 & 2 & 0 & 0 & 0 & 0 & 0 & 0 & 0 & 2 & 0 &
   0 & 0 & 0 & 0 & 0 & 0 \\
 1 & 0 & 1 & 1 & 1 & 1 & 0 & 0 & 2 & 2 & 2 & 0 & 0 & 0 & 0 & 0 & 2 & 2 & 2 & 0 & 0 & 0 & 0 & 0 & 0 & 6 & 0 & 0 & 0 & 0 & 0 & 0 & 2 & 2 & 2 & 0 & 0 & 0 & 0 & 0 & 2 & 2 &
   2 & 0 & 0 & 0 & 0 & 0 \\
 0 & 1 & 0 & 0 & 2 & 1 & 1 & 0 & 0 & 0 & 2 & 0 & 2 & 0 & 2 & 0 & 0 & 2 & 0 & 2 & 0 & 2 & 0 & 0 & 0 & 0 & 6 & 0 & 0 & 0 & 0 & 0 & 0 & 0 & 2 & 0 & 2 & 0 & 2 & 0 & 0 & 2 &
   0 & 2 & 0 & 2 & 0 & 0 \\
 0 & 1 & 0 & 0 & 1 & 2 & 1 & 0 & 0 & 2 & 0 & 2 & 0 & 2 & 0 & 0 & 0 & 0 & 2 & 0 & 2 & 0 & 2 & 0 & 0 & 0 & 0 & 6 & 0 & 0 & 0 & 0 & 0 & 2 & 0 & 2 & 0 & 2 & 0 & 0 & 0 & 0 &
   2 & 0 & 2 & 0 & 2 & 0 \\
 0 & 1 & 2 & 1 & 0 & 2 & 1 & 0 & 0 & 0 & 2 & 0 & 2 & 2 & 0 & 0 & 0 & 0 & 2 & 2 & 0 & 2 & 0 & 0 & 0 & 0 & 0 & 0 & 6 & 0 & 0 & 0 & 0 & 0 & 2 & 0 & 2 & 2 & 0 & 0 & 0 & 0 &
   2 & 2 & 0 & 2 & 0 & 0 \\
 0 & 1 & 1 & 2 & 2 & 0 & 1 & 0 & 0 & 0 & 2 & 2 & 0 & 2 & 0 & 0 & 0 & 0 & 2 & 0 & 2 & 2 & 0 & 0 & 0 & 0 & 0 & 0 & 0 & 6 & 0 & 0 & 0 & 0 & 2 & 2 & 0 & 2 & 0 & 0 & 0 & 0 &
   2 & 0 & 2 & 2 & 0 & 0 \\
 0 & 0 & 1 & 1 & 1 & 1 & 0 & 1 & 0 & 0 & 0 & 0 & 0 & 2 & 2 & 2 & 0 & 0 & 0 & 0 & 0 & 2 & 2 & 2 & 0 & 0 & 0 & 0 & 0 & 0 & 6 & 0 & 0 & 0 & 0 & 0 & 0 & 2 & 2 & 2 & 0 & 0 &
   0 & 0 & 0 & 2 & 2 & 2 \\
 0 & 0 & 0 & 0 & 0 & 0 & 1 & -2 & 0 & 0 & 0 & 0 & 0 & 0 & 0 & 2 & 0 & 0 & 0 & 0 & 0 & 0 & 0 & 2 & 0 & 0 & 0 & 0 & 0 & 0 & 0 & 6 & 0 & 0 & 0 & 0 & 0 & 0 & 0 & 2 & 0 & 0 &
   0 & 0 & 0 & 0 & 0 & 2 \\
 -2 & 0 & 2 & 0 & 0 & 0 & 0 & 0 & -1 & 1 & 1 & 0 & 0 & 0 & 0 & 0 & 2 & 2 & 0 & 0 & 0 & 0 & 0 & 0 & 2 & 2 & 0 & 0 & 0 & 0 & 0 & 0 & 6 & 0 & 0 & 0 & 0 & 0 & 0 & 0 & 0 & -2
   & 2 & 0 & 0 & 0 & 0 & 0 \\
 2 & 0 & 0 & -2 & 2 & 2 & 0 & 0 & 1 & -1 & 1 & 1 & 0 & 1 & 0 & 0 & 2 & 0 & 2 & 0 & 0 & 0 & 0 & 0 & 0 & 2 & 0 & 2 & 0 & 0 & 0 & 0 & 0 & 6 & 0 & 0 & 0 & 0 & 0 & 0 & -2 & 2
   & 0 & 0 & 2 & 0 & 2 & 0 \\
 0 & 0 & 2 & 2 & 2 & 0 & 2 & 0 & 1 & 1 & 1 & 1 & 2 & 2 & 1 & 0 & 0 & 2 & 2 & 2 & 0 & 2 & 0 & 0 & 0 & 2 & 2 & 0 & 2 & 2 & 0 & 0 & 0 & 0 & 6 & 0 & 0 & 0 & 0 & 0 & 2 & 0 &
   2 & 0 & 2 & 2 & 0 & 0 \\
 0 & 2 & 0 & 0 & 2 & 0 & 0 & 0 & 0 & 1 & 1 & -1 & 0 & 2 & 0 & 0 & 0 & 0 & 2 & 0 & 2 & 0 & 0 & 0 & 0 & 0 & 0 & 2 & 0 & 2 & 0 & 0 & 0 & 0 & 0 & 6 & 0 & 0 & 0 & 0 & 0 & 0 &
   0 & 0 & 0 & 2 & 2 & 0 \\
 0 & 0 & 0 & 0 & 0 & 2 & 2 & 0 & 0 & 0 & 2 & 0 & -1 & 1 & 1 & 0 & 0 & 0 & 0 & 2 & 0 & 2 & 0 & 0 & 0 & 0 & 2 & 0 & 2 & 0 & 0 & 0 & 0 & 0 & 0 & 0 & 6 & 0 & 0 & 0 & 0 & 2 &
   2 & 0 & 0 & 0 & 0 & 0 \\
 0 & 2 & 2 & 2 & 0 & 2 & 0 & 0 & 0 & 1 & 2 & 2 & 1 & 1 & 1 & 1 & 0 & 0 & 2 & 0 & 2 & 2 & 2 & 0 & 0 & 0 & 0 & 2 & 2 & 2 & 2 & 0 & 0 & 0 & 0 & 0 & 0 & 6 & 0 & 0 & 0 & 0 &
   2 & 2 & 0 & 2 & 0 & 2 \\
 0 & 0 & -2 & 0 & 2 & 2 & 0 & 2 & 0 & 0 & 1 & 0 & 1 & 1 & -1 & 1 & 0 & 0 & 0 & 0 & 0 & 2 & 0 & 2 & 0 & 0 & 2 & 0 & 0 & 0 & 2 & 0 & 0 & 0 & 0 & 0 & 0 & 0 & 6 & 0 & 0 & 2
   & 0 & 2 & 0 & 0 & 2 & -2 \\
 0 & 0 & 0 & 2 & 0 & 0 & 0 & -2 & 0 & 0 & 0 & 0 & 0 & 1 & 1 & -1 & 0 & 0 & 0 & 0 & 0 & 0 & 2 & 2 & 0 & 0 & 0 & 0 & 0 & 0 & 2 & 2 & 0 & 0 & 0 & 0 & 0 & 0 & 0 & 6 & 0 & 0
   & 0 & 0 & 0 & 2 & -2 & 0 \\
 -2 & 0 & 0 & 2 & 0 & 0 & 0 & 0 & 2 & 2 & 0 & 0 & 0 & 0 & 0 & 0 & -1 & 1 & 1 & 0 & 0 & 0 & 0 & 0 & 2 & 2 & 0 & 0 & 0 & 0 & 0 & 0 & 0 & -2 & 2 & 0 & 0 & 0 & 0 & 0 & 6 & 0
   & 0 & 0 & 0 & 0 & 0 & 0 \\
 2 & 0 & -2 & 0 & 2 & 2 & 0 & 0 & 2 & 0 & 2 & 0 & 0 & 0 & 0 & 0 & 1 & -1 & 1 & 1 & 0 & 1 & 0 & 0 & 0 & 2 & 2 & 0 & 0 & 0 & 0 & 0 & -2 & 2 & 0 & 0 & 2 & 0 & 2 & 0 & 0 & 6
   & 0 & 0 & 0 & 0 & 0 & 0 \\
 0 & 0 & 2 & 2 & 0 & 2 & 2 & 0 & 0 & 2 & 2 & 2 & 0 & 2 & 0 & 0 & 1 & 1 & 1 & 1 & 2 & 2 & 1 & 0 & 0 & 2 & 0 & 2 & 2 & 2 & 0 & 0 & 2 & 0 & 2 & 0 & 2 & 2 & 0 & 0 & 0 & 0 &
   6 & 0 & 0 & 0 & 0 & 0 \\
 0 & 2 & 0 & 0 & 0 & 2 & 0 & 0 & 0 & 0 & 2 & 0 & 2 & 0 & 0 & 0 & 0 & 1 & 1 & -1 & 0 & 2 & 0 & 0 & 0 & 0 & 2 & 0 & 2 & 0 & 0 & 0 & 0 & 0 & 0 & 0 & 0 & 2 & 2 & 0 & 0 & 0 &
   0 & 6 & 0 & 0 & 0 & 0 \\
 0 & 0 & 0 & 0 & 2 & 0 & 2 & 0 & 0 & 0 & 0 & 2 & 0 & 2 & 0 & 0 & 0 & 0 & 2 & 0 & -1 & 1 & 1 & 0 & 0 & 0 & 0 & 2 & 0 & 2 & 0 & 0 & 0 & 2 & 2 & 0 & 0 & 0 & 0 & 0 & 0 & 0 &
   0 & 0 & 6 & 0 & 0 & 0 \\
 0 & 2 & 2 & 2 & 2 & 0 & 0 & 0 & 0 & 0 & 2 & 0 & 2 & 2 & 2 & 0 & 0 & 1 & 2 & 2 & 1 & 1 & 1 & 1 & 0 & 0 & 2 & 0 & 2 & 2 & 2 & 0 & 0 & 0 & 2 & 2 & 0 & 2 & 0 & 2 & 0 & 0 &
   0 & 0 & 0 & 6 & 0 & 0 \\
 0 & 0 & 0 & -2 & 2 & 2 & 0 & 2 & 0 & 0 & 0 & 0 & 0 & 2 & 0 & 2 & 0 & 0 & 1 & 0 & 1 & 1 & -1 & 1 & 0 & 0 & 0 & 2 & 0 & 0 & 2 & 0 & 0 & 2 & 0 & 2 & 0 & 0 & 2 & -2 & 0 & 0
   & 0 & 0 & 0 & 0 & 6 & 0 \\
 0 & 0 & 2 & 0 & 0 & 0 & 0 & -2 & 0 & 0 & 0 & 0 & 0 & 0 & 2 & 2 & 0 & 0 & 0 & 0 & 0 & 1 & 1 & -1 & 0 & 0 & 0 & 0 & 0 & 0 & 2 & 2 & 0 & 0 & 0 & 0 & 0 & 2 & -2 & 0 & 0 & 0
   & 0 & 0 & 0 & 0 & 0 & 6 \\
\end{array}
\right)
$$
}
\end{landscape}


\begin{thebibliography}{99}
\footnotesize

\bibitem{CIZ}  Cappelli A.,  Itzykson  C.  and   Zuber J. -B.,  The ADE classification of minimal and $A_{1}^{(1)}$ conformal invariant theories, {\it Commun.  Math.  Phys.},  {13},  pp 1--26, (1987).  
\bibitem{ConwaySloane} Conway J. and Sloane N.J.A., Sphere Packings, Lattices and Groups (3rd ed.), Springer, (1999).
\bibitem{RC:periodicquivers} Coquereaux R., {Quantum McKay correspondence and global dimensions for fusion and module-categories associated with Lie groups}, {\it Jour. of Algebra}, {398}, pp 258-283, (2014).
\bibitem{CoquereauxSchieberJMP} Coquereaux R. and Schieber G., {Orders and dimensions for sl2 or sl3 module-categories and boundary conformal field theories on a torus},   {\it J. of Mathematical Physics} {\bf 48} (2007) 043511;
\url{http://arxiv.org/abs/math-ph/0610073}
\bibitem{CoquereauxSU3Maroc}Coquereaux R., Hammaoui D., Schieber G., Tahri  E.H., {Comments about quantum symmetries of SU(3) graphs},  {\it Journal of Geometry and Physics} 57 pp 269-292 (2006).
\bibitem{RCsiteWebFusionGraphs} Coquereaux R., Fusion graphs, \url{http://www.cpt.univ-mrs.fr/~coque/quantumfusion/FusionGraphs.html}
\bibitem{CoquereauxZuberNuclPhys} Coquereaux R and  Zuber J.-B.,  {On some properties of SU(3) Fusion Coefficients.},
Contribution to Mathematical Foundations of Quantum Field Theory, special issue in memory of Raymond Stora, 33 pp., {\it Nucl. Phys. B.}, DOI: 10.1016/j.nuclphysb.2016.05.029, (2016).
\bibitem{YellowBook} Di Francesco  P., Matthieu P. and  Senechal   D., Conformal field theory, Springer, (1997).
\bibitem{DiFrancescoZuber}  Di Francesco P. and   Zuber J.-B., SU(N) lattice integrable models associated with graphs,   {\it Nucl.  Phys.},  {B 338},  pp 602--646, (1990).  
\bibitem{Dorey:CoxeterElement} Dorey P., Partition Functions, Intertwiners and the Coxeter Element. arXiv:hep-th/9205040.  {\it Int. J. Mod. Phys} {A8}, pp 193-208, (1993).
\bibitem{EvansPughSU3} Evans D. E.  and Pugh M., Ocneanu cells and Boltzmann weights for the SU(3) ADE graphs. {\it M\"unster J. of Math.} {2},  pp 95-142 (2009)
\bibitem{Finkelberg} Finkelberg, M., An equivalence of fusion categories, {\it Geom. Funct. Anal.}  6 (1996), 249-267.
\bibitem{Huang} Y.-Z. Huang, Vertex operator algebras, the Verlinde conjecture, and modular tensor categories, {\it Proc. Natl. Acad. Sci. USA}, 102 (2005), 5352Ð5356.
\bibitem{Kac:book}  Kac V., Infinite dimensional Lie algebras, Cambridge University Press, Cambridge (1990).
\bibitem{KazhdanLusztig} Kazhdan D. and  Lusztig G., Tensor structures arising from affine Lie algebras, III,    {\it J. Amer. Math. Soc.},  {7}, pp 335--381, (1994).
\bibitem{KirilovOstrik} Kirillov A. and Ostrik V.,  On q-analog of McKay correspondence and ADE classification of SL2 conformal  field theories, {\it Adv. in Math.},  {171- 2}, pp 183--227, (2002).
\bibitem{Magma} Wieb Bosma, John Cannon, and Catherine Playoust, {\sl The Magma algebra system. I. The user language}, {\it J. Symbolic Comput.}, {\bf 24} (1997), 235--265, \url{http://magma.maths.usyd.edu.au}
\bibitem{Mathematica} Wolfram Research, Inc., Mathematica, Champaign, IL (2010).
\bibitem{Ocneanu:paths} Ocneanu  A., Paths on Coxeter diagrams: from Platonic solids and singularities to minimal models and subfactors,  Notes by Goto S.,  {\it Fields Institute Monographs}, Eds. Rajarama Bhat et al, (1999).  
\bibitem{Ocneanu:Bariloche}  Ocneanu  A., The Classification of  subgroups of quantum SU(N), in  ``Quantum symmetries in theoretical physics and mathematics'', Bariloche 2000, Eds. Coquereaux R., Garc\'{\i}a A. and Trinchero~R., {\it  AMS Contemporary Mathematics}, {294}, pp 133--160, (2000). 
\bibitem{Ocneanu:MSRI} Ocneanu  A., Higher Coxeter systems, \\  http://www.msri.org/publications/ln/msri/2000/subfactors/ocneanu, (2000). 
\bibitem{Ocneanu:posters} Ocneanu  A.,   Poster communications. 
\bibitem{Ocneanu:WIP} Ocneanu  A., work in progress.
\bibitem{Ostrik}   Ostrik  V., Module categories, weak Hopf algebras and modular invariants, {\it Transform.  groups},  {8}, no 2, pp 177--206, (2003).
\bibitem{PleskenPohst} Plesken W. and Pohst M., Constructing integral lattices with prescribed minimum,  {\it Mathematics of Computation}, Vol 45, No 171, pp 209-221, and supplement S5-S16. 
 \bibitem{Steinberg} Steinberg R., {Finite reflection groups},  {\it Trans. Amer. Math. Soc.} {91} pp 493-504, (1959).
\bibitem{Zagier:modularforms} Zagier D.B., Elliptic Modular Forms and Their Applications, in {\it  `The 1-2-3 of Modular forms'},  Lectures at a Summer School in Nordfjordeid, Norway, Springer (2008).
\bibitem{DezaGrishukhin} Deza M. and Grishukhin V., Delaunay Polytopes of Cut Lattices, {\it Linear Algebra and Its Applications}, 226-228:667-685 (1995).
\end{thebibliography}
 \end{document}